\documentclass[10pt]{article}  
\usepackage{amsmath}
\usepackage{amssymb}
\usepackage{graphicx}
\usepackage{amsthm}       \newtheorem{thm}{Theorem}[section]
\newtheorem{prop}[thm]{Proposition}
\newtheorem{lem}[thm]{Lemma}
\newtheorem{cor}[thm]{Corollary}  \theoremstyle{definition}
\newtheorem{df}[thm]{Definition}   \theoremstyle{definition}
\newtheorem{ques}[thm]{Question}
\newtheorem{prob}[thm]{Problem}
\newtheorem{rem}[thm]{Remark}                \theoremstyle{plain}
 \theoremstyle{definition}
\newtheorem{ex}[thm]{Example}   \def\CC{\Bbb{C}}
\def\RR{\Bbb{R}}  
\def\CCI{\hat{\CC}}        \def\NN{\Bbb{N}}

\def\G{\Gamma}

\def\B1{{\rm\kern.32em\vrule    width.12em       height1.4ex
depth-.05ex\kern-.28em 1}}

\begin{document}
\title
{Dynamics of 
postcritically bounded 
polynomial 
semigroups I: connected components of the Julia sets
\footnote{Published in Discrete and Continuous Dynamical Sistems Series A,
Vol. 29, No. 3, 2011, 1205--1244.  
2000 Mathematics Subject Classification. 
37F10, 30D05. Keywords: Complex dynamical systems, rational semigroup, polynomial semigroup, 
random iteration, random complex dynamical systems,  
Julia set, fractal geometry, iterated function systems, surrounding order.}}
\author{Hiroki Sumi\\ 
 Department of Mathematics,  
Graduate School of Science\\ 
Osaka University \\ 1-1, \ Machikaneyama,\ Toyonaka,\ Osaka,\ 560-0043,\ 
Japan\\ E-mail: sumi@math.sci.osaka-u.ac.jp\\ 
http://www.math.sci.osaka-u.ac.jp/$\sim $sumi/} 
\date{November 3, 2010}
\maketitle 
\begin{abstract}
We investigate the dynamics of 
semigroups generated by a family of polynomial maps on the Riemann sphere 
such that the postcritical set in the complex plane is bounded. 
The Julia set of such a semigroup may not be connected in general. 
We show that for such a 
polynomial semigroup, 
if $A$ and $B$ are two connected components of the Julia set, 
then one of $A$ and $B$ surrounds the other. From this, 
it is shown that each connected component of the Fatou set is 
either simply or doubly connected. 
Moreover, we show that the Julia set of such a semigroup 
is uniformly perfect. 
An upper estimate of the cardinality of the set of all 
connected components of the Julia set of such a semigroup is given. 
By using this, we give 
a criterion for the Julia set 
to be connected.   
Moreover, we show that for any $n\in \Bbb{N} \cup \{ \aleph _{0}\} ,$  
there exists a finitely generated polynomial semigroup with bounded 
planar postcritical set 
such that the cardinality of the set of all connected components of the 
Julia set is equal to $n.$ 
Many new phenomena of polynomial semigroups that do not occur in the usual dynamics of polynomials  
are found and systematically investigated.
\end{abstract}
\section{Introduction}
 The theory of complex dynamical systems, which has 
 its origin in the important 
 work of Fatou and Julia in the 1910s, 
 has been investigated by many people and discussed in depth.  
In particular, since D. Sullivan showed the famous 
``no wandering domain theorem'' using 
Teichm\"{u}ller theory 
in the 1980s, 
this subject has 
attracted 
many researchers 
from a
wide area. 
For a general 
reference on complex dynamical systems, 
see Milnor's textbook \cite{M} or Beardon's textbook \cite{Be1}.    
 
 There are several 
areas
 in which we deal with 
  generalized notions of 
classical iteration theory of rational functions.   
One of them is the theory of 
dynamics of rational semigroups 
  (semigroups generated by a family of holomorphic maps on the 
  Riemann sphere $\CCI $), and another one is 
 the theory of   
 random dynamics of holomorphic maps on the Riemann sphere. 

In this paper, we will discuss 
the dynamics of rational semigroups. 

 A {\bf rational semigroup} is a semigroup 
generated by a family of non-constant rational maps on 
$\CCI $, where $\CCI $ denotes the Riemann sphere,
 with the semigroup operation being  
functional composition (\cite{HM1}). A 
{\bf polynomial semigroup} is a 
semigroup generated by a family of non-constant 
polynomial maps.
Research on the dynamics of
rational semigroups was initiated by
A. Hinkkanen and G. J. Martin (\cite{HM1,HM2}),
who were interested in the role of the
dynamics of polynomial semigroups while studying
various one-complex-dimensional
moduli spaces for discrete groups,
and
by F. Ren's group(\cite{ZR, GR}), 
 who studied 
such semigroups from the perspective of random dynamical systems.
Moreover, the research 
on 
rational semigroups is related to 
that 
on 
``iterated function systems" in 
fractal geometry.  
In fact, 
the Julia set of a rational semigroup generated by a 
compact family has 
`` backward self-similarity" 
(cf. Lemma~\ref{hmslem}-\ref{bss}). 
For 
other 
research 
on rational semigroups, see 
\cite{Sta1, Sta2, Sta3, SY, SSS, 
SS, SU1, SU2, SU3, SU5}, and \cite{S1}--\cite{Sstar2}. 

 The research 
on the  
dynamics of rational semigroups is also directly related to 
that 
on the 
random dynamics of holomorphic maps. 
The first 
study 
in this 
direction was 
by Fornaess and Sibony (\cite{FS}), and 
much research has followed. 
(See \cite{Br, Bu1, Bu2, 
BBR,GQL, S9, S10, S12, S8, Sstar1, Sstar2}.)   

 We remark that the complex dynamical systems 
 can be used to describe some mathematical models. For 
 example, the behavior of the population 
 of a certain species can be described as the 
 dynamical system of a polynomial 
 $f(z)= az(1-z)$ 
 such that $f$ preserves the unit interval and 
 the postcritical set in the plane is bounded 
 (cf. \cite{D}). It should also be remarked that 
 according to the change of the natural environment, 
 some species have several strategies to survive in the nature. 
From this point of view, 
 it is very important to consider the random 
 dynamics of such polynomials (see also Example~\ref{realpcbex}). 
For the random dynamics of polynomials on the unit interval, 
see \cite{Steins}. 
 
 We shall give some definitions 
for the 
dynamics of rational semigroups: 
\begin{df}[\cite{HM1,GR}] 
Let $G$ be a rational semigroup. We set
\[ F(G) := \{ z\in \CCI \mid G \mbox{ is normal in a neighborhood of  $z$} \} ,
\mbox{ and } J(G)  := \CCI \setminus  F(G) .\] \(  F(G)\) is  called the
{\bf Fatou set}  of  $G$ and \( J(G)\)  is  called the {\bf 
Julia set} of $G$. 
We 
let 
$\langle h_{1},h_{2},\ldots \rangle $ 
denote 
the 
rational semigroup generated by the family $\{ h_{i}\} .$
The Julia set of the semigroup generated by 
a single map $g$ is denoted by 
$J(g).$ 
\end{df}

\begin{df}\ 
\begin{enumerate}
\item 
For each rational map $g:\CCI \rightarrow \CCI $, 
we set 
$CV(g):= \{ \mbox{all critical values of }
g: \CCI \rightarrow \CCI \} .$ 
Moreover, for each polynomial map $g:\CCI \rightarrow \CCI $, 
we set $CV^{\ast }(g):= CV(g)\setminus \{ \infty \} .$ 
\item 
Let $G$ be a rational semigroup.
We set 
$ P(G):=
\overline{\bigcup _{g\in G} CV(g)} \ (\subset \CCI ). 
$ 
This is called the {\bf postcritical set} of $G.$
Furthermore, for a polynomial semigroup $G$,\ we set 
$P^{\ast }(G):= P(G)\setminus \{ \infty \} .$ This is 
called the {\bf planar postcritical set}
(or {\bf finite postcritical set}) 
 of $G.$
We say that a polynomial semigroup $G$ is 
{\bf postcritically bounded} if 
$P^{\ast }(G)$ is bounded in $\CC .$ 
\end{enumerate}
\end{df}
\begin{rem}
\label{pcbrem}
Let $G$ be a rational semigroup 
generated by a family $\Lambda $ of rational maps. 
Then, we have that 
$P(G)=\overline{\bigcup _{g\in G\cup \{ Id\} }\ g(\bigcup _{h\in \Lambda }CV(h))}$, 
where Id denotes the identity map on $\CCI $,  
and that $g(P(G))\subset P(G)$ for each $g\in G.$  
From this formula, one can figure out how the set 
$P(G)$ (resp. $P^{\ast }(G)$) spreads in $\CCI $ (resp. $\CC $). 
In fact, in Section~\ref{Const}, using the above formula, 
we present a way to construct examples of postcritically bounded 
polynomial semigroups (with some additional properties). Moreover, 
from the above formula, one may, in the finitely generated case, 
use a computer to see if a polynomial semigroup $G$ is postcritically bounded much in the same way 
as one verifies the boundedness of the critical orbit for the maps $f_{c}(z)=z^{2}+c.$   
\end{rem}
\begin{ex}
\label{realpcbex}
Let 
$\Lambda := \{ h(z)=cz^{a}(1-z)^{b}\mid 
a,b\in \NN  ,\ c>0,\  
c(\frac{a}{a+b})^{a}(\frac{b}{a+b})^{b}$ $\leq 1\} $ 
and let $G$ be the polynomial semigroup generated by 
$\Lambda .$ 
Since for each $h\in \Lambda $, 
$h([0,1])\subset [0,1]$ and 
$CV^{\ast }(h)\subset [0,1]$, 
it follows that each subsemigroup $H$ of $G$ is postcritically 
bounded. 
\end{ex}
\begin{rem}
\label{pcbound}
It is well-known that for a polynomial $g$ with 
$\deg (g)\geq 2$, 
$P^{\ast }(\langle g\rangle )$ is bounded in $\CC $ if and only if 
$J(g)$ is connected (\cite[Theorem 9.5]{M}).
\end{rem}
As mentioned in Remark~\ref{pcbound}, 
 the planar postcritical set is one 
piece of important information 
regarding the dynamics of polynomials. 
Concerning 
the theory of iteration of quadratic polynomials, 
we have been investigating the famous ``Mandelbrot set''.
    
When investigating the dynamics of polynomial semigroups, 
it is natural for us to 
discuss the relationship between 
  the planar postcritical set and the 
  figure of the Julia set.
The first question in this 
regard 
is: 
\begin{ques}
Let $G$ be a polynomial semigroup such that each 
element $g\in G$ is of degree at least two.
Is $J(G)$ necessarily connected when $P^{\ast }(G)$ is 
bounded in $\CC $?
\end{ques}
The answer is {\bf NO.}
\begin{ex}[\cite{SY}]
Let $G=\langle z^{3}, \frac{z^{2}}{4}\rangle .$ 
Then $P^{\ast }(G) =\{ 0\} $ 
(which is bounded in $\CC $)
and $J(G)$ is disconnected ($J(G)$ is a Cantor set 
of round circles). Furthermore,\ 
according to 
\cite[Theorem 2.4.1]{S5},  
it can be shown that 
a small					 
perturbation $H$ of $G$ 
 still satisfies that 
 $P^{\ast }(H) $ is 
 bounded in $\CC $  and that $J(H)$ is disconnected. 
 ($J(H)$ is a 
Cantor set of quasi-circles with uniform dilatation.)
\end{ex}
\begin{ques}
What happens if $P^{\ast }(G) $ is bounded in $\CC $ 
and $J(G)$ is disconnected? 
\end{ques}
\begin{prob}
Classify postcritically bounded polynomial semigroups.
\end{prob}
In this paper, we show that if $G$ is a postcritically bounded 
polynomial semigroup with disconnected Julia set, then 
$\infty \in F(G)$ (cf. Theorem~\ref{mainth2}-\ref{mainth2-2}), and 
for any two connected components of $J(G)$, one of them surrounds 
the other. This implies that 
there exists an intrinsic total order $``\leq "$ 
(called the 
``surrounding 
order") 
in the space ${\cal J}_{G}$ of 
connected components of $J(G)$, and that 
every connected component of $F(G)$ is either simply 
or doubly connected (cf. Theorem~\ref{mainth1}). 
 Moreover, for such a semigroup $G$, we show  
that the interior of ``the smallest filled-in Julia set'' $\hat{K}(G)$ 
is not empty, and that there exists a maximal element and a 
minimal element 
in the space ${\cal J}_{G}$ endowed with the order $\leq $  
(cf. Theorem~\ref{mainth2}).  
From these results, we obtain the result 
 that for a postcritically bounded polynomial semigroup $G$, 
 the Julia set $J(G)$ is uniformly perfect, 
 even if $G$ is not generated by a compact family of polynomials 
 (cf. Theorem~\ref{mainupthm}). 

 Moreover, 
 we utilize Green's functions with pole at infinity
to show that for a postcritically bounded 
 polynomial semigroup $G$, the cardinality of 
the set of all connected components of $J(G)$ is less than or equal to 
 that of $J(H)$, where $H$ is the ``real affine semigroup'' 
 associated with $G$ (cf. Theorem~\ref{polyandrathm1}). 
From this result, we obtain a sufficient condition for the Julia set 
of a postcritically bounded polynomial semigroup to be connected 
(cf. Theorem~\ref{polyandrathm2}). 
In particular, 
we show that if a postcritically bounded polynomial semigroup $G$ 
is generated by a family of quadratic polynomials, then $J(G)$ is connected 
(cf. Theorem~\ref{polyandrathm3}). 
 The proofs of the results in 
this and the previous paragraphs 
 are not straightforward. 
In fact, we first prove 
(1) that 
for any two connected components of $J(G)$ that 
 are included in $\CC $, one of them surrounds the other; 
 next, using (1) and the theory of Green's functions, we prove 
(2) that 
the cardinality of the set of all connected components of 
 $J(G)$ is less than or equal to that of $J(H)$, where 
 $H$ is the associated real affine semigroup; and finally, 
 using (2) and (1), we prove 
(3) that 
$\infty \in F(G)$, 
 int$(\hat{K}(G))\neq \emptyset $, and other results in the previous 
 paragraph.     

 Moreover, we show that for any $n\in \NN \cup \{ \aleph _{0}\} $, 
 there exists a finitely generated, postcritically bounded,  
 polynomial semigroup $G$ such that the cardinality of the set of 
 all connected components of $J(G)$ is equal to $n$ 
 (cf. Proposition~\ref{fincomp}, Proposition~\ref{countprop} and 
 Proposition~\ref{countcomp}). 
A sufficient condition for the cardinality of the set of all connected components 
of 
a Julia set 
to be equal to $\aleph _{0}$ is also given 
(cf. Theorem~\ref{countthm}). To obtain these results, we 
use 
the fact that 
the map induced by any element of a semigroup on the space of connected components of the Julia set preserves the order $\leq $ (cf. Theorem~\ref{mainth1}).   
  Note that 
this is in contrast to the dynamics of 
a single rational map $h$ or a non-elementary 
 Kleinian group, where it is known that either the Julia set is connected, or 
 the Julia set has uncountably many connected components.  
 Furthermore, in Section~\ref{Const} and Section~\ref{Poly}, we 
provide 
a way of constructing examples of 
 postcritically bounded polynomial semigroups 
 with 
some additional properties (disconnectedness of Julia set, 
semi-hyperbolicity, hyperbolicity, etc.) 
(cf. Proposition~\ref{Constprop}, Theorem~\ref{shshfinprop}, 
Theorem~\ref{sphypopen}). 
For example, by Proposition~\ref{Constprop}, 
there exists a $2$-generator 
postcritically bounded polynomial semigroup $G=\langle h_{1},h_{2}\rangle $ 
with disconnected 
Julia set such that $h_{1}$ has a Siegel disk.
        
As we see in Example~\ref{realpcbex} and Section~\ref{Const}, 
it is not difficult to construct many examples,  
it is not difficult to verify the hypothesis ``postcritically 
bounded'', 
and the class of postcritically bounded polynomial semigroups is 
very wide.  

 Throughout the paper, we will see many new phenomena in polynomial 
 semigroups that do not occur in 
 the usual dynamics of polynomials. Moreover, these new phenomena are 
 systematically investigated.
 
 In Section~\ref{Main}, we present the main results 
 of this paper. We give some tools in Section~\ref{Tools}. 
 The proofs of the main results are given in Section~\ref{Proofs}. 

\ 

There are many applications of the results of postcritically 
bounded polynomial semigroups in many directions. 
In the sequel \cite{S11}, by using the results in this paper, 
we investigate the fiberwise (sequencewise) and random dynamics of polynomials and the Julia sets. 
We present a sufficient condition for a fiberwise Julia set to be of measure zero, 
a sufficient condition for a fiberwise Julia set to be a Jordan curve, 
a sufficient condition  for a fiberwise Julia set to be a quasicircle, 
and a sufficient condition for a fiberwise Julia set to be a Jordan curve which is not a quasicircle.  
Moreover, using uniform fiberwise quasiconformal surgery on a fiber bundle, 
we show that for a $G\in {\cal G}_{dis}$, 
there exist families of uncountably many mutually disjoint quasicircles with uniform 
dilatation which are parameterized by the Cantor set, densely inside $J(G).$ 
In the sequel \cite{S12}, we classify hyperbolic or semi-hyperbolic postcritically bounded compactly generated 
polynomial semigroups, 
in terms of the random complex dynamics. It is shown that in one of the classes, 
for almost every sequence $\gamma $, the Julia set $J_{\gamma }$ of $\gamma $ is a Jordan curve but not 
a quasicircle, the unbounded component of $\CCI \setminus J_{\gamma }$ is a John domain, and 
the bounded component of $\CC \setminus J_{\gamma }$ is not a John domain. 
Moreover, in \cite{S12,S11}, we find many examples with this phenomenon. 
Note that this phenomenon does not hold 
in the usual iteration dynamics of a single polynomial map $g$ with $\deg (g)\geq 2.$  
 In the sequel \cite{S8, Snew}, we 
 investigate the 
Markov process on $\CCI $ associated with 
the random dynamics of polynomials and 
we  consider the probability $T_{\infty }(z)$ 
of tending to $\infty \in \CCI $ 
starting with the initial value $z\in \CCI .$ 
Applying many results of this paper, 
it will be shown in \cite{Snew} that 
if the associated polynomial semigroup $G$ 
is postcritically bounded and the Julia set is 
disconnected, then the function $T_{\infty }$ defined on $\CCI $ 
has many interesting properties which are 
similar to those of the Cantor function. 
In fact, under certain conditions, $T_{\infty }$ is continuous on $\CCI  $ 
and varies precisely on the Julia set, of which Hausdorff dimension is strictly less than two.  
(For example, if we consider the random dynamics generated by two 
polynomials $h_{1}:=g_{1}^{2}, h_{2}:=g_{2}^{2},$ where $g_{1}(z):=z^{2}-1,g_{2}(z):=z^{2}/4$, then 
$T_{\infty }$ is continuous on $\CCI $ and $T_{\infty }$ varies precisely on the Julia set (Figure~\ref{fig:dcgraph}) of 
the semigroup generated by $h_{1},h_{2}$. See \cite{S8,S9}.)  
Such a kind of ``singular functions on the 
complex plane'' appear very naturally in  
random dynamics of polynomials, and the 
results of this paper (for example, 
the results on the space of all connected 
components of a Julia set) are the keys to 
investigating that. 
(The above results have been announced in \cite{S9, S10, Sstar1}.)  

 Moreover, as illustrated before, 
 it is very important for us to recall that 
 the complex dynamics can be applied to describe some mathematical models. 
 For example, the behavior of the population of a 
 certain species can be described as the dynamical systems
  of a polynomial $h$ such that $h$ preserves the unit interval  
 and the postcritical set in the plane is bounded. 
 When one considers such a model, 
 it is very natural to consider the random dynamics of 
 polynomial with bounded postcritical set in the plane 
 (see Example~\ref{realpcbex}).    

 In the sequel \cite{SS}, we give some 
further results on postcritically 
 bounded polynomial semigroups, by using many results in this paper and \cite{S11, S12}. 
 Moreover, in the sequel \cite{S7}, 
 we define a new kind of cohomology theory, in order to 
 investigate the action of finitely generated semigroups (iterated function systems), and 
we  apply it to the study of 
the dynamics of postcritically bounded finitely generated polynomial semigroups $G$.  
In particular, by using this new cohomology theory, we can 
describe the space ${\cal J}_{G}$ of connected components of Julia sets of $G$, 
we can give some estimates on the cardinality of ${\cal J}_{G}$, 
and we can give a sufficient condition for the cardinality of the space 
of connected components of the Fatou set of $G$ to be infinity. 
In \cite{S8,Sstar2,Scoop}, we investigate the random complex dynamics and the dynamics of 
transition operator, by developing the theory of random complex dynamics and that of dynamics of 
rational semigroups, simultaneously. It is shown that regarding the random dynamics of complex polynomials, 
generically the chaos of the averaged system disappears due to the cooperation of the generators, even though
 each map itself in the system has a chaotic part. We call this phenomenon ``cooperation principle''.   
Moreover, we see that under certain conditions, in the limit state, 
complex analogues of singular functions (continuous functions on $\CCI $ which 
vary only on the Julia set of associated rational semigroup $G$) naturally appear. 
The above function $T_{\infty }$ is a typical example of this complex analogue of singular function. 
   
\ 

\noindent {\bf Acknowledgement:} 
The author 
thanks 
R. Stankewitz for many 
valuable comments.  
  \section{Main results}
\label{Main}
In this section we present the statements of the main results.
Throughout this paper, we deal with semigroups $G$ that 
might not 
be generated by a compact family of polynomials.
The proofs are given in Section~\ref{Proofs}.

\subsection{Space of connected components of a Julia set, surrounding order}
\label{concompsec}
We present some results 
concerning the connected components of the 
Julia set of a postcritically bounded polynomial semigroup.
 The proofs are 
given in Section~\ref{pfconcompsec}.

The following theorem generalizes \cite[Theorem 1]{SY}. 
\begin{thm}
\label{mainth0}
Let $G$ be a rational semigroup 
generated by a family $\{ h_{\lambda }\} _{\lambda \in 
\Lambda }.$ 
Suppose 
that  
there exists a connected component 
$A$ of $J(G)$ such that $\sharp A> 1$ and 
$\bigcup _{\lambda \in \Lambda }
J(h_{\lambda })\subset A.$ Moreover, 
suppose that for any $\lambda \in \Lambda $ such that 
$h_{\lambda }$ is a M\"{o}bius transformation of finite order,\ 
we have $h_{\lambda }^{-1}(A)\subset A.$
Then,\ 
$J(G)$ is connected.
\end{thm} 
\begin{df}
We set 
Rat : $=\{ h:\CCI \rightarrow \CCI \mid 
h \mbox { is a non-constant rational map}\} $
endowed with the topology induced by uniform convergence on $\CCI $ 
with respect to the spherical distance.   
We set 
Poly :$=\{ h:\CCI \rightarrow \CCI 
\mid h \mbox{ is a non-constant polynomial}\} $ endowed with 
the relative topology from Rat.   
Moreover, we set 
Poly$_{\deg \geq 2}
:= \{ g\in \mbox{Poly}\mid \deg (g)\geq 2\} $ 
endowed with the relative topology from 
Rat.  
\end{df}
\begin{rem}
Let $d\geq 1$,  $\{ p_{n}\} _{n\in \NN }$ a 
sequence of polynomials of degree $d$, 
and $p$ a polynomial.  
Then, $p_{n}\rightarrow p$ in Poly if and only if 
the coefficients converge appropriately and $p$ is of degree $d.$ 
\end{rem}
\begin{df} 
Let ${\cal G} $ be the set of all polynomial semigroups 
$G$ with the following 
properties:
\begin{itemize}
\item 
 each element of $G$ is of degree 
at least two, and  
\item  $P^{\ast }(G)$ is 
bounded in $\CC $, i.e., $G$ is postcritically bounded.
\end{itemize}   
Furthermore, we set 
${\cal G}_{con}=
\{ G\in {\cal G}\mid 
J(G)\mbox{ is connected}\} $ and 
${\cal G}_{dis}=
\{ G\in {\cal G}\mid 
J(G)\mbox{ is disconnected}\}.$ 
\end{df}  
\noindent {\bf Notation:}
For a polynomial semigroup $G$,\ 
we denote by 
${\cal J}={\cal J}_{G}$ the set of all 
connected components $J$ 
of $J(G)$ such that $J\subset \CC .$   
Moreover, we denote by 
$\hat{{\cal J}}=\hat{{\cal J}}_{G}$ the set of all connected components 
of $J(G).$ 
\begin{rem}
\label{hatjcptrem}
If a polynomial semigroup $G$ is generated by a compact set 
in Poly$_{\deg \geq 2}$, then 
$\infty \in F(G)$ and thus ${\cal J}=\hat{{\cal J}}.$ 
\end{rem}
\begin{df}
For any connected sets $K_{1}$ and 
$K_{2}$ in $\CC ,\ $  ``$K_{1}\leq K_{2}$'' indicates that 
$K_{1}=K_{2}$, or $K_{1}$ is included in 
a bounded component of $\CC \setminus K_{2}.$ Furthermore, 
``$K_{1}<K_{2}$'' indicates $K_{1}\leq K_{2}$ 
and $K_{1}\neq K_{2}.$ Note that 
``$ \leq $'' is a partial order in 
the space of all non-empty compact connected 
sets in $\CC .$ This ``$\leq $" is called 
the {\bf surrounding order.} 
\end{df}

\begin{figure}[htbp]
\caption{The Julia set of $G=\langle g_{1}^{2},g_{2}^{2}\rangle $, 
where $g_{1}(z):=z^{2}-1,\ g_{2}(z):=\frac{z^{2}}{4}$. 
$G\in {\cal G}_{dis}, $ $G$ is hyperbolic, and 
$\sharp (\hat{{\cal J}}_{G})>\aleph _{0}.$ }
\ \ \ \ \ \ \ \ \ \ \ \ \ \ \ \ \ \ \ \ \ \ \ \ \ \ \ \ \ \ \ \ \ \ \ \ \ \ \ \ \ \ \ 
\ \ \ \ \ \ \ 
\includegraphics[width=3cm,width=3cm]{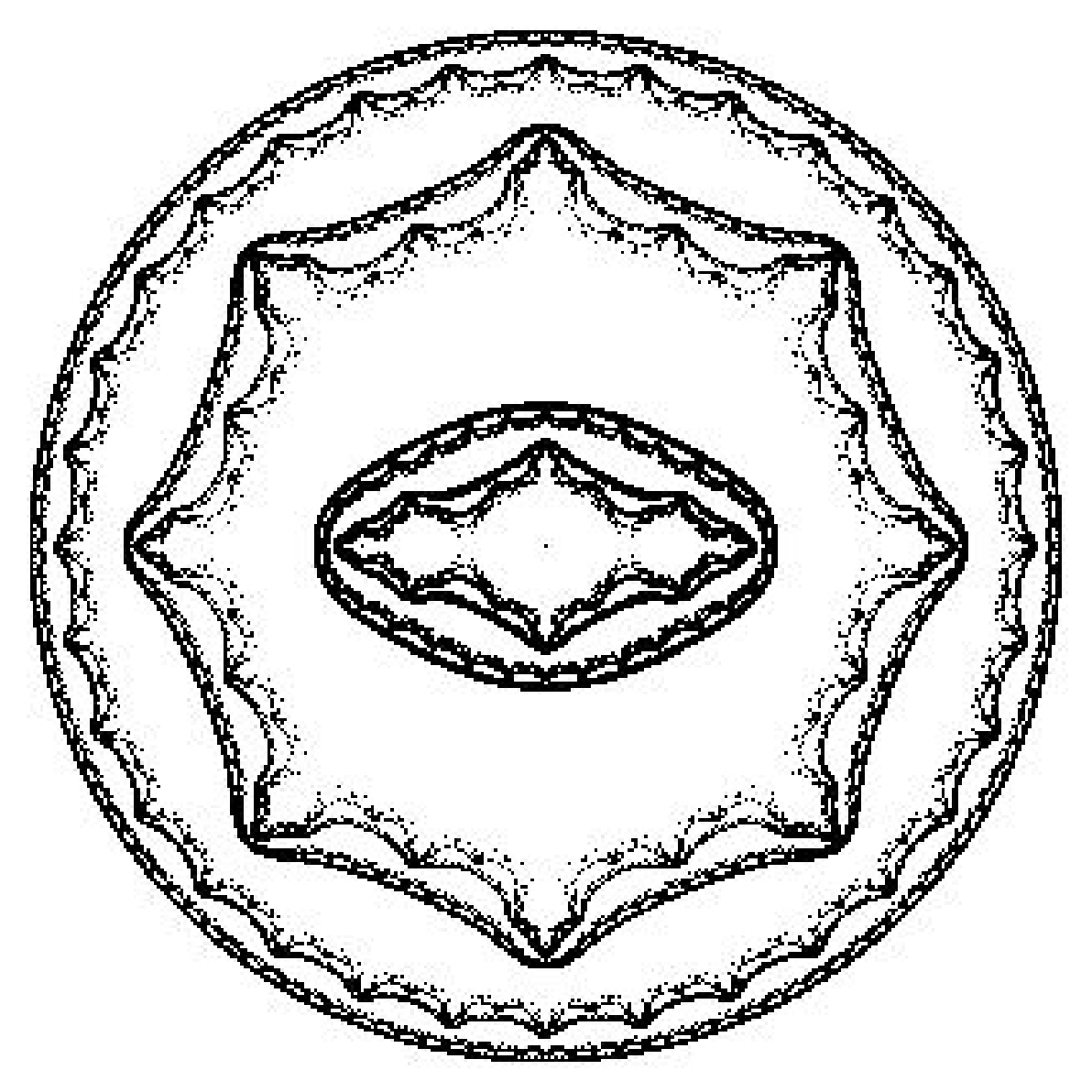}
\label{fig:dcgraph}
\end{figure}

\begin{thm}
\label{mainth1}
Let $G\in {\cal G}$ (possibly generated by a non-compact 
family). 
Then we have all of the following. 
\begin{enumerate}
\item \label{mainth1-1}
$({\cal J},\ \leq )$ is totally ordered. 
\item \label{mainth1-2}
Each connected component of 
$F(G)$ is either simply or doubly connected. 
\item \label{mainth1-3}
For any $g\in G$ and any connected component 
$J$ of $J(G)$,\ we have that  
$g^{-1}(J)$ is connected. 
Let $g^{\ast }(J)$ be the connected component of 
$J(G)$ containing $g^{-1}(J).$ 
If $J\in {\cal J}$, then 
$g^{\ast }(J)\in {\cal J}.$    
If $J_{1},J_{2}\in {\cal J} $ and $J_{1}\leq J_{2},\ $ then 
$g^{-1}(J_{1})\leq g^{-1}(J_{2})$ 
and $g^{\ast }(J_{1})\leq g^{\ast }(J_{2}).$
\end{enumerate}
\end{thm} 
For the figures of the Julia sets of semigroups $G\in {\cal G}_{dis}$, 
see figure~\ref{fig:dcgraph} and figure~\ref{fig:3mapcountjulia2}.

\subsection{Upper estimates of $\sharp (\hat{{\cal J}})$}
\label{Upper}
Next, we present some results on the space $\hat{{\cal J}}$
and some results on upper estimates of 
$\sharp (\hat{{\cal J}}).$ The proofs are given in 
Section~\ref{Proof of Upper} and Section~\ref{Proof of Properties}.
\begin{df}
\ 
\begin{enumerate}
\item For a polynomial $g$, we denote by 
$a(g)\in \CC $ the coefficient of 
the highest degree term of $g.$ 
\item 
We set 
RA $:=\{ ax+b\in \RR [x]\mid a,b\in \RR ,a\neq 0\} $ 
endowed with the topology such that, 
$a_{n}x+b_{n}\rightarrow ax+b$ if and only if 
$a_{n}\rightarrow a$ and $b_{n}\rightarrow b.$ 
 The space RA is a semigroup 
with the semigroup operation being functional composition.
Any subsemigroup of RA will be called a {\em real affine semigroup}.
We define a map $\Psi :$ Poly $\rightarrow $ RA as follows: 
For a polynomial $g\in $ Poly,   
we set $\Psi (g)(x):= \deg (g)x+\log | a(g)|.$ 

 Moreover, for a polynomial semigroup $G$, 
 we set $\Psi (G):= \{ \Psi (g)\mid g\in G\}  $ ($\subset $ RA). 

\item 
We set $\hat{\RR }:= \RR \cup \{ \pm \infty \} $ endowed with 
the topology such that 
$\{ (r,+\infty ]\} _{r\in \RR }$ makes a fundamental neighborhood system of 
$+\infty $, and such that $\{ [-\infty ,r)\} _{r\in \RR }$ makes a 
fundamental neighborhood system of $-\infty .$ 
For a real affine semigroup $H$, we set 
$$M(H):= \overline{ \{ x\in \RR \mid \exists h\in H,  h(x)=x, |h'(x)|>1\} } 
\ (\subset \hat{\RR })
,$$
where the closure is taken in the space $\hat{\RR }.$   
Moreover, we denote by ${\cal M}_{H}$ the set of all connected components of 
$M(H).$ 
\item 
We denote by $\eta : $ RA $\rightarrow $ Poly the natural embedding 
defined by $\eta (x\mapsto ax+b)=(z\mapsto az+b)$, where 
$x\in \RR $ and $z\in \CC .$  
\item 
We define a map $\Theta :$ Poly $\rightarrow $ Poly  as follows. 
For a polynomial $g$, we set 
$\Theta (g)(z)=a(g)z^{\deg (g)}.$ 
Moreover, for a polynomial semigroup $G$, 
we set $\Theta (G):= \{ \Theta (g)\mid g\in G\} .$
\end{enumerate} 

\end{df}
\begin{rem}
\ 
\begin{enumerate}
\item 
The map $\Psi :$ Poly $\rightarrow $ RA is a semigroup homomorphism. 
That is, we have $\Psi (g\circ h)=\Psi (g)\circ \Psi (h).$
Hence, for a polynomial semigroup $G$, 
the image $\Psi (G)$ is a real affine semigroup. 
Similarly, the map $\Theta :$ Poly $\rightarrow $ Poly 
is a semigroup homomorphism. Hence, 
for a polynomial semigroup $G$, 
the image $\Theta (G)$ is a polynomial semigroup.  
\item 
The maps $\Psi :$ Poly $\rightarrow $ RA,  
$\eta :$ RA $\rightarrow $ Poly, 
and $\Theta :$ Poly $\rightarrow $ Poly are continuous.
\end{enumerate}
\end{rem}
\begin{df}
For any connected sets $M_{1}$ and $M_{2}$ in $\hat{\RR  }$, 
``$M_{1}\leq _{r} M_{2}$'' indicates that $M_{1}=M_{2}$, or 
each $(x,y)\in M_{1}\times M_{2}$ satisfies 
$x<y.$ Furthermore, 
``$M_{1}< _{r} M_{2}$'' indicates $M_{1}\leq _{r}M_{2}$ and 
$M_{1}\neq M_{2}.$ 
\end{df}
\begin{rem}
The above ``$\leq _{r}$" is a partial order in the space of 
non-empty connected subsets of $\hat{\RR }.$ 
Moreover, for each real affine semigroup $H$, 
$({\cal M}_{H},\leq _{r})$ is totally ordered.
\end{rem}
The following theorem gives us some upper estimates of $\sharp (\hat{{\cal J}}_{G}).$ 
\begin{thm}
\label{polyandrathm1}
\ 
\begin{enumerate}
\item \label{polyandrathm1-1}
Let $G$ be a polynomial semigroup in ${\cal G}.$ 
Then, we have 
$\sharp (\hat{{\cal J}}_{G})\leq \sharp ({\cal M}_{\Psi (G)}).$ 
More precisely, there exists an injective map 
$\tilde{\Psi }: \hat{{\cal J}}_{G}\rightarrow 
{\cal M}_{\Psi (G)}$ such that if 
$J_{1},J_{2}\in {\cal J}_{G}$ and $J_{1}<J_{2}$, 
then 
$\tilde{\Psi }(J_{1})<_{r}\tilde{\Psi }(J_{2}).$ 
\item \label{polyandrathm1-2}
If $G\in {\cal G}_{dis}$, then we have that 
$M(\Psi (G))\subset \RR $ and $M(\Psi (G))=J(\eta (\Psi (G))).$  
\item \label{polyandrathm1-3}
Let $G$ be a polynomial semigroup in ${\cal G}.$ Then, 
$\sharp (\hat{{\cal J}}_{G})\leq 
\sharp ({\hat{\cal J}}_{\eta (\Psi (G))}).$
\end{enumerate}
\end{thm}
\begin{cor}
\label{polyandracor}
Let $G$ be a polynomial semigroup in ${\cal G}.$ 
Then, we have 
$\sharp ({\hat{\cal J}}_{G})\leq \sharp ({\hat{\cal J}}_{\Theta (G)}).$ 
More precisely, there exists an injective map 
$\tilde{\Theta }: \hat{{\cal J}}_{G}\rightarrow 
\hat{{\cal J}}_{\Theta (G)}$ such that 
if $J_{1},J_{2}\in {\cal J}_{G}$ and $J_{1}<J_{2}$, 
then $\tilde{\Theta }(J_{1})\in {\cal J}_{\Theta (G)}$, 
 $\tilde{\Theta }(J_{2})\in {\cal J}_{\Theta (G)}$, and 
 $\tilde{\Theta }(J_{1})<\tilde{\Theta }(J_{2}).$  
\end{cor}
The following three theorems give us sufficient conditions for the Julia set of a $G\in {\cal G}$ 
to be connected. 
\begin{thm}
\label{polyandrathm2}
Let $G=\langle h_{1},\ldots ,h_{m}\rangle $ be a finitely generated 
polynomial semigroup in ${\cal G}.$ 
For each $j=1,\ldots ,m$, 
let $a_{j}$ be the coefficient of the highest degree term 
of polynomial $h_{j}.$ Let $\alpha := 
\min  _{j=1,\ldots ,m}\{ \frac{-1}{\deg (h_{j})-1}\log |a_{j}|\} $ 
and $\beta := 
\max _{j=1,\ldots ,m}\{ \frac{-1}{\deg (h_{j})-1}\log |a_{j}|\} .$ 
We set $[\alpha ,\beta ]:= \{ x\in \RR \mid \alpha\leq x\leq \beta \} .$  
If 
$[\alpha ,\beta ]\subset \bigcup _{j=1}^{m}\Psi (h_{j})^{-1}([\alpha ,\beta ])$,  
then $J(G)$ is connected. 
\end{thm}
\begin{thm}
\label{polyandrathm3}
Let $G$ be a 
polynomial semigroup in ${\cal G}$ generated by a 
(possibly non-compact) 
family of 
polynomials of degree two. Then, $J(G)$ is connected.
\end{thm}
\begin{thm}
\label{polyandrathm4}
Let $G$ be a polynomial semigroup in ${\cal G}$ generated by 
a (possibly non-compact) family $\{ h_{\lambda }\} _{\lambda \in \Lambda }$ of 
polynomials. 
Let $a_{\lambda }$ be the coefficient of the highest degree term of 
the polynomial $h_{\lambda }.$ 
Suppose that for any  $\lambda ,\xi \in \Lambda $, 
we have $(\deg (h_{\xi })-1)\log |a_{\lambda }|=
(\deg (h_{\lambda })-1)\log |a_{\xi }|.$ Then, 
$J(G)$ is connected.
\end{thm}
\begin{rem}
\label{r:intmev}
In \cite{S7}, a new cohomology theory for (backward) self-similar systems (iterated function 
systems) was 
introduced by the author of this paper. By using this new cohomology theory,  
for a postcritically bounded finitely generated polynomial semigroup $G$,  
we can describe the space of connected components of $G$ and 
we can give 
some estimates on $\sharp  ({\cal J}_{G})$ and $\sharp ({\cal M}_{\Psi (G)}).$  
\end{rem}

\subsection{Properties of ${\cal J}$}
\label{Properties}
In this section, we present some results on 
${\cal J}.$ The proofs are given in Section~\ref{Proof of Properties}.
\begin{df}
 For a polynomial semigroup $G$,\ we set 
$$ \hat{K}(G):=\{ z\in \CC 
\mid \bigcup _{g\in G}\{ g(z)\} \mbox{ is bounded in }\CC \} $$ 
and call $\hat{K}(G)$ the {\bf smallest filled-in Julia set} of 
$G.$ 
For a polynomial $g$, we set $K(g):= \hat{K}(\langle g\rangle ).$ 

\end{df}
\noindent {\bf Notation:} 
For a set $A\subset \CCI $, we denote by int$(A)$ the set of 
all  
interior points of $A.$ 
\begin{prop}
\label{fcprop}
Let $G\in {\cal G}.$  If $U$ is a connected component 
of $F(G)$ such that $U\cap \hat{K}(G)\neq \emptyset $, 
then $U\subset $ {\em int}$(\hat{K}(G))$ and $U$ is 
simply connected. Furthermore, 
we have $\hat{K}(G)\cap F(G)=$ {\em int}$(\hat{K}(G)).$

\end{prop}
%
%
%
%
%
%
%
%
%
\noindent {\bf Notation:} 
For a polynomial semigroup $G$ with 
$\infty \in F(G)$, we denote by 
$F_{\infty }(G)$ the connected component of $F(G)$ containing 
$\infty .$ Moreover, for a polynomial $g$ with 
$\deg (g)\geq 2$, we set $F_{\infty }(g):= 
F_{\infty }(\langle g\rangle ).$ 

\ 

The following theorem is the key to obtaining further results of postcritically bounded 
polynomial semigroups in this paper, 
and those of related random dynamics of polynomials in the sequel \cite{S11,Snew}. 
We remark that   
Theorem~\ref{mainth2}-\ref{mainth2-4} generalizes \cite[Theorem 2]{SY}. 
\begin{thm}
\label{mainth2}
Let $G\in {\cal G} _{dis}$ (possibly generated by a non-compact family). Then, under the above notation,\ 
we have the following.
\begin{enumerate}
\item \label{mainth2-2}
We have that $\infty \in F(G)$ (thus ${\cal J}=\hat{{\cal J}}$) and the  
 connected component 
$F_{\infty }(G)$ of $F(G)$ containing $\infty $ 
is simply connected. 
Furthermore,\ 
the element $J_{\max }=J_{\max}(G)\in {\cal J}$  
containing $\partial F_{\infty }(G)$ 
is the unique element of ${\cal J}$ satisfying that 
$J\leq J_{\max }$ for each 
$J\in {\cal J}.$  
\item 
\label{mainth2-3}
There exists a unique element 
$J_{\min }=J_{\min }(G)\in {\cal J}$ such that 
$J_{\min }\leq J$ for 
  each element $J\in {\cal J}. $
  Furthermore, let $D$ be the unbounded 
  component of $\CC \setminus J_{\min }. $ 
Then,     
$ P^{\ast }(G) \subset \hat{K}(G)\subset 
  \CC \setminus D $ and 
  $\partial \hat{K}(G)\subset J_{\min }. $
\item \label{mainth2-3b}
If $G$ is generated by a family 
$\{ h_{\lambda }\} _{\lambda \in \Lambda },$ 
then there exist two elements $\lambda _{1}$ 
and $\lambda _{2}$ of $\Lambda $ satisfying: 
  \begin{enumerate}
  \item there exist two elements $J_{1}$ and 
  $J_{2}$ of ${\cal J} $ with  
  $J_{1}\neq J_{2}$ such that $J(h_{\lambda _{i}})
  \subset J_{i}$ for each $i=1,2$;  
  \item $J(h_{\lambda _{1}})\cap J_{\min }=\emptyset $; 
  \item for each $n\in \NN $, 
  we have $h_{\lambda _{1}}^{-n}(J(h_{\lambda _{2}}))
  \cap J(h_{\lambda _{2}})=\emptyset $ and 
  $h_{\lambda _{2}}^{-n}(J(h_{\lambda _{1}}))
  \cap J(h_{\lambda _{1}})=\emptyset $; and 
  \item $h_{\lambda _{1}}$ has an attracting 
  fixed point $z_{1}$ in $\CC $,    
{\em int}$(K(h_{\lambda _{1}}))$ consists of 
  only one immediate attracting basin 
  for $z_{1}$,   
  and 
  $K(h_{\lambda _{2}})\subset $ {\em int}$(K(h_{\lambda _{1}})).$ 
  Furthermore, $z_{1}\in $ {\em int}$(K(h_{\lambda _{2}})).$ 
  \end{enumerate}
\item 
\label{mainth2ast1}
For each $g\in G$ with 
  $J(g)\cap J_{\min }=\emptyset $,\ 
  we have that $g$ has an attracting 
  fixed point $z_{g}$ in $\CC $, {\em int}$(K(g))$ consists of 
  only one immediate attracting basin for $z_{g}$,   
  and  $J_{\min }\subset $ {\em int}$(K(g)).$ 
  Note that 
it is not necessarily true that $z_{g}=z_{f}$ when 
$g,f\in G$ are such that $J(g)\cap J_{\min }=\emptyset $ 
and $J(f)\cap J_{\min }=\emptyset $ (see Proposition~\ref{fincomp}).  
\item 
\label{mainth2-4}
We have that  $\mbox{{\em int}}(\hat{K}(G))\neq 
\emptyset .$ Moreover,  

  \begin{enumerate}
  \item \label{mainth2-4-1}
  $\CC \setminus J_{\min }$ is 
  disconnected, $\sharp J\geq 2$ for each 
  $J\in \hat{{\cal J}}$, and  
  
  \item \label{mainth2-4-2}
  for each $g\in G$ with 
  $J(g)\cap J_{\min }=\emptyset $,\ 
  we have that $J_{\min }<g^{\ast }(J_{\min })$, 
  $g^{-1}(J(G))\cap J_{\min }=\emptyset $,  
  $g(\hat{K}(G)\cup J_{\min })\subset $ {\em int}$(\hat{K}(G))$,    
  and   
  the unique attracting fixed point $z_{g}$  
  of $g$ in $\CC $ belongs to 
  $\mbox{{\em int}}(\hat{K}(G)).$ 
  \end{enumerate} 
  
\item 
\label{mainth2ast2}
Let ${\cal A} $ be the set of all doubly connected components 
of $F(G).$ Then, $\bigcup _{A\in {\cal A}}A\subset \CC $ and 
$({\cal A},\leq )$ is totally ordered.    
\end{enumerate}
\end{thm}
We present a result on uniform perfectness of the Julia sets 
of semigroups in ${\cal G}.$ 
\begin{df}
A compact set $K$ in $\CCI $ is said to be 
uniformly perfect if 
$\sharp K\geq 2$ and there exists a constant $C>0$ such that 
each annulus $A$ that separates $K$ satisfies that 
mod $A<C$, where mod $A$ denotes the modulus of $A$ 
(See the definition in \cite{LV}). 
\end{df}  
\begin{thm}
\label{mainupthm}
\ 
\begin{enumerate}
  \item \label{mainupthm1}
  Let $G$ be a polynomial semigroup in $ {\cal G}.$ Then,  
  $J(G) $ is uniformly perfect. Moreover, if $z_{0}\in J(G)$ is a 
  superattracting fixed point of an element of $G$, then 
  $z_{0}\in $ {\em int}$(J(G)).$ 
\item 
\label{mainupthm2} 
If $G\in {\cal G}$ and $\infty \in J(G)$, then 
$G\in {\cal G}_{con}$ and $\infty \in $ {\em int}$(J(G)).$  
  \item \label{mainupthm3}
  Suppose that $G\in {\cal G}_{dis}.$   
Let $z_{1}\in J(G)\cap \CC $ be a superattracting 
fixed point of $g\in G.$   
Then 
$z_{1}\in $ {\em int}$(J_{\min })$ and 
  $J(g)\subset J_{\min }.$  
\end{enumerate}

\end{thm}
We remark that in \cite{HM2}, it was shown that there exists a rational semigroup $G$ such that 
$J(G)$ is not uniformly perfect. 

 We now present results on the Julia sets of subsemigroups of an element of ${\cal G}_{dis}.$ 
 \begin{prop}
\label{orderjprop}
Let $G\in {\cal G}_{dis}$ 
and let $J_{1}, J_{2}\in {\cal J}={\cal J}_{G}$  
with $J_{1}\leq J_{2}.$ 
%
%
Let $A_{i}$ be the unbounded component of 
$\CC \setminus J_{i}$ for each $i=1,2.$   
Then, we have the following.
\begin{enumerate}
\item \label{orderjprop1}
Let $Q_{1}=\{ g\in G\mid \exists J \in {\cal J} \mbox{with }
J_{1}\leq J,\ J(g)\subset J\} $ 
and let $H_{1}$ be the subsemigroup of $G$ generated by 
$Q_{1}.$ Then 
$J(H_{1})\subset J_{1}\cup A_{1}.$ 
\item \label{orderjprop2}
Let $Q_{2}=\{ g\in G\mid \exists J \in {\cal J} \mbox{with }
J\leq J_{2},\ J(g)\subset J\} $ 
and let $H_{2}$ be the subsemigroup of $G$ generated by 
$Q_{2}.$ Then 
$J(H_{2})\subset \CC \setminus A_{2}.$ 
\item \label{orderjprop3} 
Let $Q=\{ g\in G\mid \exists J \in {\cal J} \mbox{with }
J_{1}\leq J\leq J_{2},\ J(g)\subset J\} $ 
and let $H$ be the subsemigroup of $G$ generated by 
$Q.$ Then 
$J(H)\subset J_{1}\cup (A_{1}\setminus A_{2}).$
\end{enumerate}
\end{prop}

\begin{prop}
\label{bminprop}
Let $G$ be a 
polynomial 
semigroup generated by a compact subset 
$\G $ of {\em Poly}$_{\deg \geq 2}.$ Suppose that $G\in {\cal G}_{dis}.$ Then,    
there exists an element 
$h_{1}\in \G $ with 
 $J(h_{1})\subset J_{\max } $ and 
there exists an element 
$h_{2}\in \G $ with 
$J(h_{2})\subset J_{\min }.$ 
\end{prop}

\subsection{Finitely generated polynomial semigroups $G\in {\cal G}_{dis}$ such that 
$2\leq \sharp (\hat{{\cal J }}_{G})\leq \aleph _{0}$}
\label{Poly}
In this section, we present some results on various 
finitely generated polynomial semigroups $G\in {\cal G}_{dis}$ such that 
$2\leq \sharp (\hat{{\cal J}}_{G})\leq \aleph _{0}.$ 
The proofs are given in Section~\ref{Proof of Poly}.

It is well-known that for a rational map $g$ 
with $\deg (g)\geq 2$, if $J(g)$ is 
disconnected, then $J(g)$ has uncountably 
many connected components 
(See \cite{M}). 
Moreover, if $G$ is a non-elementary 
Kleinian group with disconnected 
Julia set (limit set), then  
$J(G)$ has uncountably many connected components. 
However,  
for general rational semigroups, 
we have the following examples.
\begin{thm}
\label{fcthm}
Let $G$ be a polynomial semigroup in ${\cal G}$ 
generated by a (possibly non-compact) family $\G $ 
in {\em Poly}$_{\deg \geq 2}.$  
Suppose that there exist mutually distinct 
elements $J_{1},\ldots , J_{n}\in \hat{{\cal J}}_{G}$ 
such that for each $h\in \G $ and each $j\in \{ 1,\ldots ,n\} $,  
there exists an element $k\in \{ 1,\ldots ,n\} $ with 
$h^{-1}(J_{j})\cap J_{k}\neq \emptyset .$ 
Then, we have $\sharp (\hat{{\cal J}}_{G})=n.$   
\end{thm}

\begin{prop}
\label{fincomp}
For any $n\in \NN  $ with $n>1$, there exists a 
finitely generated 
polynomial semigroup $G_{n}=\langle h_{1},\ldots ,h_{2n}\rangle $ 
in ${\cal G}$ satisfying 
$\sharp (\hat{{\cal J}}_{G_{n}})=n.$ 
In fact, let $0<\epsilon <\frac{1}{2}$ and 
we set for each $j=1,\ldots ,n,\ a_{j}(z):=\frac{1}{j} z^{2}$ 
and $\beta _{j}(z):=\frac{1}{j}(z-\epsilon )^{2}+\epsilon .$ 
Then, for any sufficiently large $l\in \NN $, there exists an open neighborhood 
$V$ of $(\alpha _{1}^{l},\ldots ,\alpha _{n}^{l}, 
\beta _{1}^{l},\ldots ,\beta _{n}^{l})$ in  
 {\em (Poly)}$^{2n}$ such that for any 
 $(h_{1},\ldots ,h_{2n})\in V$,  
 the semigroup $G=\langle h_{1},\ldots ,h_{2n}\rangle $ satisfies that 
 $G\in {\cal G}$ and 
 $\sharp (\hat{{\cal J}}_{G})=n.$    

\end{prop}
\begin{thm}
\label{countthm}
Let $G=\langle h_{1},\ldots ,h_{m}\rangle \in {\cal G}_{dis}$ be a 
polynomial semigroup with $m\geq 3.$ Suppose 
that there exists an element $J_{0}\in \hat{{\cal J}}$ 
such that $\bigcup _{j=1}^{m-1}J(h_{j})\subset J_{0}$, and 
such that for each $j=1,\ldots ,m-1, $ 
we have $h_{j}^{-1}(J(h_{m}))\cap J_{0}\neq \emptyset .$ 
Then, 
we have all of the following.
\begin{enumerate}
\item $\sharp (\hat{{\cal J}})=\aleph _{0}.$ 
\item $J_{0}=J_{\min }$, or $J_{0}=J_{\max }.$ 
\item If $J_{0}=J_{\min }$, then $J_{\max }=J(h_{m})$,  
$J(G)=J_{\max}\cup \bigcup _{n\in \NN \cup \{ 0\} }
(h_{m})^{-n}(J_{\min })$, 
and for any $J\in \hat{{\cal J}}$ with $J\neq J_{\max}$, 
there exists no sequence $\{ C_{j}\} _{j\in \NN } $ 
of mutually distinct elements of $\hat{{\cal J}} $ such that 
$\min _{z\in C_{j}}d(z,J)\rightarrow 0$ as $j\rightarrow \infty .$ 

\item 
 If 
$J_{0}=J_{\max }$, then $J_{\min }=J(h_{m})$,  
$J(G)=J_{\min}\cup \bigcup _{n\in \NN \cup \{ 0\} }
(h_{m})^{-n}(J_{\max })$, and 
for any $J\in \hat{{\cal J}}$ with $J\neq J_{\min }$, 
there exists no sequence $\{ C_{j}\} _{j\in \NN }$ 
of mutually distinct elements of $\hat{{\cal J}} $ such that 
$\min _{z\in C_{j}}d(z,J)\rightarrow 0$ as $j\rightarrow \infty .$  
\end{enumerate} 
\end{thm}
\begin{prop}
\label{countprop}
There exists an open set $V$ in {\em (Poly}$_{\deg \geq 2})^{3}$ such that for any $(h_{1},h_{2},h_{3})\in V$, 
$G=\langle h_{1},h_{2},h_{3}\rangle $ satisfies that 
$G\in {\cal G}_{dis}$,\ 
$\bigcup _{j=1}^{2}J(h_{j})\subset J_{\min }(G)$,
 $J_{\max }(G)=J(h_{3})$,   
$h_{j}^{-1}(J(h_{3}))\cap J_{\min }(G)\neq \emptyset $ for each 
$j=1,2$, 
and $\sharp (\hat{{\cal J}}_{G})=\aleph _{0}.$  
\end{prop}
\begin{prop}
\label{countcomp}
There exists a $3$-generator polynomial 
semigroup $G=\langle h_{1},h_{2},h_{3}\rangle $ 
in ${\cal G}_{dis }$ such that   
$\bigcup _{j=1}^{2}(h_{j})^{-1}(J_{\max }(G))\subset J_{\min }(G)$,
 $J_{\max }(G)=J(h_{3})$,  
$ \sharp (\hat{{\cal J}}_{G})=\aleph _{0} $,\ there exists 
a superattracting fixed point $z_{0}$ of some element 
of $G$ with  $z_{0}\in J(G)$,  and {\em int}$(J_{\min }(G))\neq \emptyset .$
\end{prop}
As mentioned before, these results illustrate new phenomena which can hold in the rational semigroups,  
but cannot hold in the dynamics of a single rational map or Kleinian groups. 

\begin{figure}[htbp]
\caption{The Julia set of a $3$-generator hyperbolic 
polynomial semigroup $G\in {\cal G}_{dis}$ such that  
$\sharp (\hat{{\cal J}}_{G})=\aleph _{0}.$}    
\ \ \ \ \ \ \ \ \ \ \ \ \ \ \ \ \ \ \ \ \ \ \ \ \ \ \ \ \ \ \ \ 
\ \ \ \  \ \ \ \ \ \ \ \ \ \ \ \ \ 
\includegraphics[width=3.3cm,width=3.3cm]{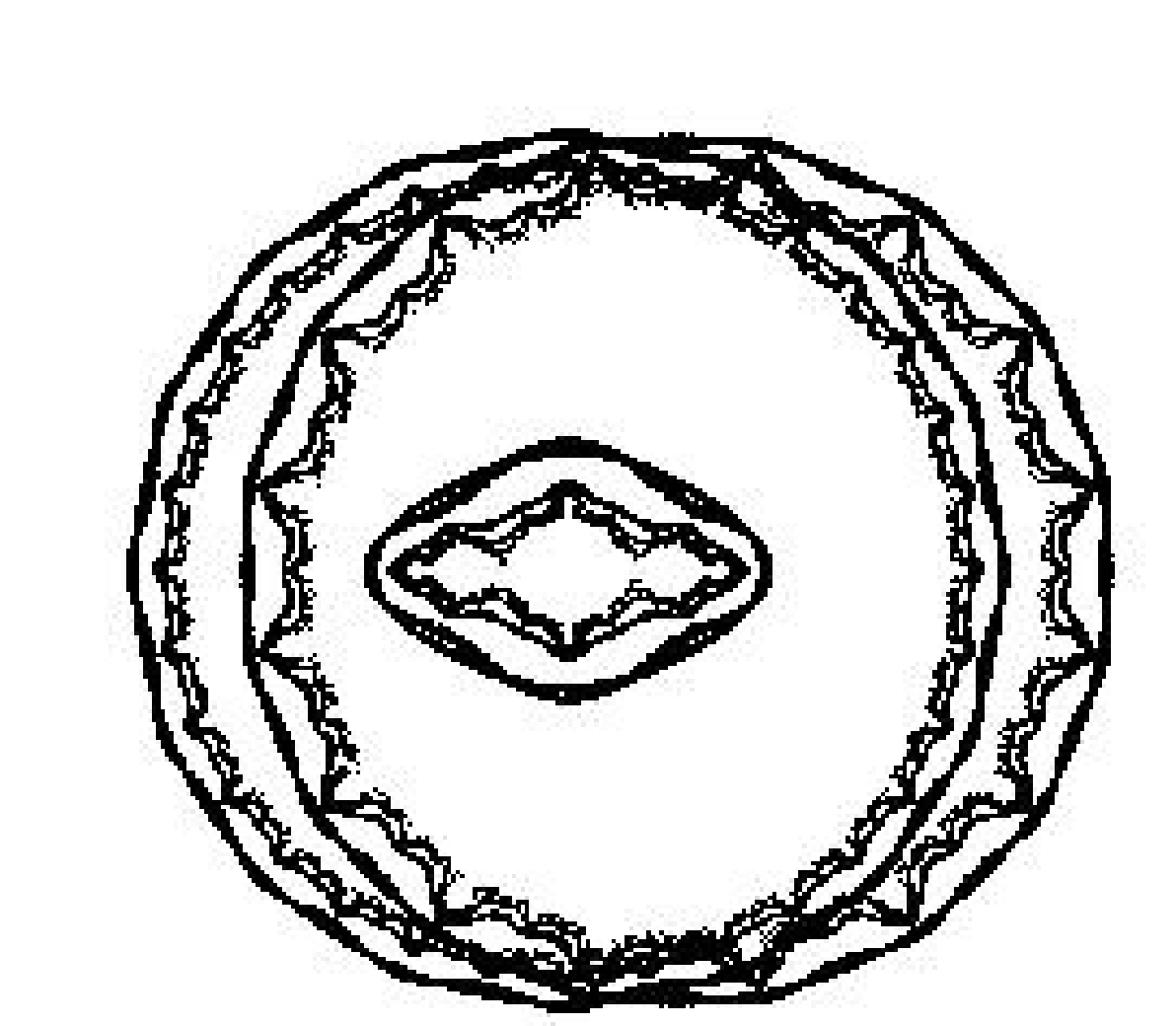}
\label{fig:3mapcountjulia2}
\end{figure}

 For the figure of the Julia set of a $3$-generator polynomial semigroup $G\in {\cal G}_{dis}$ such that  
 $\sharp \hat{{\cal J}}_{G}=\aleph _{0}$, see figure~\ref{fig:3mapcountjulia2}.
\begin{rem}
\label{r:intjf} 
In \cite{S7}, a new cohomology theory for (backward) self-similar systems (iterated function 
systems) was 
introduced by the author of this paper. By using it, 
for a finitely generated $G\in {\cal G}$, 
we can describe the space ${\cal J}_{G}$ of connected components of $J(G)$, and 
we can give 
some estimates on $\sharp  ({\cal J}_{G})$. Moreover, 
by using this new cohomology, a sufficient condition for the cardinality of 
the set of all connected components of the Fatou set of a postcritically 
bounded finitely generated polynomial semigroup $G$ to be infinity was given.

\end{rem}
\subsection{Hyperbolicity and semi-hyperbolicity}
\label{Hyperbolicity}
In this section, we present some results on hyperbolicity and 
semi-hyperbolicity. 
\begin{df}
Let $G$ 
be a 
polynomial semigroup generated by a 
subset $\G $ of Poly$_{\deg \geq 2}.$ 
Suppose 
$G\in {\cal G}_{dis}.$ Then 
we set 
$ \G_{\min }:=\{ h\in \G \mid 
J(h)\subset J_{\min }\} ,$
where $J_{\min }$ denotes the 
unique minimal element in $({\cal J},\ \leq )$ 
in Theorem~\ref{mainth2}-\ref{mainth2-3}. 
Furthermore, if $\G _{\min }\neq \emptyset $, 
let $G_{\min ,\G }$ be the subsemigroup 
of $G$ that is generated by 
$\G _{\min }.$ 
\end{df}
\begin{rem}
\label{jminrem}
Let $G$ be a polynomial semigroup generated by a compact subset 
$\G $ of Poly$_{\deg \geq 2}.$ Suppose $G\in {\cal G}_{dis}.$ Then, 
 by Proposition~\ref{bminprop},  
 we have $\G _{\min }\neq \emptyset $ and 
 $\G \setminus \G _{\min }\neq \emptyset .$ 
 Moreover, $\G _{\min }$ is a compact subset of $\G .$ For, 
 if $\{ h_{n}\} _{n\in \NN }\subset \G _{\min }$ and 
 $h_{n}\rightarrow h_{\infty }$ in $\G $, 
 then for a repelling periodic point $z_{0}\in 
 J(h_{\infty })$ of $h_{\infty }$, 
we have that $d(z_{0}, J(h_{n}))\rightarrow 0$ as $n\rightarrow \infty $, 
which implies that $z_{0}\in J_{\min }$ and thus $h_{\infty }\in \G _{\min }.$ 
\end{rem}
The following Proposition~\ref{nonminnoncpt} 
means that for a polynomial semigroup 
$G\in {\cal G}_{dis}$ generated by a 
compact subset $\G $ of Poly$_{\deg \geq 2}$, 
we rarely have the situation that 
``$\G \setminus \G _{\min }$ is not compact.''
\begin{prop}
\label{nonminnoncpt}
Let $G$ be a polynomial semigroup generated by a non-empty compact 
subset $\G $ of {\em Poly}$_{\deg \geq 2}.$ 
Suppose that $G\in {\cal G}_{dis}$ and that 
$\G \setminus \G _{\min }$ is not compact. 
Then, both of the following statements \ref{nonminnoncpt1} and  
\ref{nonminnoncpt2} hold.  
\begin{enumerate}
\item \label{nonminnoncpt1}
Let $h\in \G _{\min }.$ 
Then, $J(h)=J_{\min }(G), K(h)=\hat{K}(G)$, and  
{\em int}$(K(h))$ is a non-empty connected set.
\item \label{nonminnoncpt2}
Either 

\begin{enumerate}
\item \label{nonminnoncpt2-1}
for each $h\in \G _{\min }$, $h$ is hyperbolic and $J(h)$ is a 
quasicircle; or 
\item \label{nonminnoncpt2-2}
for each $h\in \G _{\min }$, {\em int}$(K(h))$ is an immediate 
parabolic basin of a parabolic fixed point of $h.$ 
\end{enumerate} 
\end{enumerate}
\end{prop}
\begin{df}
Let $G$ be a rational semigroup. 
\begin{enumerate}
\item 
We say that 
$G$ is hyperbolic if $P(G)\subset F(G).$ 
\item We say that $G$ is semi-hyperbolic if 
there exists a number $\delta >0$ and a 
number $N\in \NN $ such that 
for each $y\in J(G)$ and each $g\in G$, 
we have $\deg (g:V\rightarrow B(y,\delta ))\leq N$ for 
each connected component $V$ of $g^{-1}(B(y,\delta ))$, 
where $B(y,\delta )$ denotes the ball of radius $\delta $ 
with center $y$ with respect to the spherical distance, 
and $\deg (g:\cdot \rightarrow \cdot )$ denotes the 
degree of a finite branched covering. 
\end{enumerate}
\end{df}
\begin{rem}
There are many nice properties of hyperbolic or semi-hyperbolic 
rational semigroups. For example, for a finitely generated 
semi-hyperbolic rational semigroup $G$ , there exists an attractor in the Fatou set (\cite{S1, S4}),  
and the Hausdorff dimension $\dim _{H}(J(G))$ of the Julia set is less than or equal to the 
critical exponent $s(G)$ of the Poincar\'{e} series of $G$ (\cite{S4}).  
If we assume further the ``open set condition'', 
then $\dim _{H}(J(G))=s(G)$ (\cite{S6, SU3}). 
Moreover, if $G\in {\cal G}$ is generated by a compact set $\Gamma $ and if $G$ is semi-hyperbolic, 
then for each sequence $\gamma \in \Gamma ^{\NN }$, 
the basin of infinity for $\gamma $ is a John domain and the Julia set of $\gamma $ is locally connected 
(\cite{S4}). In \cite{S12}, by using the above result, we classify hyperbolic or semi-hyperbolic 
postcritically bounded compactly generated polynomial semigroups, in terms of the random complex dynamics. 
It is shown that in one of the classes, for almost every sequence $\gamma $, the Julia set $J_{\gamma }$ of $\gamma $ is a Jordan curve but not 
a quasicircle, the unbounded component of $\CCI \setminus J_{\gamma }$ is a John domain, and 
the bounded component of $\CC \setminus J_{\gamma }$ is not a John domain. 
Moreover, in \cite{S12,S11}, we find many examples with this phenomenon. 
Note that this phenomenon does not hold 
in the usual iteration dynamics of a single polynomial map $g$ with $\deg (g)\geq 2.$   
     
\end{rem}

We now present some results on semi-hyperbolic or hyperbolic polynomial semigroups 
in ${\cal G}_{dis}.$ These results are used to construct examples of 
semi-hyperbolic or hyperbolic polynomial semigroups $G\in {\cal G}_{dis}$ (see the proof of 
Proposition~\ref{Constprop}). Therefore these are important in terms of the sequel \cite{S11, S12}. 
\begin{thm}
\label{shshprop}
Let $G$ be a polynomial semigroup generated by a non-empty compact subset 
$\G $ of {\em Poly}$_{\deg \geq 2}.$ Suppose that $G\in {\cal G}_{dis}.$ 
If $G_{\min ,\G }$ is semi-hyperbolic, then $G$ is semi-hyperbolic.
\end{thm}
\begin{thm}
\label{hhprop}
Let $G$ be a polynomial semigroup generated by a 
non-empty compact subset $\G $ of 
{\em Poly}$_{\deg \geq 2}.$ Suppose that 
$G\in {\cal G}_{dis}.$ If 
$G_{\min ,\G }$ is hyperbolic and 
$(\bigcup _{h\in \G \setminus \G _{\min }}CV^{\ast }(h))\cap 
J_{\min }(G)=\emptyset $, 
then $G$ is hyperbolic.
\end{thm}
\begin{rem}
\label{hhrem}
In \cite{SS}, it will be shown that 
in Theorem~\ref{hhprop}, the condition 
$(\bigcup _{h\in \G \setminus \G _{\min }}CV^{\ast }(h))\cap 
J_{\min }(G)=\emptyset $ is necessary. For the figures of the Julia sets of 
hyperbolic polynomial semigroups $G\in {\cal G}_{dis}$, see figure~\ref{fig:dcgraph} and figure~\ref{fig:3mapcountjulia2}. 
\end{rem}
\begin{prop}
\label{p:2.312as}
Let $G$ be a polynomial semigroup generated by a non-empty compact 
subset $\G $ of {\em Poly}$_{\deg \geq 2}.$ 
Suppose that $G\in {\cal G}_{dis}$ and that 
$\G \setminus \G _{\min }$ is not compact. 
Suppose that statement~\ref{nonminnoncpt2-1} in Theorem~\ref{nonminnoncpt}
 holds. Then, both of the following statements hold. 
\begin{enumerate}
\item We have that $G_{\min ,\G }$ is hyperbolic and $G$ is semi-hyperbolic.
\item Suppose further that $(\bigcup _{h\in \G \setminus \G _{\min }}CV^{\ast }(h))\cap 
J_{\min }(G)=\emptyset $. Then $G$ is hyperbolic. 
\end{enumerate} 
\end{prop}

\subsection{Construction of examples}
\label{Const}
In this section, we present a way 
to construct examples of semigroups $G$ in ${\cal G}_{dis}$ (with some additional properties). 
These examples are important in terms of the sequel \cite{S11, S12}.  
\begin{prop}
\label{Constprop}
Let $G$ be a 
polynomial semigroup generated by 
a compact subset $\G $ of {\em Poly}$_{\deg \geq 2}.$ 
Suppose that $G\in {\cal G}$ and  
{\em int}$(\hat{K}(G))\neq \emptyset .$ 
Let $b\in $ {\em int}$(\hat{K}(G)).$ 
Moreover, let $d\in \NN $ be any positive integer such that 
$d\geq 2$, and such that 
$(d, \deg (h))\neq (2,2)$ for each $h\in \G .$ 
Then, there exists a number $c>0$ such that 
for each $a\in \CC $ with $0<|a|<c$, 
there exists a compact neighborhood $V$ of 
$g_{a}(z)=a(z-b)^{d}+b$ in {\em Poly}$_{\deg \geq 2}$ satisfying 
that for any non-empty subset $V'$ of $V$,  
the polynomial semigroup 
$H_{\G, V'} $ generated by the family $\G \cup V'$ 
belongs to ${\cal G}_{dis}$, $\hat{K}(H_{\G, V'})=\hat{K}(G)$ 
 and $(\G \cup V')_{\min }\subset \G .$ 
Moreover, in addition to the assumption above, 
if $G$ is semi-hyperbolic (resp. hyperbolic), 
then the above $H_{\G, V'}$ is semi-hyperbolic (resp. hyperbolic).   
\end{prop}
\begin{rem}
By Proposition~\ref{Constprop}, 
there exists a $2$-generator polynomial semigroup 
$G=\langle h_{1},h_{2}\rangle $ in ${\cal G}_{dis}$ such that 
$h_{1}$ has a Siegel disk.  
\end{rem}
\begin{df}
Let $d\in \NN $ with $d\geq 2.$ 
We set 
${\cal Y}_{d}:=\{ h\in \mbox{Poly} \mid \deg (h)=d\} $ 
endowed with the relative topology from Poly.
\end{df}
\begin{thm}
\label{shshfinprop}
Let $m\geq 2$ and let $d_{2},\ldots ,d_{m}\in \NN $ be such that  
$d_{j}\geq 2$ for each $j=2,\ldots ,m.$ Let 
$h_{1}\in {\cal Y}_{d_{1}}$ with {\em int}$(K(h_{1}))\neq \emptyset $ be  
such that $\langle h_{1}\rangle \in {\cal G}.$ 
Let $b_{2},b_{3},\ldots ,b_{m}\in $ {\em int}$(K(h_{1})).$ 
Then, both of the following statements hold. 
\begin{enumerate}
\item \label{shshfinprop1}
Suppose that $\langle h_{1}\rangle $ is 
semi-hyperbolic (resp. hyperbolic). 
Then, there exists a number $c>0$ such that 
for each $(a_{2},a_{3},\ldots ,a_{m})\in \CC ^{m-1}$ with 
$0<|a_{j}|<c$ ($j=2,\ldots ,m$), 
setting $h_{j}(z)=a_{j}(z-b_{j})^{d_{j}}+b_{j}$ ($j=2,\ldots ,m$), 
the polynomial semigroup 
$G=\langle h_{1},\ldots ,h_{m}\rangle $ satisfies that 
$G\in {\cal G}$, $\hat{K}(G)=K(h_{1})$ and $G$ is semi-hyperbolic (resp. hyperbolic). 
\item \label{shshfinprop2}
Suppose that $\langle h_{1}\rangle $ is 
semi-hyperbolic (resp. hyperbolic). Suppose also that 
either (i) there exists a $j\geq 2$ with $d_{j}\geq 3$, or 
(ii) $\deg(h_{1})=3$, $b_{2}=\cdots =b_{m}.$ Then, there exist 
$a_{2},a_{3},\ldots ,a_{m}>0$ such that setting 
$h_{j}(z)=a_{j}(z-b_{j})^{d_{j}}+b_{j}$ ($j=2,\ldots ,m$), 
the polynomial semigroup $G=\langle h_{1},h_{2},\ldots ,h_{m}\rangle $ 
satisfies that $G\in {\cal G}_{dis}$, $\hat{K}(G)=K(h_{1})$ and 
$G$ is semi-hyperbolic (resp. hyperbolic).  

\end{enumerate} 
\end{thm}
\begin{df}
Let $m\in \NN .$ We set 
\begin{itemize}
\item ${\cal H}_{m}:= 
\{ (h_{1},\ldots ,h_{m})\in 
(\mbox{Poly}_{\deg \geq 2 })^{m}\mid 
\langle h_{1},\ldots ,h_{m}\rangle  \mbox{ is hyperbolic}\}$, 
\item 
${\cal B}_{m}:= 
 \{ (h_{1},\ldots ,h_{m})\in 
(\mbox{Poly}_{\deg \geq 2 })^{m}\mid 
\langle h_{1},\ldots ,h_{m}\rangle \in {\cal G} \}$, and 
\item 
${\cal D}_{m}:= 
\{ (h_{1},\ldots ,h_{m})\in 
(\mbox{Poly}_{\deg \geq 2 })^{m}\mid 
J(\langle h_{1},\ldots ,h_{m}\rangle ) \mbox{ is disconnected}\} $.  
\end{itemize}
Moreover, let $\pi _{1}: (\mbox{Poly}_{\deg \geq 2 })^{m}\rightarrow 
\mbox{Poly}_{\deg \geq 2 }$ be the projection defined by 
$\pi (h_{1},\ldots ,h_{m})=h_{1}.$ 
\end{df}
\begin{thm}
\label{sphypopen}
Under the above notation, all of the following statements hold.
\begin{enumerate}
\item \label{sphypopen1}
${\cal H}_{m}, {\cal H}_{m}\cap {\cal B}_{m},$ 
${\cal H}_{m}\cap {\cal D}_{m}$, 
and ${\cal H}_{m}\cap {\cal B}_{m}\cap {\cal D}_{m}$ 
are open in $(\mbox{{\em Poly}}_{\deg \geq 2})^{m}.$ 
\item \label{sphypopen2}
Let $d_{1},\ldots ,d_{m}\in \NN $ be such that  
$d_{j}\geq 2$ for each $j=1,\ldots ,m.$\\  
Then, 
$\pi _{1}: {\cal H}_{m}\cap {\cal B}_{m}\cap 
({\cal Y}_{d_{1}}\times \cdots \times 
{\cal Y}_{d_{m}})\rightarrow 
{\cal H}_{1}\cap {\cal B}_{1}\cap {\cal Y}_{d_{1}}$ is 
surjective. 
\item \label{sphypopen3} 
Let $d_{1},\ldots ,d_{m}\in \NN $ be such that  
$d_{j}\geq 2$ for each $j=1,\ldots ,m$ and such that  
$(d_{1},\ldots ,d_{m})\neq (2,2,\ldots ,2).$ Then, 
 $\pi _{1}: {\cal H}_{m}\cap {\cal B}_{m}\cap {\cal D}_{m}\cap 
({\cal Y}_{d_{1}}\times \cdots \times 
{\cal Y}_{d_{m}})\rightarrow 
{\cal H}_{1}\cap {\cal B}_{1}\cap {\cal Y}_{d_{1}}$ is 
surjective. 
\end{enumerate}
\end{thm}
\begin{rem}
Combining Proposition~\ref{Constprop}, 
Theorem~\ref{shshfinprop}, and Theorem~\ref{sphypopen}, 
we can construct many examples of semigroups $G$ 
in ${\cal G} $ (or ${\cal G}_{dis}$) with some 
additional properties (semi-hyperbolicity, hyperbolicity, etc.). 
\end{rem}

\section{Tools}
\label{Tools}
To show the main results, we need 
some tools in this section.
\subsection{Fundamental properties of rational semigroups}
{\bf Notation:} 
For a rational semigroup $G$, we set 
$E(G):=\{ z\in \CCI \mid \sharp (\bigcup _{g\in G}g^{-1}(\{ z\} ))<\infty \} .$ 
This is called the exceptional set of $G.$ \\ 

The following Lemma~\ref{hmslem} and Theorem~\ref{repdense} will be used in the proofs of the main results. 
\begin{lem}[\cite{HM1,GR,S3,S1}]
\label{hmslem}
Let $G$ be a rational semigroup.\
\begin{enumerate}
\item 
\label{invariant}
For each $h\in G,\ $ we have 
$h(F(G))\subset F(G)$ and $h^{-1}(J(G))\subset J(G).$ Note that we do not 
have that the equality holds in general.
\item
\label{bss}
If $G=\langle h_{1},\ldots ,h_{m}\rangle $, then 
$J(G)=h_{1}^{-1}(J(G))\cup \cdots \cup h_{m}^{-1}(J(G)).$ 
More generally, if $G$ is generated by a compact subset 
$\G $ of {\em Rat}, then 
$J(G)=\bigcup _{h\in \G}h^{-1}(J(G)).$ 
(We call this property of the Julia set of a compactly generated rational semigroup ``backward self-similarity." )
\item
\label{Jperfect}
If \( \sharp (J(G)) \geq 3\) ,\ then \( J(G) \) is a 
perfect set.\ 
\item
\label{egset}
If $\sharp (J(G))\geq 3$ ,\ then 
$ \sharp (E(G)) \leq 2. $
\item
\label{o-set}
If a point \( z\) is not in \( E(G),\ \) then 
 \( J(G)\subset \overline{\bigcup _{g\in G}g^{-1}(\{ z\} )} .\) In particular 
if a point\ \( z\)  belongs to \ \( J(G)\setminus E(G), \) \ then
$ \overline{\bigcup _{g\in G}g^{-1}(\{ z\})}=J(G). $
\item
\label{backmin}
If \( \sharp (J(G))\geq 3 \) ,\ 
then
\( J(G) \) is the smallest closed backward invariant set containing at least
three points. Here we say that a set $A$ is backward invariant under $G$ if
for each \( g\in G,\ g^{-1}(A)\subset A.\ \)
\end{enumerate}
\end{lem}
\begin{thm}[\cite{HM1,GR,S3}]
\label{repdense}
Let $G$ be a rational semigroup.
If $\sharp (J(G))\geq 3$, then \\ 
$ J(G)=\overline{\{ z\in \CCI \mid 
\exists g \in G,\ g(z)=z,\ |m(g,z)|>1 \} } $, 
where $m(g,z)$ denotes the multiplier of $g$ at $z$ (\cite{Be1}).  
In particular,\ $J(G)=\overline{\bigcup _{g\in G}J(g)}.$
\end{thm}
\begin{rem}
If a rational semigroup $G$ contains an element $g$ with $\deg (g)\geq 2$, 
then $\sharp (J(g))\geq 3$, which implies that $\sharp (J(G))\geq 3.$  
\end{rem}
 \begin{lem}
\label{fibconncor2}
Let $G=\langle h_{1},h_{2}\rangle \in {\cal G}.$ Then, 
$h_{1}^{-1}(J(h_{2}))$ is connected.
\end{lem} 
\begin{proof}
Since $h_{2}\in G\in {\cal G}$, $F_{\infty }(h_{2})$ is simply connected. 
Since $G\in {\cal G}$, 
there exists no finite critical value of $h_{1}$ in $F_{\infty }(h_{2}).$ 
By the Riemann-Hurwitz formula, it follows that 
$h_{1}^{-1}(F_{\infty }(h_{2}))$ is connected and simply connected. 
Thus $\partial (h_{1}^{-1}(F_{\infty }(h_{2})))=h_{1}^{-1}(J(h_{2}))$ is connected.   
\end{proof}
\begin{df}
Let 
$G$ be a polynomial semigroup. 
Let $p\in \CC $ and 
$\epsilon >0.$ 
We set \\ 
${\cal F}_{G,p,\epsilon }:= 
\{\alpha :D(p,\epsilon )\rightarrow \CC \mid 
\alpha \mbox{ is a well-defined branch of }
g^{-1}, g\in G \} .$  
\end{df}
\begin{lem}
\label{invnormal2}
Let $\Gamma $ be a non-empty compact subset of {\em Poly}$_{\deg \geq 2}$ and 
let $G$ be a polynomial semigroup generated by $\G .$  
Let $R>0,\epsilon >0$, and \\ 
${\cal F}:=
\{ \alpha \circ \beta :D(0,1)\rightarrow \CC 
\mid 
\beta :D(0,1)\cong D(p,\epsilon ),\ 
\alpha :D(p,\epsilon )\rightarrow \CC ,\ 
\alpha \in {\cal F}_{G,p,\epsilon },\ p\in D(0,R)\} .$ 
Then, ${\cal F} $ is normal in $D(0,1).$ 
\end{lem}
\begin{proof}
Since $\Gamma $ is a non-empty compact subset of Poly$_{\deg \geq 2}$, 
there exists a ball $B$ around $\infty $ 
with $B\subset \CCI \setminus D(0,R+\epsilon )$ such that 
for each $h\in \Gamma $, $h(B)\subset B.$ 
Let $p\in D(0,R).$ Then, for each 
$\alpha \in {\cal F}_{G,p,\epsilon }$, 
$\alpha (D(p,\epsilon ))\subset \CCI \setminus B.$ 
Hence, ${\cal F} $ is normal in $D(0,1).$ 
\end{proof}

\subsection{A lemma from general topology}
\begin{lem}[\cite{N}]
\label{nadlem}
Let $X$ be a compact metric space 
and let $f:X\rightarrow X$ be a continuous 
open map. 
Let $A$ be a compact connected subset of 
$X.$ Then for each connected component 
$B$ of $f^{-1}(A)$, we have 
$f(B)=A.$ 

\end{lem}
\section{Proofs of the main results}
\label{Proofs} 
In this section, we demonstrate the main results.
\subsection{Proofs of results in \ref{concompsec}}
\label{pfconcompsec}
In this section, we demonstrate the results in \ref{concompsec}. \\ 

\noindent {\bf Proof of Theorem~\ref{mainth0}:}
First, we show the following:\\ 
Claim: For any $\lambda \in \Lambda $, 
$h_{\lambda }^{-1}(A)\subset A.$ 

 To show  the claim, 
let $\lambda \in \Lambda $ with $J(h_{\lambda })\neq \emptyset $
 and 
let $B$ be a connected component of 
$h_{\lambda }^{-1}(A).$ Then 
by Lemma~\ref{nadlem}, 
$h_{\lambda }(B)=A.$
Combining this with 
$h_{\lambda }^{-1}(J(h_{\lambda }))=
J(h_{\lambda })$, we obtain 
$B\cap J(h_{\lambda })\neq \emptyset .$ 
Hence $B\subset A.$ 
This means that 
$h_{\lambda }^{-1}(A)\subset A$ for each 
$\lambda \in \Lambda $ with $J(h_{\lambda })\neq \emptyset .$ 
Next, let $\lambda \in \Lambda $ with 
$J(h_{\lambda })=\emptyset .$ 
Then $h_{\lambda }$ is either identity or an 
elliptic M\"{o}bius transformation. 
By hypothesis and Lemma~\ref{hmslem}-\ref{invariant},  
 we obtain $h_{\lambda }^{-1}(A)\subset A.$ 
Hence, we have shown the claim.

 Combining the above claim with 
$\sharp A\geq 3$, 
by 
Lemma~\ref{hmslem}-\ref{backmin} 
we obtain $J(G)\subset A.$ Hence $J(G)=A$ and 
$J(G)$ is connected.
\qed  

\

\noindent {\bf Notation:}
We denote by $d$ the spherical distance on $\CCI .$ 
Given $A\subset \CCI $ and $z\in \CCI $, 
we set $d(z,A):=\inf \{ d(z,w)\mid w\in A\} .$ 
Given $A\subset \CCI $ and $\epsilon >0$,\ we set 
$B(A,\epsilon ):= \{ a\in \CCI \mid d(a,A)<\epsilon \} .$
Furthermore, given $A\subset \CC $, $z\in \CC $, 
and $\epsilon >0$, 
we set 
$d_{e}(z,A):=\inf \{ |z-w| \mid w\in A\} $ and 
$D(A,\epsilon ):= \{ a\in \CC \mid d_{e}(a,A)<\epsilon \} .$  

\ 

 We need the following lemmas to prove the main results.
\begin{lem}
\label{appjlem}
Let $G\in {\cal G}$ and let   
$J$ be a connected component of 
$J(G)$, $z_{0}\in J$ a point, and 
$\{ g_{n}\} _{n\in \NN }$ a sequence in  
$G$ such that 
$d(z_{0},J(g_{n}))\rightarrow 
0$ as $n\rightarrow \infty .$ 
Then 
$\sup\limits _{z\in J(g_{n})}d(z,J)\rightarrow 
0$ as $n\rightarrow \infty .$ 
\end{lem}
\begin{proof}
Suppose 
 there exists a connected component 
 $J'$ of $J(G)$ with $J'\neq J$ 
 and a subsequence $\{ g_{n_{j}}\} _{j\in \NN }$ 
 of $\{ g_{n}\} _{n\in \NN }$ such that 
 $\min\limits _{z\in J(g_{n_{j}})}
 d(z,J')\rightarrow 0$ as 
 $j\rightarrow \infty .$ 
 Since $J(g_{n_{j}})$ is compact and 
 connected for each $j$, 
 we may assume, passing to a subsequence, that there exists a non-empty 
 compact connected subset $K$ of $\CCI $ such that 
 $J(g_{n_{j}})\rightarrow K$ as 
 $j\rightarrow \infty $, with respect to 
 the Hausdorff metric. 
 Then $K\cap J \neq \emptyset $ and 
 $K\cap J'\neq \emptyset .$ Since 
 $K\subset J(G)$ and $K$ is connected, 
 it contradicts $J'\neq J.$
\end{proof}
\begin{lem}
\label{appj2lem}
Let $G\in {\cal G}.$ Then given $J\in {\cal J}$ and 
$\epsilon >0$,\ there exists an element $g\in G$ such that 
$J(g)\subset B(J,\epsilon ).$  

\end{lem}
\begin{proof}
We take a point $z\in J.$ Then, by Theorem~\ref{repdense}, 
there exists a sequence $\{ g_{n}\} _{n\in \NN }$ in $G$ such that 
$d(z,J(g_{n}))\rightarrow 0$ as $n\rightarrow \infty .$ 
By Lemma~\ref{appjlem}, we conclude that there exists 
an $n\in \NN $ such that 
$J(g_{n})\subset B(J,\epsilon ).$  
\end{proof}
\begin{lem}
\label{inftyj1}
Let $G$ be a polynomial semigroup. 
Suppose that $J(G)$ is disconnected, and  
 $\infty \in J(G).$ Then,\ 
 the connected component $A$ of 
 $J(G)$ containing $\infty $ is equal to 
 $\{ \infty \} .$ 
\end{lem}
\begin{proof}
 By Lemma~\ref{nadlem}, we obtain  
 $g^{-1}(A)\subset A$ for each 
 $g\in G.$ Hence, if $\sharp A\geq 3$, then 
 $J(G)\subset A$,\ by Lemma~\ref{hmslem}-\ref{backmin}. 
 Then $J(G)=A$ and it causes a contradiction, since 
 $J(G)$ is disconnected.  
\end{proof}

We now demonstrate Theorem~\ref{mainth1}.\\ 
\noindent {\bf Proof of Theorem~\ref{mainth1}:}
First, we show statement \ref{mainth1-1}.
Suppose the statement is false. Then, 
there exist elements $J_{1},J_{2}\in {\cal J}$ 
such that $J_{2}$ is included in the 
unbounded component $A_{1}$  
of $\CC \setminus J_{1}$, and such that 
$J_{1}$ is included in the unbounded 
component $A_{2}$  of 
$\CC \setminus J_{2}.$ 
Then we can find an $\epsilon >0$ such that 
$\overline{B(J_{2},\epsilon )}$ is 
included in the unbounded component 
of $\CC \setminus \overline{B(J_{1},\epsilon )}$, 
and such that 
$\overline{B(J_{1},\epsilon )}$ is 
included in the unbounded component of 
$\CC \setminus \overline{B(J_{2},\epsilon )}.$  
By Lemma~\ref{appj2lem}, for each $i=1,2$, 
 there exists an element $g_{i}\in G$ such that  
 $J(g_{i})\subset B(J_{i},\epsilon ).$ 
 This implies that 
 $J(g_{1})\subset A_{2}'$ and 
 $J(g_{2})\subset A_{1}'$, where 
 $A_{i}'$ denotes the unbounded 
 component of $\CC \setminus J(g_{i}).$  
 Hence we obtain 
 $K(g_{2})\subset A_{1}'.$ 
 Let  $v$ be a critical value 
 of $g_{2}$ in $\CC .$ Since 
 $P^{\ast }(G) $ is 
 bounded in $\CC $, we have 
 $v\in K(g_{2}).$  It implies  
 $v\in A_{1}'.$ 
  Hence $g_{1}^{l}(v)\rightarrow \infty $ as $l\rightarrow \infty .$ 
  However, this implies a 
 contradiction since 
 $P^{\ast }(G) $ is 
 bounded in $\CC .$ 
 Hence we have shown statement \ref{mainth1-1}.

 Next, we show statement \ref{mainth1-2}. 
 Let $F_{1}$ be a connected component of 
 $F(G).$ 
 Suppose that there exist three connected 
 components $J_{1},J_{2}$ and $J_{3}$ of 
 $J(G)$ such that they are mutually disjoint and 
 such that 
 $\partial F_{1}\cap J_{i}\neq \emptyset $ 
 for each $i=1,2,3.$ 
 Then, by statement \ref{mainth1-1}, 
 we may assume that we have 
 either (1): $J_{i}\in {\cal J}$ for each 
  $i=1,2,3$ and $J_{1}<J_{2}<J_{3}$, or 
  (2): $J_{1},J_{2}\in {\cal J},\ J_{1}<J_{2}$, and 
  $\infty \in J_{3}.$ 
  Each of these cases implies that 
  $J_{1}$ is included in a bounded component 
  of $\CC \setminus J_{2}$ and 
  $J_{3}$ is included in the unbounded component of 
  $\CCI \setminus J_{2}.$ However, it 
  causes a contradiction, 
  since $\partial F_{1}\cap J_{i}\neq \emptyset $ 
 for each $i=1,2,3.$ Hence, we have shown that 
 we have either \\  
 Case I: $\sharp \{ J:\mbox{component of }J(G)\mid 
 \partial F_{1}\cap J\neq \emptyset \} =1$ or \\ 
 Case II: $\sharp \{ J:\mbox{component of }J(G)\mid 
 \partial F_{1}\cap J\neq \emptyset \} =2.$

  Suppose that we have Case I. Let 
 $J_{1}$ be the connected component of $J(G)$ such that 
 $\partial F_{1}\subset J_{1}.$ 
 Let $D_{1}$ be the connected component of 
 $\CCI \setminus J_{1}$ containing $F_{1}.$ 
 Since $\partial F_{1}\subset J_{1},$
  we have $\partial F_{1}\cap D_{1}=\emptyset .$ 
  Hence, we have $F_{1}=D_{1}.$ Therefore, 
  $F_{1}$ is simply connected. 

  Suppose that we have Case II. 
 Let $J_{1}$ and $J_{2}$ be the two connected components 
 of $J(G)$ such that 
 $J_{1}\neq J_{2}$ and 
 $\partial F_{1}\subset J_{1}\cup J_{2}.$ 
 Let $D$ be the connected component of 
 $\CCI \setminus ( J_{1}\cup J_{2})$ 
 containing $F_{1}.$ 
 Since $\partial F_{1}\subset J_{1}\cup J_{2},$
 we have $\partial F_{1}\cap D=\emptyset .$ Hence, 
 we have $F_{1}=D.$ Therefore, $F_{1}$ is 
 doubly connected. Thus, we have shown 
statement \ref{mainth1-2}.

 We now show statement \ref{mainth1-3}.
 Let $g\in G$ be an element and 
$J$ a connected component 
of $J(G).$ 
 Suppose that $g^{-1}(J)$ is disconnected. 
 Then, by Lemma~\ref{nadlem}, 
 there exist at most finitely 
 many connected components $C_{1},\ldots ,C_{r}$
 of $g^{-1}(J)$ with $r\geq 2.$  
 Then there exists a positive number 
 $\epsilon $ such that 
 denoting by $B_{j}$ the connected 
 component of $g^{-1}(B(J,\epsilon ))$ 
 containing $C_{j}$ for each 
 $j=1,\ldots ,r$, 
 $\{ B_{j}\} $ are mutually 
 disjoint. By 
 Lemma~\ref{nadlem}, 
 we see that, for each 
 connected component $B$ of 
 $g^{-1}(B(J,\epsilon )), $
  $g(B)=B(J,\epsilon )$ 
 and $B\cap C_{j}\neq \emptyset $ for 
 some $j.$ Hence we get that 
 $g^{-1}(B(J,\epsilon ))=
 \bigcup _{j=1}^{r}B_{j}$ 
  (disjoint union) and 
  $g(B_{j})=B(J,\epsilon )$ for each 
  $j.$  
By Lemma~\ref{appj2lem}, there exists an element 
$h\in G$ such that $J(h)\subset B(J,\epsilon ).$  
Then it follows that 
$g^{-1}(J(h))\cap B_{j}\neq \emptyset 
$ for each $j=1,\ldots ,r.$ Moreover, 
we have $g^{-1}(J(h))\subset 
g^{-1}(B(J,\epsilon ))=\bigcup _{j=1}^{r}B_{j}.$
On the other hand, by 
Lemma~\ref{fibconncor2}, 
we have that $g^{-1}(J(h))$ is connected. 
This is a contradiction. 
Hence, we have shown that, for each $g\in G$ and 
each connected component $J$ of $J(G)$,\   
$g^{-1}(J)$ is connected. 

 By Lemma~\ref{inftyj1}, we get that 
if $J\in {\cal J}$, then $g^{\ast }(J)\in {\cal J}.$ 
 Let $J_{1}$ and $J_{2}$ be two elements 
 of ${\cal J}$ such that 
 $J_{1}\leq J_{2}.$ 
 Let $U_{i}$ be the unbounded component 
 of $\CC \setminus J_{i}$, for each 
 $i=1,2.$ 
 Then 
 $U_{2}\subset U_{1}.$ 
  Let $g\in G$ be an element. 
 Then $g^{-1}(U_{2})\subset g^{-1}(U_{1}).$ 
 Since $g^{-1}(U_{i})$ is the unbounded 
 connected component of $\CC \setminus 
 g^{-1}(J_{i})$ for each $i=1,2$, 
 it follows that 
 $g^{-1}(J_{1})\leq g^{-1}(J_{2}).$ 
 Hence $g^{\ast }(J_{1})\leq g^{\ast }(J_{2})$, 
 otherwise 
$g^{\ast }(J_{2})<g^{\ast }(J_{1})$, 
and 
 it contradicts $g^{-1}(J_{1})\leq g^{-1}(J_{2}).$    
\qed \\ 
\subsection{Proofs of results in \ref{Upper}}
\label{Proof of Upper}
In this section, we prove the results in Section \ref{Upper}, 
except Theorem~\ref{polyandrathm1}-\ref{polyandrathm1-2} and 
Theorem~\ref{polyandrathm1}-\ref{polyandrathm1-3}, which will be 
proved in Section~\ref{Proof of Properties}.

 To demonstrate Theorem~\ref{polyandrathm1}, we need the 
following lemmas.
\begin{lem}
\label{j1aj2lem}
Let $G$ be a polynomial semigroup in ${\cal G}_{dis}.$ 
Let $J_{1},J_{2}\in \hat{{\cal J}}$ be two elements with 
$J_{1}\neq J_{2}.$ 
Then, we have the following.
\begin{enumerate}
\item \label{j1aj2lem1}
If $J_{1},J_{2}\in {\cal J}$ and 
$J_{1}<J_{2}$, then 
there exists a doubly connected component $A$ of $F(G)$ such that 
$J_{1}<A<J_{2}.$ 
\item \label{j1aj2lem2}
If $\infty \in J_{2}$, then 
there exists a doubly connected component $A$ of $F(G)$ such that 
$J_{1}<A.$ 
\end{enumerate}
\end{lem}
\begin{proof}
First, we show statement \ref{j1aj2lem1}. 
Suppose that $J_{1},J_{2}\in {\cal J}$ and $J_{1}<J_{2}.$ 
We set $B=\bigcup _{J\in {\cal J}, J_{1}\leq J\leq J_{2}}J.$ 
Then, $B$ is a closed disconnected set. 
Hence, there exists a multiply connected component 
$A'$ of $\CCI \setminus B.$ Since $A'$ is multiply connected, 
we have that 
$A'$ is included in the unbounded component of $\CCI \setminus J_{1}$,  
and that 
$A'$ is included in a bounded component of $
\CCI \setminus J_{2}.$ 
This 
implies that $A'\cap J(G)=\emptyset .$ 
Let $A$ be the connected component of $F(G)$ such that  
$A'\subset A.$ 
Since $B\subset J(G)$, we have 
$F(G)\subset \CCI \setminus B.$  
Hence, $A$ must be equal to $A'.$ 
Since $A'$ is multiply connected, 
Theorem~\ref{mainth1}-\ref{mainth1-2} implies 
that $A=A'$ is doubly connected. 
Let $J$ be the connected component $J(G)$ such that 
$J<A$ and $J\cap \partial A\neq \emptyset .$  Then, since 
$A'=A$ is included in the unbounded component of 
$\CCI \setminus J_{1},$ 
we have that $J$ does not meet any bounded component 
of $\CC \setminus J_{1}.$ Hence, we obtain $J_{1}\leq J$, 
which implies that $J_{1}\leq J<A.$ 
Therefore, 
$A$ is a doubly connected component of $F(G)$ such that 
$J_{1}<A<J_{2}.$ Hence, we have shown statement \ref{j1aj2lem1}.

 Next, we show statement \ref{j1aj2lem2}. 
 Suppose that $\infty \in J_{2}.$ 
We set $B=(\bigcup _{J\in {\cal J}, J_{1}\leq J}J )\cup J_{2}.$ 
Then, $B$ is a disconnected closed set. 
Hence, there exists a multiply connected component 
$A'$ of $\CCI \setminus B.$ By the same method as that of proof 
of statement \ref{j1aj2lem1}, we see that 
$A'$ is equal to a doubly connected component $A$ of 
$F(G)$ such that $J_{1}<A.$ 
 Hence, we have shown statement \ref{j1aj2lem2}. 
\end{proof}
\begin{lem}
\label{cptmsetlem}
Let $H_{0}$ be a real affine semigroup 
generated by a compact set $C$ in {\em RA}. 
Suppose that each element $h\in C$ is of the form 
$h(x)=b_{1}(h)x+b_{2}(h)$, where 
$b_{1}(h), b_{2}(h)\in \RR $, 
$|b_{1}(h)|>1.$   
Then, for any subsemigroup $H$ of $H_{0}$, we have 
$M(H)=J(\eta (H))\subset \RR .$ 

\end{lem}
\begin{proof}
From the assumption, there exists a number $R>0$ such that 
for each $h\in C$, 
$\eta (h)(B(\infty ,R))\subset B(\infty ,R).$  
Hence, we have $B(\infty ,R)\subset F(\eta (H))$, which 
implies that $J(\eta (H))$ is a bounded subset of $\CC .$  
We consider the following cases:\\ 
Case 1: $\sharp (J(\eta (H)))\geq 3.$ \\ 
Case 2: $\sharp (J(\eta (H)))\leq 2.$ 

 Suppose that we have case 1. Then, 
from Theorem~\ref{repdense}, it follows that 
$M(H)=J(\eta (H))\subset \RR .$ 

 Suppose that we have case 2. 
Let $b(h)$ be the unique fixed point of $h\in H$ 
in $\RR .$ From the hypothesis, 
we have that for each $h\in H$, 
$b(h)\in J(\eta (H)).$ Since we assume $\sharp (J(\eta (H)))\leq 2$, 
Lemma~\ref{hmslem}-\ref{invariant} implies that 
 there exists a point $b\in \RR $ such that 
for each $h\in H$, we have $b(h)=b.$ Then 
any element $h\in H$ is of the form 
$h(x)=c_{1}(h)(x-b)+c_{2}(h),$ where $c_{1}(h), c_{2}(h)\in \RR , 
|c_{1}(h)|>1.$ Hence, $M(H)=\{ b\} \subset  J(\eta (H)).$ 
Suppose that there exists a point $c$ in $J(\eta (H))\setminus 
\{ b\} .$ Since $J(\eta (H))$ is a bounded set of $\CC $, 
and since we have $h^{-1}(J(\eta (H)))\subset J(\eta (H))$ for each $h\in H$ 
(Lemma~\ref{hmslem}-\ref{invariant}), 
we get that $h^{-1}(c)\in J(\eta (H))\setminus (\{ b\} \cup \{ c\} ) $, for each element $h\in H.$ 
This implies that $\sharp (J(\eta (H)))\geq 3$, which is a contradiction.
Hence, we must have that 
$J(\eta (H))=\{ b\} =M(H).$     
\end{proof}
We need the notion of Green's functions, in order to demonstrate 
Theorem~\ref{polyandrathm1}.
\begin{df}
Let $D$ be a domain in $\CCI $ with 
$\infty \in D.$ 
We denote by $\varphi (D,z)$ Green's function on 
$D$ with pole at $\infty .$ 
By definition, this is the unique function on $D\cap \CC $ with the 
properties:
\begin{enumerate}
\item $\varphi (D,z)$ is harmonic and positive in $D\cap \CC $;
\item $\varphi (D,z)-\log |z|$ is bounded in a neighborhood of 
$\infty $; and 
\item there exists a Borel subset $A$ of $\partial D$ such that 
the logarithmic capacity of $(\partial D)\setminus A$ is zero and such that 
for each $\zeta \in A$, we have $\varphi (D,z)\rightarrow 0$ as $z\rightarrow \zeta.$ 
\end{enumerate} 
\end{df}
\begin{rem}
\label{r:Greenf}
\ 
\begin{enumerate}
\item The limit $\lim\limits _{z\rightarrow \infty }(\varphi (D,z)-\log |z|)$ 
exists and this is called {\em Robin's constant} of $D.$ 
\item If $D$ is a simply connected domain with $\infty \in D$ and $\sharp (\CCI \setminus D)>1$,  
then we have $\varphi (D,z)=-\log |\psi (z)|$, where 
$\psi :D\rightarrow \{ z\in \CC \mid |z|<1\} $ denotes 
a biholomorphic map with $\psi (\infty )=0.$ 
\item \label{r:Greenf3}
It is well-known that for any $g\in $ Poly$_{\deg \geq 2},$  
\begin{equation}
\label{polygreeneq}
\varphi (F_{\infty }(g),z)=\log |z|+
\frac{1}{\deg (g)-1}\log |a(g)|+o(1)\ \mbox{ as }z\rightarrow \infty .
\end{equation}
(See \cite[p147]{Ste}.) Note that the point 
$-\frac{1}{\deg (g)-1}\log |a(g)|\in \RR $ is the unique fixed point 
of $\Psi (g)$ in $\RR .$  
\end{enumerate}
\end{rem}
\begin{lem}
\label{greenk1k2}
Let $K_{1}$ and $K_{2}$ be two non-empty connected compact sets in 
$\CC $ such that $K_{1}<K_{2}$ and $\sharp K_{1}\neq 1.$  
Let $A_{i}$ denote the unbounded component of 
$\CCI \setminus K_{i}$, for each $i=1,2.$ 
Then, we have 
$\lim _{z\rightarrow \infty }(\log |z|-\varphi (A_{1},z)) 
<\lim _{z\rightarrow \infty }(\log |z|-\varphi (A_{2},z)).$
\end{lem}
\begin{proof}
The function 
$\phi (z):= 
\varphi (A_{2},z)-\varphi (A_{1},z)=
(\log |z|-\varphi (A_{1},z))-(\log |z|-\varphi (A_{2},z))$ 
is harmonic on $A_{2}\cap \CC .$ This $\phi $ is bounded 
around $\infty .$ Hence $\phi $ extends to a 
harmonic function on $A_{2}.$ 
Moreover, since $K_{1}<K_{2}$, 
we have $\limsup _{z\rightarrow \partial A_{2}}\phi (z)$ $<0.$ 
From the maximum principle, 
it follows that 
$\phi (\infty )<0.$ Therefore, the statement of our lemma holds. 
\end{proof}
In order to demonstrate Theorem~\ref{polyandrathm1}-\ref{polyandrathm1-1}, 
we will prove the following lemma.  
 (Theorem~\ref{polyandrathm1}-\ref{polyandrathm1-2} and 
 Theorem~\ref{polyandrathm1}-\ref{polyandrathm1-3} will be proved 
 in Section~\ref{Proof of Properties}.) 
\begin{lem}
\label{part1-1lem}
Let $G$ be a polynomial semigroup in ${\cal G}.$ Then, there exists an 
injective map $\tilde{\Psi }:
\hat{{\cal J}}_{G}\rightarrow 
{\cal M}_{\Psi (G)}$ such that: 
\begin{enumerate}
\item  
if $J_{1},J_{2}\in {\cal J}_{G}$ and 
$J_{1}<J_{2}$, 
then $\tilde{\Psi }(J_{1})<_{r}\tilde{\Psi }(J_{2})$;  
\item if $J\in \hat{{\cal J}}_{G}$ and $\infty \in J$, then 
$+\infty \in \tilde{\Psi }(J)$; and 
\item if $J\in {\cal J}_{G}$, then 
$\tilde{\Psi }(J)\subset \hat{\RR }\setminus \{ +\infty \} .$   
\end{enumerate}
\end{lem}  
\begin{proof}
We first show the following claim.\\  
\noindent Claim 1: 
In addition to the assumption of Lemma~\ref{part1-1lem}, 
if we have $\infty \in F(G)$, then 
$M(\Psi (G))\subset \hat{\RR } \setminus \{ +\infty \} .$ 
 
 To show this claim, let $R>0$ be a number such that 
 $J(G)\subset D(0,R).$ Then, 
 for any $g\in G$, we have 
 $K(g)<\partial D(0,R).$  
By Lemma~\ref{greenk1k2}, we get that 
there exists a constant $C>0$ such that for each 
$g\in G$, $\frac{-1}{\deg (g)-1}\log |a(g)|\leq C.$ 
Hence, it follows that $M(\Psi (G))\subset 
[-\infty ,C].$ Therefore, we have shown Claim 1.

We now prove the statement of the lemma in the case $G\in {\cal G}_{con} .$ 
If $\infty \in F(G)$, then claim 1 implies that 
$ M({\Psi (G)}) \subset \hat{\RR }\setminus \{ +\infty \} $ and the statement of the lemma holds. 
We now suppose $\infty \in J(G).$  
We put $L_{g}:= \max _{z\in J(g) } |z|$ for each $g\in G.$  
Moreover, for each non-empty compact subset $E$ of $\CC $, we denote by Cap $(E)$ the logarithmic capacity of $E.$ 
We remark that Cap$(E)=\exp (\lim _{z\rightarrow \infty }(\log |z|-\varphi (D_{E},z)))$, where $D_{E}$ denotes the 
connected component of $\CCI \setminus E$ containing $\infty .$   
We may assume that $0\in P^{\ast }(G).$ Then, by \cite{A}, we have 
Cap$(J(g))\geq $ Cap $([0,L_{g}])\geq L_{g}/4$ for each $g\in G.$  
Combining this with $\infty \in J(G)$, Theorem~\ref{repdense}, 
and Remark~\ref{r:Greenf}-\ref{r:Greenf3},  
we obtain $+\infty \in M_{\Psi (G)}$ and defining 
$\tilde{\Psi }(J(G))$ to be the connected component of ${\cal M}_{\Psi (G)}$ containing $+\infty $, 
the statement of the lemma holds. 
  
We now prove the statement of the lemma in the case $G\in {\cal G}_{dis}.$  
Let $\{ J_{\lambda }\} _{\lambda \in \Lambda }$
be the set $\hat{{\cal J}}_{G}$ of all connected components 
of $J(G).$ 
By Lemma~\ref{appj2lem}, 
for each $\lambda \in \Lambda $ and each $n\in \NN $, 
there exists an element $g_{\lambda ,n}\in G$ such that 
\begin{equation}
\label{polyandrathm1pf0a}
J(g_{\lambda ,n})\subset B(J_{\lambda },\frac{1}{n}).
\end{equation} 
We have that the fixed point of 
$\Psi (g_{\lambda ,n})$ in $\RR $ is 
equal to $\frac{-1}{\deg (g_{\lambda ,n})-1}\log |a(g_{\lambda ,n})|.$ 
We may assume that 
$\frac{-1}{\deg (g_{\lambda ,n})-1}\log |a(g_{\lambda ,n})|\rightarrow 
\alpha _{\lambda }$ as $n\rightarrow \infty $, 
where $\alpha _{\lambda }$ is an element of $\hat{\RR }.$ 
For each $\lambda \in \Lambda $, let $B_{\lambda }\in 
{\cal M}_{\Psi (G)}$ be the element with 
$\alpha _{\lambda }\in B_{\lambda }.$ 
Let $\tilde{\Psi }(J_{\lambda })=B_{\lambda }$ for each $\lambda \in \Lambda .$ 
We will show the following claim.\\ 
Claim 2:  
If $\lambda ,\xi $ are two elements in $\Lambda $ with 
$\lambda \neq \xi $, then 
$B_{\lambda }\neq B_{\xi }.$ 
Moreover, if $J_{\lambda }, J_{\xi }\in {\cal J}_{G}$ and 
$J_{\lambda }<J_{\xi }$, then $B_{\lambda }<_{r}B_{\xi }. $
Furthermore, if $J_{\xi }\in \hat{{\cal J}}_{G}$ with 
$\infty \in J_{\xi }$, then $+\infty \in B_{\xi }.$ 
 
 To show this claim,  let $\lambda $ and $\xi $ be two elements in 
 $\Lambda $ with $\lambda \neq \xi .$ 
 We have the following two cases:\\ 
Case 1: $J_{\lambda },J_{\xi }\in {\cal J}_{G}$ and 
$J_{\lambda }<J_{\xi }.$ \\ 
Case 2:  $J_{\lambda }\in {\cal J}_{G}$ and $\infty \in J_{\xi } .$ 
(Note: in this case, by Lemma~\ref{inftyj1}, we have $J_{\xi }=\{ \infty \} .$) 

 Suppose that we have case 1. 
 By Lemma~\ref{j1aj2lem}, 
 there exists a doubly connected component $A$ of $F(G)$ such that 
\begin{equation}
 \label{polyandrathm1pf0b}
 J_{\lambda }<A<J_{\xi }.
\end{equation} 
 Let $\zeta _{1}$ and $\zeta _{2}$ be two Jordan curves in 
 $A$ such that they are not null-homotopic in $A$, and such that 
 $\zeta _{1}<\zeta _{2}.$   
For each $i=1,2,$ let $A_{i}$ be the unbounded component 
of 
$\CCI \setminus \zeta _{i}.$ 
Moreover, we set 
$\beta _{i}:=\lim _{z\rightarrow \infty }
(\log |z|-\varphi (A_{i},z))$, for each $i=1,2.$ 
By Lemma~\ref{greenk1k2}, we have $\beta _{1}<\beta _{2}.$ 
Let $g\in G$ be any element. 
By (\ref{polyandrathm1pf0a}) and (\ref{polyandrathm1pf0b}), 
there exists an $m\in \NN $ such that 
$J(g_{\lambda ,m})<\zeta _{1}.$ 
Since $P^{\ast }(G)\subset K(g_{\lambda ,m})$, 
it follows that $P^{\ast }(G)$ is included in the bounded component 
of $\CC \setminus \zeta _{1}.$ 
Hence, 
we see that  
\begin{equation}
\label{polyandrathm1pf0}
\mbox{either }J(g)<\zeta _{1}, \mbox{ or }
\zeta _{2}<J(g). 
\end{equation}
From Lemma~\ref{greenk1k2}, it follows that 
either 
$\frac{-1}{\deg (g)-1}\log |a(g)|<\beta _{1}$, or 
$\beta _{2}<\frac{-1}{\deg (g)-1}\log |a(g)|.$ 
This implies that 
\begin{equation}
\label{polyandrathm1pf1}
M(\Psi (G))\subset \hat{\RR }\setminus (\beta _{1},\beta _{2}),
\end{equation}
where $(\beta _{1},\beta _{2}):=\{ x\in \RR \mid \beta _{1}<x<\beta _{2}\}.$ 
Moreover, combining (\ref{polyandrathm1pf0a}), (\ref{polyandrathm1pf0b}), 
and (\ref{polyandrathm1pf0}), we get that  
there exists a number $n_{0}\in \NN $ such that 
for each $n\geq n_{0}$, 
$J(g_{\lambda ,n})<\zeta _{1}<\zeta _{2}<J(g_{\xi ,n}).$ 
From Lemma~\ref{greenk1k2}, it follows that 
\begin{equation}
\label{polyandrathm1pf2}
\frac{-1}{\deg (g_{\lambda ,n})-1}\log |a(g_{\lambda ,n})|
<\beta _{1}<\beta _{2}<
\frac{-1}{\deg (g_{\xi ,n})-1}\log |a(g_{\xi ,n})|,
\end{equation} 
for each $n\geq n_{0}.$ 
By (\ref{polyandrathm1pf1}) and (\ref{polyandrathm1pf2}), 
we obtain $B_{\lambda }<_{r} B_{\xi }.$ 

 We now suppose that we have case 2.
Then, by Lemma~\ref{j1aj2lem}, 
there exists a doubly connected component $A$ of $F(G)$ such that 
$J_{\lambda }<A.$ Continuing the same argument as that of case 1, 
we obtain $B_{\lambda }\neq B_{\xi }.$ 
In order to show $+\infty \in B_{\xi }$, 
let $R$ be any number such that $P^{\ast }(G)\subset D(0,R).$ 
Since $P^{\ast }(G)\subset K(g)$ for each $g\in G$, 
combining it with (\ref{polyandrathm1pf0a}) and 
Lemma~\ref{inftyj1}, we see that there exists an $n_{0}=n_{0}(R)$ 
such that for each $n\geq n_{0}$, $D(0,R)<J(g_{\xi ,n}).$ 
From Lemma~\ref{greenk1k2}, it follows that 
$\frac{-1}{\deg (g_{\xi ,n})-1}\log |a(g_{\xi ,n})|
\rightarrow +\infty .$ Hence, $+\infty \in B_{\xi }.$ 
Therefore, we have shown 
Claim 2. 
%

 Combining Claims 1 and 2, the statement of the lemma follows. 
 
 Therefore, we have proved Lemma~\ref{part1-1lem}. 
\end{proof}

\ 

We now demonstrate Theorem~\ref{polyandrathm1}-\ref{polyandrathm1-1}.\\ 
{\bf Proof of Theorem~\ref{polyandrathm1}-\ref{polyandrathm1-1}:} 
From Lemma~\ref{part1-1lem}, Theorem~\ref{polyandrathm1}-\ref{polyandrathm1-1}
follows.
\qed 

\ 

We now demonstrate Corollary~\ref{polyandracor}.\\ 
{\bf Proof of Corollary~\ref{polyandracor}:}
By Theorem~\ref{repdense}, 
we have 
$J(\Theta (G))=
\overline{ \bigcup _{h\in \Theta (G)}J(h) }=
\overline{\bigcup _{g\in G}J(\Theta (g))}$, 
where  the closure is taken in $\CCI .$  
Since $J(\Theta (g))=
\{ z\in \CC \mid |z|=|a(g)|^{-\frac{1}{\deg(g)-1}}\} $,  
we obtain 
\begin{equation}
\label{polyandracorpf1}
J(\Theta (G))=
\overline{\bigcup _{g\in G}\{ z\in \CC 
\mid |z|=|a(g)|^{-\frac{1}{\deg(g)-1}}\} },
\end{equation}
where the closure is taken in $\CCI .$ 
Hence, we see that 
$\sharp ({\hat{\cal J}}_{\Theta (G)})$ 
is equal to the cardinality of 
the set of all connected components of 
$J(\Theta (G))\cap [0,+\infty ].$
Moreover, let $\psi :[0,+\infty ]\rightarrow \hat{\RR }$ 
be the homeomorphism defined by 
$\psi (x):=\log (x)$ for $x\in (0,+\infty )$, 
$\psi (0):=-\infty $, and $\psi (+\infty )=+\infty .$ 
Then, 
(\ref{polyandracorpf1}) implies that, 
the map $\psi : 
[0,\infty ]\rightarrow \hat{\RR }$, 
maps $J(\Theta (G))\cap [0,+\infty ]$
onto $M(\Psi (\Theta (G))).$ 
For any $J\in \hat{{\cal J}}_{\Theta (G)}$, 
let $\tilde{\psi }(J)\in {\cal M}_{\Psi (\Theta (G))}=
{\cal M}_{\Psi (G)}$ be the element such that 
$\psi (J\cap [0,+\infty ])=\tilde{\psi }(J).$ 
Then,  
the map $\tilde{\psi }:\hat{{\cal J}}_{\Theta (G)}\rightarrow 
{\cal M}_{\Psi (\Theta (G))}={\cal M}_{\Psi (G)}$ is a 
bijection, and moreover, 
for any $J_{1},J_{2}\in {\cal J}_{\Theta (G)}$, 
we have that $J_{1}<J_{2}$ if and only if 
$\tilde{\psi }(J_{1})<_{r}\tilde{\psi }(J_{2}).$ 
Furthermore, for any $J\in \hat{{\cal J}}_{\Theta (G)}$, 
$\infty \in J$ if and only if $+\infty \in \tilde{\psi}(J).$ 
Let $\tilde{\Theta }:\hat{{\cal J}}_{G}\rightarrow 
\hat{{\cal J}}_{\Theta (G)}$ be the map defined by 
$\tilde{\Theta }=\tilde{\psi }^{-1}\circ \tilde{\Psi }$, 
where $\tilde{\Psi }: \hat{{\cal J}}_{G}\rightarrow {\cal M}_{\Psi (G)}$ is the 
map in Lemma~\ref{part1-1lem}.  
Then, by Lemma~\ref{part1-1lem}, 
$\tilde{\Theta }:\hat{{\cal J}}_{G}\rightarrow 
\hat{{\cal J}}_{\Theta (G)}$ is injective, and moreover, 
if $J_{1},J_{2}\in {\cal J}_{G}$ and $J_{1}<J_{2}$, then 
$\tilde{\Theta }(J_{1})\in {\cal J}_{\Theta (G)}$, 
$\tilde{\Theta }(J_{2})\in {\cal J}_{\Theta (G)}$, 
and $\tilde{\Theta }(J_{1})<\tilde{\Theta }(J_{2}).$ 

 Thus, we have proved Corollary~\ref{polyandracor}.
\qed  

\ 

 We now demonstrate Theorem~\ref{polyandrathm2}.

\noindent {\bf Proof of Theorem~\ref{polyandrathm2}:}
We have that for any $j=1,\ldots ,m$, 
$(\Psi (h_{j}))^{-1}(x)=\frac{1}{\deg (h_{j})}(x-\log |a_{j}|)
=\frac{1}{\deg (h_{j})}(x-\frac{-1}{\deg (h_{j})-1}\log |a_{j}|)+
\frac{-1}{\deg (h_{j})-1}\log |a_{j}|$, 
where $x\in \RR .$ 
Hence, it is easy to see that 
$\bigcup _{j=1}^{m}(\Psi (h_{j}))^{-1}([\alpha ,\beta ])
\subset [\alpha ,\beta ].$ 
From the assumption, it follows that 
\begin{equation}
\label{polyandrathm2pf1}
\bigcup _{j=1}^{m}(\Psi (h_{j}))^{-1}([\alpha ,\beta ])=
[\alpha ,\beta ].
\end{equation} 
Moreover, by Lemma~\ref{hmslem}-\ref{bss}, 
we have 
\begin{equation}
\label{polyandrathm2pf2}
\bigcup _{j=1}^{m}(\eta (\Psi (h_{j})))^{-1}(J(\eta (\Psi (G))))=J(\eta (\Psi (G))).
\end{equation}
Furthermore, by Lemma~\ref{cptmsetlem}, $J(\eta (\Psi (G)))$ is a compact 
subset of $\RR .$ 
Applying \cite[Theorem 2.6]{F}, it follows that 
$J(\eta (\Psi (G)))=[\alpha ,\beta ].$ 
Combined with Lemma~\ref{cptmsetlem}, we obtain 
$M(\Psi (G))=[\alpha ,\beta ].$ Hence, $M(\Psi (G))$ is connected. 
Therefore, from Theorem~\ref{polyandrathm1}-\ref{polyandrathm1-1}, 
it follows that 
$J(G)$ is connected.
\qed 

\ 

 We now demonstrate Theorem~\ref{polyandrathm3}.\\ 
{\bf Proof of Theorem~\ref{polyandrathm3}:}
Let $C$ be a set of polynomials of degree two such that $C$ generates $G.$ 
 Suppose that $J(G)$ is disconnected. Then, 
by Theorem~\ref{mainth0}, 
there exist two elements $h_{1},h_{2}\in C$ 
such that the semigroup 
$H=\langle h_{1},h_{2}\rangle $ satisfies that 
$J(H)$ is disconnected. 
For each $j=1,2$, let $a_{j}$ be the coefficient of 
the highest degree term of polynomial $h_{j}.$ 
Let $\alpha := 
\min _{j=1,2}\{ \frac{-1}{\deg (h_{j})-1}\log |a_{j}|\} $ 
and 
$\beta :=\max _{j=1,2}\{ \frac{-1}{\deg (h_{j})-1}\log |a_{j}|\}.$ 
Then we have that $\alpha =\min _{j=1,2}\{- \log |a_{j}|\} $ 
and $\beta =\max _{j=1,2}\{- \log |a_{j}|\} .$ 
Since 
$\Psi (h_{j})^{-1}(x)=\frac{1}{2}(x-\log |a_{j}|)
=\frac{1}{2}(x-(-\log |a_{j}|))+(-\log |a_{j}|)$ for each $j=1,2$, 
we obtain  
$[\alpha ,\beta ]= 
\bigcup _{j=1}^{2}(\Psi (h_{j}))^{-1}([\alpha ,\beta ]).$ 
Hence, by Theorem~\ref{polyandrathm2}, 
it must be true that 
$J(H)$ is connected. However, this is a contradiction. 
Therefore, $J(G)$ must be connected.
\qed 

\ 

We now demonstrate Theorem~\ref{polyandrathm4}. \\ 
{\bf Proof of Theorem~\ref{polyandrathm4}:}
For each $\lambda \in \Lambda $, let 
$b_{\lambda }$ be the fixed point of 
$\Psi (h_{\lambda })$ in $\RR .$ 
It is easy to see that 
$b_{\lambda }=\frac{-1}{\deg (h_{\lambda })-1}\log |a_{\lambda }|$, 
for each $\lambda \in \Lambda .$ 
From the assumption, 
it follows that 
there exists a point $b\in \RR $ such that 
for each $\lambda \in \Lambda $, 
$b_{\lambda }=b.$ 
This implies that 
for any element $g\in G$, 
the fixed point $b(g)\in \RR $ of 
$\Psi (g)$ in $\RR $ is equal to $b.$ 
Hence, we obtain  
$M(\Psi (G))=\{ b\} .$ 
Therefore, $M(\Psi (G))$ is connected. 
From Theorem~\ref{polyandrathm1}-\ref{polyandrathm1-1}, 
it follows that 
$J(G)$ is connected.
\qed 

\ 
\subsection{Proofs of results in \ref{Properties}}
\label{Proof of Properties}
In this section, we prove the results in \ref{Properties}, 
Theorem~\ref{polyandrathm1}-\ref{polyandrathm1-2} and 
Theorem~\ref{polyandrathm1}-\ref{polyandrathm1-3}. 

 In order to demonstrate Theorem~\ref{mainth2},  
Theorem~\ref{polyandrathm1}-\ref{polyandrathm1-2}, and 
Theorem~\ref{polyandrathm1}-\ref{polyandrathm1-3}, we need the following lemma.
\begin{lem}
\label{gdisinflem}
If $G\in {\cal G}_{dis}$, then 
$\infty \in F(G).$ 
\end{lem}
\begin{proof}
Suppose that $G\in {\cal G}_{dis}$ and  
$\infty \in J(G).$ We will deduce a contradiction. 
By Lemma~\ref{inftyj1}, 
the element $J\in \hat{{\cal J}}_{G}$ with 
$\infty \in J$ satisfies that 
$J=\{ \infty \} .$ 
Hence, by Lemma~\ref{appj2lem}, 
for each $n\in \NN $, 
there exists an element $g_{n}\in G$ such that 
$J(g_{n})\subset B(\infty ,\frac{1}{n}).$ 
Let $R>0$ be any number which is sufficiently large so that 
$P^{\ast }(G)\subset B(0,R).$ Since 
we have that $P^{\ast }(G)\subset K(g)$ for each $g\in G$, 
it must hold that there exists a number $n_{0}=n_{0}(R)\in \NN $ 
such that for each $n\geq n_{0}$, 
$B(0,R)<J(g_{n}).$ From Lemma~\ref{greenk1k2}, 
it follows that 
$\lim _{z\rightarrow \infty }
(\log|z|-\varphi (F_{\infty }(g_{n}),z))\rightarrow +\infty 
$ as $n\rightarrow \infty .$
Hence, we see that 
$\frac{-1}{\deg (g_{n})-1}\log |a(g_{n})|\rightarrow 
+\infty $, as $n\rightarrow \infty .$ 
This implies that 
\begin{equation}
\label{typeproppf0}
|a(g_{n})|^{-\frac{1}{\deg(g_{n})-1}}\rightarrow \infty , \mbox{ as }
n\rightarrow \infty .
\end{equation}
Furthermore, by Theorem~\ref{polyandrathm1}-\ref{polyandrathm1-1}, 
we must have that $M(\Psi (G))$ is disconnected. 

 We now consider the polynomial semigroup 
 $H=\{ z\mapsto |a(g)|z^{\deg(g)}\mid g\in G\} \in {\cal G}.$
By Theorem~\ref{repdense}, 
we have 
$J(H)=\overline{\bigcup _{h\in H}J(h)}.$ 
Since the Julia set of polynomial $|a(g)|z^{\deg(g)}$ 
is equal to $\{ z\in \CC \mid |z|=|a(g)|^{-\frac{1}{\deg(g)-1}}\} $, 
it follows that 
\begin{equation}
\label{typeproppf1}
J(H)=\overline{\bigcup _{g\in G}\{ z\in \CC \mid 
|z|=|a(g)|^{-\frac{1}{\deg(g)-1}}\} }, 
\end{equation}
where the closure is taken in $\CCI .$ 
Moreover, $J(\Theta (G))=J(H).$ 
Combining it with (\ref{typeproppf0}), 
(\ref{typeproppf1}), and 
Corollary~\ref{polyandracor}, 
we see that 
\begin{equation}
\label{gdisinflempf1}
\infty \in J(H),  
 \mbox{ and } J(H) \mbox{ is disconnected.}
\end{equation} 
Let $\psi :[0,+\infty ]\rightarrow \hat{\RR }$ be 
the homeomorphism as in the proof of 
Corollary~\ref{polyandracor}. 
By (\ref{typeproppf1}), we have 
\begin{equation}
\label{gdisinflempf2}
\psi (J(H)\cap [0,+\infty ])= M(\Psi (H))=M(\Psi (G)). 
\end{equation}
Moreover, by Lemma~\ref{hmslem}-\ref{invariant}, we have 
\begin{equation}
\label{gdisinflempf3}
h(F(H)\cap [0,+\infty ])\subset F(H)\cap [0,+\infty ], 
\mbox{ for each } h\in H.
\end{equation}
Furthermore, we have that  
\begin{equation}
\label{gdisinflempf4}
\psi \circ h=\Psi (h)\circ \psi \mbox{ on } [0,+\infty ], 
\mbox{ for each } h\in H. 
\end{equation}
Combining (\ref{gdisinflempf2}), (\ref{gdisinflempf3}), 
and (\ref{gdisinflempf4}), we see that 
\begin{equation}
\label{gdisinflempf5}
\Psi (h)(\hat{\RR }\setminus M(\Psi (H)))\subset 
(\hat{\RR }\setminus M(\Psi (H))), \mbox{ for each } 
h\in H.
\end{equation}
By Lemma~\ref{inftyj1} and (\ref{gdisinflempf1}), 
we get that the connected component $J$ of 
$J(H)$ containing $\infty $ satisfies that 
\begin{equation}
\label{gdisinflempf6}
J=\{ \infty \} .
\end{equation} 
Combined with Lemma~\ref{appj2lem}, we see that 
for each $n\in \NN $, there exists an element $h_{n}\in H $ such that 
\begin{equation}
\label{gdisinflempf7}
J(h_{n})\subset B(\infty ,\frac{1}{n}).
\end{equation}
 Combining (\ref{typeproppf1}), (\ref{gdisinflempf2}), 
 (\ref{gdisinflempf6}), and (\ref{gdisinflempf7}),  
we obtain the following claim.\\ 
Claim 1:   
$+\infty $ is a non-isolated point of $M(\Psi (H))$ and 
the connected component of $M(\Psi (H))$ containing $+\infty $ is 
equal to $\{ +\infty \} .$ 

 Let $h\in H$ be an element. Conjugating $G$ by some 
 linear transformation, we may assume that 
$h$ is of the form $h(z)=z^{s}, s\in \NN , s>1.$ 
Hence $\Psi (h)(x)=sx, s>1.$  
Since $0$ is a fixed point of $\Psi (h)$, we have that 
$0\in M(\Psi (H)).$ By Claim 1, 
there exists $c_{1},c_{2}\in [0,+\infty )$ with  
$c_{1}<c_{2}$ such that the open interval $I=(c_{1},c_{2})$ is a 
connected component of $\hat{\RR }\setminus M(\Psi (H)).$ 
 We now show the following claim. \\ 
Claim 2:  
Let $Q=(r_{1},r_{2})\subset (0, + \infty )$ be any connected open interval in  
$\hat{\RR } \setminus M(\Psi (H))$, where 
$0\leq r_{1}<r_{2}<+\infty .$ Then, we have $r_{2}\leq sr_{1}.$ 

 To show this claim, suppose that $sr_{1}< r_{2}.$ Then, 
it implies that $\bigcup _{n\in \NN \cup \{ 0\} }\Psi (h)^{n}(Q)=(r_{1},+\infty ).$  
However, by (\ref{gdisinflempf5}), 
we have $\bigcup _{n\in \NN \cup \{ 0\} }\Psi (h)^{n}(Q)\subset 
\hat{\RR }\setminus M(\Psi (H))$, which implies that 
the connected component $Q'$ of $\hat{\RR }\setminus M(\Psi (H))$ containing $Q$ satisfies that  
$Q'\supset (r_{1},+\infty ).$ 
This contradicts Claim 1. 
Hence, we obtain Claim 2.

 By Claim 2, we obtain $c_{1}>0.$ Let
 $c_{3}\in (0,c_{1})$ be a number so that 
 $c_{2}-c_{3}>s(c_{1}-c_{3}).$ 
Since $c_{1}\in M(\Psi (H))$, there exists an element 
$c\in (c_{3}, c_{1}]$ and an element $h_{1}\in H$ such that 
$\Psi (h_{1})(c)=c$ and $(\Psi (h_{1}))'(c)>1.$ 
Since $c_{2}-c_{3} >s(c_{1}-c_{3})$, we obtain 
\begin{equation}
\label{gdisinflempf8}
c_{2}-c>s(c_{1}-c).
\end{equation}
Let $t:=(\Psi (h_{1}) )'(c)>1.$      
Then, for each $n\in \NN $, 
we have 
$(\Psi (h_{1}) )^{n}(I)=
(t^{n}(c_{1}-c)+c, t^{n}(c_{2}-c)+c).$ 
From Claim 2 and (\ref{gdisinflempf5}), it follows that 
$t^{n}(c_{2}-c)+c\leq s(t^{n}(c_{1}-c)+c)$, for each $n\in \NN .$  
Dividing both sides by $t^{n}$ and then letting $n\rightarrow \infty $, 
we obtain $c_{2}-c\leq s(c_{1}-c).$ However, this contradicts
(\ref{gdisinflempf8}). 
Hence, we must have that 
$\infty \in F(G).$ Thus, we have proved Lemma~\ref{gdisinflem}.
\end{proof} 

\ 

We now demonstrate Proposition~\ref{fcprop}.\\ 
\noindent {\bf Proof of Proposition~\ref{fcprop}:}
Let $U$ be a connected component of $F(G)$ with $U\cap \hat{K}(G)\neq 
\emptyset .$ 
Let $g\in G$ be an element. 
Then we have $\hat{K}(G)\cap F(G)\subset $ int$(K(g)).$ 
Since $h(F(G))\subset F(G)$ and $h(\hat{K}(G)\cap F(G))\subset 
\hat{K}(G)\cap F(G)$ for each $h\in G,$ 
it follows that $h(U)\subset 
$ int$(K(g))$ for each $h\in G.$ Hence 
$U\subset $ int$(\hat{K}(G)).$ From this, 
it is easy to see 
$\hat{K}(G)\cap F(G)=$ int$(\hat{K}(G)).$ 
By the maximum principle, 
we see that $U$ is simply connected.
\qed \\ 

We now demonstrate Theorem~\ref{mainth2}. 

\noindent {\bf Proof of Theorem~\ref{mainth2}:}
%

 First, we show statement \ref{mainth2-2}. 
By Lemma~\ref{gdisinflem}, we have that 
$\infty \in F(G).$ 
 Let $J\in {\cal J}$ be an element 
 such that $\partial F_{\infty }(G)\cap J\neq \emptyset .$ 
 Let $D$ be the unbounded component of 
 $\CCI \setminus J.$ Then 
 $F_{\infty }(G)\subset D$ and $D$ is 
 simply connected. We show 
 $F_{\infty }(G)=D.$ Otherwise, 
 there exists an element $J_{1}\in {\cal J}$ such that  
 $J_{1} \neq J$ and $J_{1}\subset D.$ 
 By Theorem~\ref{mainth1}-\ref{mainth1-1}, we have 
 either $J_{1}< J$ or $J< J_{1}.$ 
 Hence, it follows that $J< J_{1}$ and we have 
 that $J$ is included in a bounded component 
 $D_{0}$ of $\CC \setminus J_{1}.$ Since 
 $F_{\infty }(G)$ is included in the unbounded 
 component $D_{1}$ of $\CCI \setminus J_{1},$ 
 it contradicts $\partial F_{\infty }(G)\cap J\neq 
 \emptyset .$ Hence, $F_{\infty }(G)=D$ and 
 $F_{\infty }(G)$ is simply connected.  

 Next, let $J_{\max }$ be the element of 
 ${\cal J}$ with $\partial F_{\infty }(G)\subset J_{\max }$, and 
 suppose that there exists an element 
 $J\in {\cal J}$ such that 
 $J_{\max }<J.$ Then $J_{\max }$ is 
 included in a bounded component of 
 $\CC \setminus J.$ On the other hand, 
 $F_{\infty }(G)$ is included in the 
 unbounded component of $\CCI \setminus J.$ 
 Since $\partial F_{\infty }(G)\subset J_{\max },$ 
 we have a contradiction. Hence, 
 we have shown that $J\leq J_{\max }$ for each 
 $J\in {\cal J}.$ 

 Therefore, we have shown statement \ref{mainth2-2}.

 Next, we show statement \ref{mainth2-3}. 
 Since $\emptyset \neq P^{\ast }(G) \subset 
 \hat{K}(G),$ we have $\hat{K}(G)\neq \emptyset .$ 
 By Proposition~\ref{fcprop}, we have 
 $\partial \hat{K}(G)\subset J(G).$ 
 Let $J_{1}$ be a connected component 
 of $J(G)$  
 with $J_{1}\cap 
 \partial \hat{K}(G)\neq \emptyset .$ By Lemma~\ref{inftyj1}, $J_{1}\in {\cal J}.$  
 Suppose that there exists an element 
 $J\in {\cal J}$ such that 
 $J<J_{1}.$ Let $z_{0}\in J$ be a point. 
 By Theorem~\ref{repdense}, there exists a 
 sequence $\{ g_{n}\} _{n\in \NN }$ in $G$ such that 
 $d(z_{0},J(g_{n}))\rightarrow 0$ as $n\rightarrow \infty .$ 
 Then by Lemma~\ref{appjlem}, 
 $\sup\limits _{z\in J(g_{n})}d(z,J)\rightarrow 
0$ as $n\rightarrow \infty .$ 
 Since $J_{1}$ is included in the unbounded component 
 of $\CC \setminus J$, it follows that 
 for a large $n\in \NN ,$ 
 $J_{1}$ is included in the unbounded component 
 of $\CC \setminus J(g_{n}).$ 
 However, this causes a contradiction, since 
 $J_{1}\cap \hat{K}(G)\neq \emptyset .$ 
 Hence, by Theorem~\ref{mainth1}-\ref{mainth1-1}, it must hold that 
 $J_{1}\leq J$ for each $J\in {\cal J}.$ 
 This argument shows that 
 if $J_{1}$ and $J_{2}$ are two connected 
 components of $J(G)$ such that 
 $J_{i}\cap \partial \hat{K}(G)\neq \emptyset $ for each 
 $i=1,2$, then $J_{1}=J_{2}.$ Hence,\ 
 we conclude that there exists a unique minimal element 
 $J_{\min }$ in $({\cal J},\leq )$ and   
 $\partial \hat{K}(G)\subset J_{\min }.$ 

 Next, let $D$ be the unbounded component of 
 $\CC \setminus J_{\min }.$ 
 Suppose $D\cap \hat{K}(G)\neq \emptyset .$ 
 Let $x\in D\cap \hat{K}(G)$ be a point. 
 By Theorem~\ref{repdense} and Lemma~\ref{appjlem}, 
 there exists a sequence $\{ g_{n}\} _{n\in \NN }$ in $G$ 
 such that 
 $\sup\limits _{z\in J(g_{n})}d(z,J_{\min})\rightarrow 
0$ as $n\rightarrow \infty .$ 
Then, for a large $n\in \NN $, 
$x$ is in the unbounded component of 
$\CC \setminus J(g_{n}).$ 
However, this is a contradiction, 
since $g_{n}^{l}(x)\rightarrow \infty 
$ as $l\rightarrow \infty $,  
and $x\in \hat{K}(G).$
 Hence, we have shown statement \ref{mainth2-3}.

 Next, we show statement \ref{mainth2-3b}. 
 By Theorem~\ref{mainth0}, there exist 
 $\lambda _{1},\lambda _{2}\in \Lambda $ and 
 connected components $J_{1},J_{2}$ of $J(G)$ 
 such that 
 $J_{1}\neq J_{2}$ and 
 $J(h_{\lambda _{i}})\subset J_{i}$ for each $i=1,2.$ 
 By Lemma~\ref{inftyj1}, we have $J_{i}\in {\cal J}$ for each 
 $i=1,2.$ Then $J(h_{\lambda _{1}})\cap 
 J(h_{\lambda _{2}})=\emptyset .$ 
 Since $P^{\ast }(G) $ is bounded 
 in $\CC $, we may assume 
 $J(h_{\lambda _{2}})< J(h_{\lambda _{1}}).$ 
 Then we have $ K(h_{\lambda _{2}})\subset 
 $ int$(K(h_{\lambda _{1}})) $ and 
 $J_{2}<J_{1}.$ By statement \ref{mainth2-3}, 
 $J_{1}\neq J_{\min }.$ Hence 
 $J(h_{\lambda _{1}})\cap J_{\min }=\emptyset .$  
 Since $P^{\ast }(G)$ is bounded 
 in $\CC $, we have that 
 $ K(h_{\lambda _{2}})$ is connected. 
 Let $U$ be the connected component of 
 int$(K(h_{\lambda _{1}}))$ containing 
 $K(h_{\lambda _{2}}).$ 
 Since $P^{\ast }(G) \subset 
 K(h_{\lambda _{2}}),$ it follows that 
 there exists 
 an attracting fixed point $z_{1}$ of 
 $h_{\lambda _{1}}$ in $K(h_{\lambda _{2}})$ and  
 $U$ is the immediate attracting basin 
 for $z_{1}$ with respect to the 
 dynamics of $h_{\lambda _{1}}.$  
  Furthermore, 
 by Lemma~\ref{fibconncor2}, 
$ h_{\lambda _{1}}^{-1}(J(h_{\lambda _{2}}))$ 
is connected. Therefore,  
$h_{\lambda _{1}}^{-1}(U)=U.$ Hence,  
int$(K(h_{\lambda _{1}}))=U.$ 

 Suppose that there exists an $n\in \NN $ such that 
 $h_{\lambda _{1}}^{-n}(J(h_{\lambda _{2}}))
  \cap J(h_{\lambda _{2}})\neq \emptyset .$ 
  Then, by Lemma~\ref{fibconncor2}, $A:=\bigcup _{s\geq 0 } h_{\lambda _{1}}^{-ns}
  (J(h_{\lambda _{2}})) $ is connected and 
  its closure $\overline{A}$ contains 
  $J(h_{\lambda _{1}}).$ Hence 
  $J(h_{\lambda _{1}})$ and $J(h_{\lambda _{2}})$ 
  are included in the same connected component of $J(G).$ 
  This is a contradiction. 
  Therefore, for each $n\in \NN $, 
  we have $h_{\lambda _{1}}^{-n}(J(h_{\lambda _{2}}))
  \cap J(h_{\lambda _{2}})=\emptyset .$ Similarly, 
  for each $n\in \NN $, 
  we have $h_{\lambda _{2}}^{-n}(J(h_{\lambda _{1}}))
  \cap J(h_{\lambda _{1}})=\emptyset .$ 
  Combining $h_{\lambda _{1}}^{-1}(J(h_{\lambda _{2}}))
  \cap J(h_{\lambda _{2}})=\emptyset $ 
   with $z_{1}\in K(h_{\lambda _{2}})$, 
  we obtain $z_{1}\in $ int$(K(h_{\lambda _{2}})).$ 
Hence, we have proved statement \ref{mainth2-3b}. 
   
 We now prove statement \ref{mainth2ast1}. Let $g\in G$ be an element with 
 $J(g)\cap J_{\min }=\emptyset .$ 
We show the following:\\ 
Claim 2: $J_{\min }<J(g).$ 

 To show the claim, suppose 
 that $J_{\min }$ is included in 
 the unbounded component $U$ of 
 $\CC \setminus J(g).$ 
Since $\emptyset \neq \partial \hat{K}(G)\subset  J_{\min }$, 
it follows that $\hat{K}(G)\cap U\neq \emptyset .$ However, 
this is a contradiction.
Hence, 
 we have shown Claim 2. 

 Combining Claim 2, Theorem~\ref{repdense} and 
 Lemma~\ref{appjlem}, we get that 
 there exists an element $h_{1}\in G$ 
 such that $J(h_{1})<J(g).$  
 From an argument which we have used in the proof of 
 statement \ref{mainth2-3b},  
 it follows that $g$ has an attracting fixed 
 point $z_{g}$ in $\CC $ and 
 int$(K(g))$ consists of only one immediate  
 attracting basin 
 for $z_{g}.$ 
  Hence, we have shown 
  statement \ref{mainth2ast1}.

Next, we show statement \ref{mainth2-4}.
Suppose that int$(\hat{K}(G))=\emptyset .$ We will deduce a 
contradiction.  If int$(\hat{K}(G))=\emptyset $, then 
by Proposition~\ref{fcprop}, we obtain 
$F(G)\cap \hat{K}(G)=\emptyset .$ 
By statement \ref{mainth2-3b}, 
there exist two elements $g_{1}$ and 
$g_{2}$ of $G$ and two elements 
$J_{1}$ and $J_{2}$ of ${\cal J}$ 
such that $J_{1}\neq J_{2}$, 
such that $J(g_{i})\subset J_{i}$ for each 
$i=1,2$, such that $g_{1}$ has an attracting fixed 
point $z_{0}$ in int$(K(g_{2}))$, and such that 
$K(g_{2})\subset $ int$(K(g_{1})).$ 
Since we assume $F(G)\cap \hat{K}(G)=\emptyset ,$ 
we have $z_{0}\in P^{\ast }(G)\subset \hat{K}(G)\subset J(G).$ 
Let $J$ be the connected component of 
$J(G)$ containing $z_{0}.$ 
We now show $J=\{ z_{0}\} .$ 
Suppose $\sharp J\geq 2.$  
Then $J(g_{1})\subset 
\overline{\bigcup _{n\geq 0}g_{1}^{-n}(J)}.$
Moreover,  
by Theorem~\ref{mainth1}-\ref{mainth1-3},  
$g_{1}^{-n}J $ is connected 
for each $n\in \NN .$ 
Since $g_{1}^{-n}(J)\cap J\neq \emptyset $ 
for each $n\in \NN $, we see that 
$\overline{\bigcup _{n\geq 0}g_{1}^{-n}(J)}$ 
is connected. Combining this with 
$z_{0}\in $ int$(K(g_{2}))$,  
$K(g_{2})\subset $ int$(K(g_{1}))$, 
$z_{0}\in J$ and 
$J(g_{1})\subset 
 \overline{\bigcup _{n\geq 0}g_{1}^{-n}(J)}$, 
 we obtain 
 $ \overline{\bigcup _{n\geq 0}g_{1}^{-n}(J)}\cap J(g_{2})
 \neq \emptyset .$ 
 Then it follows that 
 $J(g_{1})$ and $J(g_{2})$ are included 
 in the same connected component of $J(G).$ 
 This is a contradiction. Hence, we have shown 
 $J=\{ z_{0}\} .$ By statement \ref{mainth2-3}, 
 we obtain $\{ z_{0}\} =J_{\min }
 =P^{\ast }(G).$ 
Let 
$\varphi (z):= \frac{1}{z-z_{0}}$ and 
let $\tilde{G}:= \{ \varphi g\varphi ^{-1}\mid g\in G\} .$ 
Then $\tilde{G}\in {\cal G}_{dis}.$ Moreover, 
since $z_{0}\in J(G)$, we have that 
$\infty \in J(\tilde{G}).$ This contradicts Lemma~\ref{gdisinflem}. 
Therefore, we must have that int$(\hat{K}(G))\neq \emptyset .$ 
%

 Since $\partial \hat{K}(G)\subset 
 J_{\min }$ (statement \ref{mainth2-3}) 
 and $\hat{K}(G)$ is bounded, 
 it follows that $\CC \setminus J_{\min } $ 
 is disconnected and 
 $\sharp J_{\min }\geq 2.$ Hence, $\sharp J\geq 2$ 
 for each $J\in {\cal J}=\hat{{\cal J}}.$ 
 Now, let $g\in G$ be an element with 
 $J(g)\cap J_{\min }=\emptyset .$ 
 we show $J_{\min }\neq g^{\ast }(J_{\min }).$ 
 If $J_{\min }=g^{\ast }(J_{\min })$, then 
 $g^{-1}(J_{\min })\subset 
 J_{\min }.$ Since $\sharp J_{\min }\geq 3,$ 
 it follows that $J(g)\subset J_{\min }$, which  
 is a contradiction. 
 Hence, $J_{\min }\neq g^{\ast }(J_{\min })$, and so 
 $J_{\min }<g^{\ast }(J_{\min }).$ Combined with 
 Theorem~\ref{mainth1}-\ref{mainth1-3}, 
 we obtain $g^{-1}(J(G))\cap J_{\min }=\emptyset .$ 
 Since $g(\hat{K}(G))\subset \hat{K}(G)$, 
 we have $g($int$(\hat{K}(G)))\subset $ int$(\hat{K}(G)).$ 
 Suppose $g(\partial \hat{K}(G))\cap 
 \partial \hat{K}(G)\neq \emptyset .$ 
 Then, since $\partial \hat{K}(G)\subset J_{\min }$ 
 (statement \ref{mainth2-3}), 
 we obtain $g(J_{\min })\cap J_{\min }\neq 
 \emptyset .$ This implies 
 $g^{-1}(J_{\min })\cap J_{\min }\neq 
 \emptyset $, 
 which contradicts $g^{-1}(J(G))\cap J_{\min }=\emptyset .$ 
Hence, 
 it must hold 
 $g(\partial \hat{K}(G))\subset $ int$(\hat{K}(G))$, and so 
 $g(\hat{K}(G))\subset $ int$(\hat{K}(G)).$  
Moreover, since $g^{-1}(J(G))\cap J_{\min }=\emptyset $, 
we have that $g(J_{\min })$ is a connected subset of $F(G).$ 
Since $\partial \hat{K}(G)\subset J_{\min }$ 
and $g(\partial \hat{K}(G))\subset $ int$(\hat{K}(G))$, 
Proposition~\ref{fcprop} implies that $g(J_{\min })$ must be contained in int$(\hat{K}(G)).$ 

 By statement \ref{mainth2ast1}, 
 $g$ has a unique attracting fixed point 
 $z_{g}$ in $\CC .$ 
Then, $z_{g}\in P^{\ast }(G) 
\subset \hat{K}(G).$ Hence, 
$z_{g}=g(z_{g})\in g(\hat{K}(G))
\subset $ int$(\hat{K}(G)).$ 
 Hence, we have shown statement \ref{mainth2-4}.

 We now show statement \ref{mainth2ast2}. 
 Since $F_{\infty }(G)$ is simply connected 
 (statement \ref{mainth2-2}), 
we have $\bigcup _{A\in {\cal A}}A\subset \CC .$ 
 Suppose that there exist two distinct elements 
 $A_{1}$ and $A_{2}$ in ${\cal A}$ such that 
 $A_{1}$ is included in the unbounded component of 
 $\CC \setminus A_{2}$, and such that $A_{2}$ is included in the 
 unbounded component of $\CC \setminus A_{1}.$ 
For each $i=1,2$, let $J_{i}\in {\cal J}$ be the 
element that intersects the bounded component of $\CC \setminus A_{i}.$ 
Then, $J_{1}\neq J_{2}.$ 
Since $({\cal J},\leq )$ is 
totally ordered (Theorem~\ref{mainth1}-\ref{mainth1-1}), 
we may assume that $J_{1}<J_{2}.$ 
Then, it implies that $A_{1}<J_{2}<A_{2}$, which is a contradiction. 
Hence, $({\cal A},\leq )$ is totally ordered. 
Therefore, we have proved statement \ref{mainth2ast2}. 

Thus, we have proved Theorem~\ref{mainth2}.
\qed 

\ 

 We now demonstrate Theorem~\ref{mainupthm}.\\ 
{\bf Proof of Theorem~\ref{mainupthm}:} 
First, we show Theorem~\ref{mainupthm}-\ref{mainupthm1}.
If 
$G\in {\cal G}_{con}$, then $J(G)$ is uniformly perfect. 

 We now suppose that $G\in {\cal G}_{dis}.$ 
Let $A$ be an annulus separating $J(G).$
Then $A$ separates $J_{\min }$ and $J_{\max }.$ 
 Let $D$ be the unbounded component 
 of $\CC \setminus J_{\min }$ and 
 let $U$ be the connected component 
 of $\CC \setminus J_{\max } $ containing 
 $J_{\min}.$ Then 
 it follows that $A\subset U\cap D.$ 
 Since $\sharp J_{\min }>1$ and $\infty \in F(G)$ 
 (Theorem~\ref{mainth2}), 
 we get that 
 the doubly connected domain 
 $U\cap D$ satisfies mod $(U\cap D) <\infty .$ 
 Hence, we obtain 
 mod $A\leq $ mod $(U\cap D )<\infty .$ 
 Therefore, $J(G)$ is uniformly perfect.  

 If a point $z_{0}\in J(G)$ is a superattracting fixed point of 
 an element $g\in G$, then, combining uniform perfectness of $J(G)$ and 
 \cite[Theorem 4.1]{HM2}, it follows that $z_{0}\in $ int$(J(G)).$  
Thus, we have shown Theorem~\ref{mainupthm}-\ref{mainupthm1}.
%

Next, we show Theorem~\ref{mainupthm}-\ref{mainupthm2}.
If $G\in {\cal G}$ and $\infty \in J(G)$, then 
by Lemma~\ref{gdisinflem}, we obtain  
$G\in {\cal G}_{con}.$ Moreover, 
Theorem~\ref{mainupthm}-\ref{mainupthm1} 
implies that 
$\infty \in $ int$(J(G)).$ Therefore, we have shown 
Theorem~\ref{mainupthm}-\ref{mainupthm2}.

 We now show Theorem~\ref{mainupthm}-\ref{mainupthm3}.
Suppose that $G\in {\cal G}_{dis}.$ 
 Let $g\in G$ and let $z_{1}\in J(G)\cap \CC $ with 
 $g(z_{1})=z_{1}$ and $g'(z_{1})=0.$ Then, 
$z_{1}\in P^{\ast }(G)\subset \hat{K}(G).$  
By Theorem~\ref{mainth2}-\ref{mainth2-3},  
  we obtain 
 $z_{1}\in J_{\min }.$
%
Moreover, 
Theorem~\ref{mainupthm}-\ref{mainupthm1} 
implies that 
 $z_{1}\in $ int$(J(G)).$ Combining this and 
 $z_{1}\in J_{\min }$, we obtain $z_{1}\in $ int$(J_{\min }).$ 
 By Theorem~\ref{mainth2}- \ref{mainth2-4-2}, we 
 obtain  
 $J(g)\subset J_{\min }.$ 

 Hence, we have shown Theorem~\ref{mainupthm}.
\qed 

\ 

We now demonstrate Theorem~\ref{polyandrathm1}-\ref{polyandrathm1-2}.\\ 
{\bf Proof of Theorem~\ref{polyandrathm1}-\ref{polyandrathm1-2}:}
Suppose $G\in {\cal G}_{dis}. $ Then, by Lemma~\ref{gdisinflem}, 
we obtain $\infty \in F(G).$ Hence, 
there exists a number $R>0$ such that 
for each $g\in G$, $J(g)<\partial B(0,R).$ 
From Lemma~\ref{greenk1k2}, it follows that 
there exists a constant $C_{1}>0$ such that 
for each $g\in G$, 
$\frac{-1}{\deg(g)-1}\log |a(g)|< C_{1}.$ 
This implies that there exists a constant $C_{2}\in \RR $ such that 
\begin{equation}
\label{polyandrathm1-2pf1}
M(\Psi (G))\subset [-\infty ,C_{2}].
\end{equation} 
Moreover, by Theorem~\ref{mainth2}-\ref{mainth2-4}, 
we have that int$(\hat{K}(G))\neq \emptyset .$ 
Let $B$ be a closed disc in int$(\hat{K}(G)).$ 
Then it must hold that 
for each $g\in G$, $B<J(g).$ 
Hence, by Lemma~\ref{greenk1k2}, there exists a constant 
$C_{3}\in \RR $ such that 
for each $g\in G$, 
$C_{3}\leq \frac{-1}{\deg(g)-1}\log |a(g)|.$ 
Therefore, we obtain  
\begin{equation}
\label{polyandrathm1-2pf2}
M(\Psi (G))\subset [C_{3},+\infty ].
\end{equation}
Combining (\ref{polyandrathm1-2pf1}) and (\ref{polyandrathm1-2pf2}), 
we obtain $M(\Psi (G))\subset \RR .$ 
Let $C_{4}$ be a large number so that 
$M(\Psi (G))\subset D(0,C_{4}).$ 
Since for each $g\in G$, 
the repelling fixed point 
$-\frac{1}{\deg(g)-1}\log |a(g)|$ of 
$\eta (\Psi (g))$ belongs to $D(0,C_{4})\cap \RR $, 
we see that for each $z\in \CC \setminus D(0,C_{4})$, 
$|\eta (\Psi (g))(z)|=
|\deg (g)(z-\frac{-1}{\deg(g)-1}\log |a(g)|)+
\frac{-1}{\deg(g)-1}\log |a(g)||
\geq \deg (g)C_{4}-(\deg (g)-1)C_{4}=
C_{4}.$ It follows that 
$\infty \in F(\eta (\Psi (G))).$ 
Combining this and Theorem~\ref{repdense}, 
we obtain  
$M(\Psi (G))=J(\eta (\Psi (G)))$, 
if $\sharp (J(\eta (\Psi (G))))\geq 3.$ 

 Suppose that $\sharp (J(\eta (\Psi (G))))= 2.$ 
Let $g\in G$ be an element and let $b\in \RR $ be 
the unique fixed point of $\Psi (g)$ in $\RR .$ 
Then, since $\infty \in F(\eta (\Psi (G)))$, 
there exists a point $c\in (J(\eta (\Psi (G)))\cap \CC )\setminus 
\{ b\} .$ By Lemma~\ref{hmslem}-\ref{invariant}, 
$(\eta (\Psi (g)))^{-1}(c)\in J(\eta (\Psi (G)))\setminus \{ b,c\} .$ 
This contradicts $\sharp (J(\eta (\Psi (G))))= 2.$ Hence it must hold that  
 $\sharp (J(\eta (\Psi (G))))\neq  2.$

 Suppose that $\sharp (J(\eta (\Psi (G))))=1 .$   
Since $M(\Psi (G))\subset \RR $ 
and $M(\Psi (G))\cap \RR \subset J(\eta (\Psi (G)))$, 
it follows that $M(\Psi (G))=J(\eta (\Psi (G))).$ 

 Therefore, we always have that $M(\Psi (G))=J(\eta (\Psi (G))).$ 
Thus, we have proved Theorem~\ref{polyandrathm1}-\ref{polyandrathm1-2}.
\qed 

\ 

We now demonstrate Theorem~\ref{polyandrathm1}-\ref{polyandrathm1-3}.\\ 
{\bf Proof of Theorem~\ref{polyandrathm1}-\ref{polyandrathm1-3}:} 
By Theorem~\ref{polyandrathm1}-\ref{polyandrathm1-1} and 
Theorem~\ref{polyandrathm1}-\ref{polyandrathm1-2}, 
the statement holds.
\qed 

\ 

We now demonstrate Proposition~\ref{orderjprop}.\\ 
\noindent {\bf Proof of Proposition~\ref{orderjprop}:}
First, we show statement \ref{orderjprop1}. Let $g\in Q_{1}.$ 
 We show the following:\\ 
\noindent Claim 1: For any element $J_{3}\in {\cal J}$ with 
$J_{1}\leq J_{3}$, we have 
$J_{1}\leq g^{\ast }(J_{3}).$ \\ 
To show this claim, let 
$J\in {\cal J}$ be an element with $J(g)\subset J.$
We consider the following two cases;\\   
Case 1: $J\leq J_{3}$, and\\ 
Case 2: $J_{1}\leq J_{3}\leq J.$ 

 Suppose that we have Case 1. Then, $J_{1}\leq J=g^{\ast }(J)\leq 
 g^{\ast }(J_{3}).$ Hence, the statement of Claim 1 is true. 

 Suppose that we have Case 2. 
If we have  $g^{\ast }(J_{3})< J_{3}$,  
then, we have $(g^{n})^{\ast }(J_{3})\leq g^{\ast }(J_{3})<J_{3}\leq J$ 
for each $n\in \NN .$ Hence, $\inf \{ d(z,J)\mid z\in g^{-n}(J_{3}),
n\in \NN \} >0.$ However, since $J(g)\subset J$ and 
$\sharp J_{3}\geq 3$, we obtain a contradiction. 
Hence, we must have $J_{3}\leq g^{\ast }(J_{3})$, which implies 
$J_{1}\leq J_{3}\leq g^{\ast }(J_{3}).$ 
Hence, we conclude that Claim 1 holds. 

 Now, let $K_{1}:= J(G)\cap (J_{1}\cup A_{1}).$ Then, 
 by Claim 1,\ we obtain $g^{-1}(K_{1})\subset K_{1}$, 
 for each $g\in Q_{1}.$ From Lemma~\ref{hmslem}-\ref{backmin}, 
 it follows that $J(H_{1})\subset K_{1}.$  
Hence, we have shown statement \ref{orderjprop1}.

 Next, we show statement \ref{orderjprop2}. Let $g\in Q_{2}.$ 
Then, by the same method as that of the proof 
 of Claim 1, we obtain the following.\\ 
 \noindent Claim 2: For any element $J_{4}\in {\cal J}$ with 
 $J_{4}\leq J_{2}$, we have 
 $g^{\ast }(J_{4})\leq J_{2}.$ 

 Now, let $K_{2}:=J(G)\cap (\CC \setminus A_{2}).$ 
 Then, by Claim 2, we obtain $g^{-1}(K_{2})\subset 
 K_{2}$, for each $g\in Q_{2}.$ From Lemma~\ref{hmslem}-\ref{backmin}, 
 it follows that $J(H_{2})\subset K_{2}.$ Hence, we have shown 
 statement \ref{orderjprop2}.

 Next, we show statement \ref{orderjprop3}. 
 By statements \ref{orderjprop1} and \ref{orderjprop2}, 
 we obtain 
 $J(H)\subset J(H_{1})\cap J(H_{2})\subset K_{1}\cap K_{2}
 \subset (\CC \setminus A_{2})\cap (J_{1}\cup A_{1})
 \subset J_{1}\cup (A_{1}\setminus A_{2}).$ 

 Hence, we have proved Proposition~\ref{orderjprop}.   
\qed 

\ 

We now demonstrate Proposition~\ref{bminprop}.\\ 
{\bf Proof of Proposition~\ref{bminprop}:}
Suppose that for any $h\in \G $, $J(h)\cap J_{\max }=\emptyset .$ 
Then, since $\sharp J_{\max}\geq 3$ (Theorem~\ref{mainth2}-\ref{mainth2-4-1}), 
we get that for any $h\in \G $, 
$h^{-1}(J_{\max })\cap J_{\max }=\emptyset .$ 
Combining it with Theorem~\ref{mainth1}-\ref{mainth1-3}, 
it follows that for any $h\in \G $, 
$h^{-1}(J(G))\cap J_{\max }=\emptyset .$ 
However, since 
$J(G)=\bigcup _{h\in \G}h^{-1}(J(G))$ (Lemma~\ref{hmslem}-\ref{bss}), 
it causes a contradiction. Hence, there must be an element 
$h_{1}\in \G $ such that $J(h_{1})\subset J_{\max }.$

 By the same method as above, we can show that there exists an 
 element $h_{2}\in \G $ such that $J(h_{2})\subset J_{\min }.$  
\qed 

\subsection{Proofs of results in \ref{Poly}}
\label{Proof of Poly}
In this section, we prove the results in \ref{Poly}.

We now prove Theorem~\ref{fcthm}.\\ 
\noindent {\bf Proof of Theorem~\ref{fcthm}:} 
Combining the assumption and Theorem~\ref{mainth1}-\ref{mainth1-3}, 
we get that for each $h\in \G $ and each $j\in \{ 1,\ldots ,n\} $, 
there exists a $k\in \{ 1,\ldots ,n\} $ with 
$h^{-1}(J_{j})\subset J_{k}.$ 
Hence, 
\begin{equation}
\label{fcthmpfeq1}
h^{-1}(\bigcup _{j=1}^{n}J_{j})\subset \bigcup _{j=1}^{n}J_{j}, 
\mbox{ for each } h\in \G .
\end{equation}
Moreover, by Theorem~\ref{mainth2}-\ref{mainth2-4-1}, 
we obtain 
\begin{equation}
\label{fcthmpfeq2}
\sharp (\bigcup _{j=1}^{n}J_{j})\geq 3.
\end{equation}
Combining (\ref{fcthmpfeq1}), (\ref{fcthmpfeq2}), and 
Lemma~\ref{hmslem}-\ref{backmin}, 
it follows that $J(G)\subset \bigcup _{j=1}^{n}J_{j}.$ 
Hence, $J(G)=\bigcup _{j=1}^{n}J_{j}.$ 
Therefore, we have proved Theorem~\ref{fcthm}.
\qed 

 We now prove Proposition~\ref{fincomp}.

\noindent {\bf Proof of Proposition~\ref{fincomp}:} 
Let $n\in \NN $ with $n>1$ and let 
$\epsilon $ be a number with $0<\epsilon <\frac{1}{2}.$ 
For each $j=1,\ldots ,n$, 
let $\alpha _{j}(z)=\frac{1}{j}z^{2}$ and  
let $\beta _{j}(z)=\frac{1}{j}(z-\epsilon )^{2}+\epsilon .$  

For any large $l\in\NN $, 
there exists an open neighborhood $U$  
 of $\{ 0,\epsilon \} $ with  
 $U\subset \{ z\mid |z|<1\} $ and 
 a open neighborhood $V$ of 
 $(\alpha _{1}^{l},\ldots ,\alpha _{n}^{l},
 \beta _{1}^{l},\ldots ,\beta _{n}^{l})$ in 
 (Poly)$^{2n}$ such that 
 for each $(h_{1},\ldots ,h_{2n})\in V$, 
 we have $\bigcup _{j=1}^{2n}h_{j}(U)\subset U$ and 
 $\bigcup _{j=1}^{m}C(h_{j})\cap \CC \subset U$, 
 where $C(h_{j})$ denotes the set of all critical points of 
 $h_{j}.$ Then, by Remark~\ref{pcbrem}, for each $(h_{1},\ldots ,h_{2n})\in V $,  
 $\langle h_{1},\ldots ,h_{2n}\rangle \in {\cal G}.$  
If $l$ is large enough and $V$ is so small, then, 
for each $(h_{1},\ldots ,h_{2n})\in V$, the set  
 $I_{j}:=
J(h_{j})\cup 
J(h_{j+n})$ is connected, for each $j=1,\ldots ,n$, and   
we have:
\begin{equation}
\label{finpropeq1}
(h_{i})^{-1}(I_{j})\cap I_{i}\neq \emptyset , 
(h_{i+n})^{-1}(I_{j})\cap I_{i}\neq \emptyset ,
\end{equation}
for each $(i,j).$ 
Furthermore, for a closed annulus $A=\{z\mid \frac{1}{2}
\leq |z|\leq n+1\} $,   
if $l\in \NN $ is large enough and $V$ is so small, then 
for each $(h_{1},\ldots ,h_{2n})\in V$, 
$\bigcup _{j=1}^{2n}(h_{j})^{-1}(A)
\subset \mbox{int}(A)$  
 and 
$\{ (h_{j})^{-1}(A)\cup (h_{j+n})^{-1}(A)\} _{j=1}^{n}$ 
are mutually disjoint. 
Combining it with Lemma~\ref{hmslem}-\ref{backmin} and 
Lemma~\ref{hmslem}-\ref{bss}, we get 
that for each $(h_{1},\ldots ,h_{2n})\in V$,\ 
$J(\langle h_{1},\ldots ,h_{2n}\rangle)\subset A$ and 
$\{ J_{j}\} _{j=1}^{n}$ are 
mutually disjoint, where 
$J_{j}$ denotes the connected component of 
$J(\langle h_{1},\ldots ,h_{2n}\rangle )$ containing 
$I_{j}=J(h_{j})\cup J(h_{j+n}).$ 
Combining it with (\ref{finpropeq1}) and Theorem~\ref{fcthm}, 
it follows that for each $(h_{1},\ldots ,h_{2n})\in V$, 
 the polynomial semigroup 
$G=\langle h_{1},\ldots ,h_{2n}\rangle $ satisfies that  
$\sharp (\hat{{\cal J}}_{G})=n.$ 
\qed  

\ 

 To prove Theorem~\ref{countthm}, 
we need the following notation.
\begin{df} 
\ 
\begin{enumerate}
\item 
 Let $X$ be a metric space.
 Let $h_{j}:X\rightarrow X\ (j=1,\ldots ,m)$ be a 
 continuous map. 
 Let $G=\langle h_{1},\ldots ,h_{m}\rangle $ be 
 the semigroup generated by $\{ h_{j}\} .$
 A non-empty compact subset $L$ of $X$ is 
 said to be a {\bf backward self-similar set 
 with respect to $\{ h_{1},\ldots ,h_{m}\} $} 
 if  
(a) $ L=\bigcup _{j=1}^{m}h_{j}^{-1}(L)$ 
 and (b) 
$g^{-1}(z)\neq \emptyset $ for each 
 $z\in L$ and $g\in G.$
 For example, if 
 $G=\langle h_{1},\ldots ,h_{m}\rangle $
 is a  
 finitely generated rational semigroup, then  
 the Julia set $J(G)$ is 
 a backward self-similar set with respect to 
 $\{ h_{1},\ldots ,h_{m}\} .$ 
 (See Lemma~\ref{hmslem}-\ref{bss}.)
\item We set $\Sigma _{m}:=
\{ 1,\ldots ,m\} ^{\NN }.$
For each $x=(x_{1},x_{2},\ldots ,)\in \Sigma _{m} $,  
we set 
$ L_{x}:=\bigcap _{j=1}^{\infty }
h_{x_{1}}^{-1}\cdots h_{x_{j}}^{-1}(L) \ (\neq \emptyset ).$ 

\item 
For a finite word $w=(w_{1},\ldots ,w_{k})\in 
\{ 1\ldots ,m\} ^{k}$, 
we set 
$h_{w}:=h_{w_{k}}\circ \cdots \circ h_{w_{1}}.$
 
\item Under the notation of \cite[page 110--page 115]{Sp},  
for any $k\in \NN $, 
let $\Omega _{k}=\Omega _{k}(L,\{ h_{1},\ldots ,h_{m}\} )$  
be the graph (one-dimensional simplicial complex) whose vertex set is  
$\{ 1,\ldots ,m\}^{k}$ and 
that satisfies that mutually different 
$w^{1},w^{2}\in \{ 1,\ldots ,m\}^{k}$ makes 
a $1$-simplex if and only if 
$\bigcap _{j=1}^{2} h_{w^{j}}^{-1}(L)\neq \emptyset .$  
 Let $\varphi _{k}:\Omega _{k+1}\rightarrow 
\Omega _{k}$ be the simplicial map 
defined by:  
$(w_{1},\ldots ,w_{k+1})\mapsto (w_{1},\ldots ,w_{k})$   
for each 
$ (w_{1},\ldots ,w_{k+1})\in \{ 1,\ldots ,m\} ^{k+1}.$
Then $\{ \varphi _{k}:\Omega _{k+1}\rightarrow 
\Omega _{k}\} _{k\in \NN }$  makes an inverse system of 
simplicial maps. Let $|\Omega _{k}|$ be the realization (\cite{Sp}) of $\Omega _{k}.$ 
As in \cite{Sp}, we embed the vertex set $\{ 1,\ldots ,m\} ^{m}$ into $|\Omega _{k}|.$  
\item 
Let ${\cal C}(|\Omega _{k})|)$ be the 
set of all connected components of the realization 
$|\Omega _{k} |$ of $\Omega _{k}.$
Let $ \{ (\varphi _{k})_{\ast }: {\cal C}(|\Omega _{k+1}|)
\rightarrow {\cal C}(|\Omega _{k}|) \} _{k\in \NN }$ be 
the inverse system induced by 
$\{ \varphi _{k}\} _{k}.$      
\end{enumerate}
\end{df}
\noindent {\bf Notation:} 
We fix an $m\in \NN .$  
We set ${\cal W}^{\ast }:= 
\bigcup _{k=1}^{\infty } 
\{ 1,\ldots ,m\} ^{k}$ (disjoint union) and 
$\tilde{{\cal W}}:= {\cal W}^{\ast }\cup 
\Sigma _{m} $ (disjoint union). 
For an element $x\in \tilde{{\cal W}}$, 
we set $|x|= k$ if $x\in 
\{ 1,\ldots ,m\} ^{k} $, and 
$|x|=\infty $ if $x\in \Sigma _{m}.$ 
(This is called the word length of $x.$)  
For any $x\in \tilde{{\cal W}} $ and any 
$j\in \NN $ with $j\leq |x|$, we set 
$x|j:= (x_{1},\ldots ,x_{j})\in 
\{ 1,\ldots ,m\} ^{j} .$   
For any $x^{1}=(x_{1}^{1},\ldots ,
x_{p}^{1})\in {\cal W}^{\ast }$ 
and any $x^{2}=(x_{1}^{2},x_{2}^{2},\ldots )\in \tilde{{\cal W}}$, 
we set $
x^{1}x^{2}:=(x_{1}^{1},\ldots ,x_{p}^{1},
x_{1}^{2},x_{2}^{2},\ldots )\in \tilde{{\cal W}}.$

 To prove Theorem~\ref{countthm}, 
we need the following lemmas.
\begin{lem}
\label{phionto} 
Let $L$ be a backward self-similar set with respect to 
$\{ h_{1},\ldots ,h_{m}\} .$ 
Then, for each $k\in \NN $, the map 
$|\varphi _{k}|: |\Omega _{k+1}|\rightarrow 
|\Omega _{k}|$ induced from 
$\varphi _{k}:\Omega _{k+1}\rightarrow \Omega _{k}$ 
is surjective. 
In particular, 
$ (\varphi _{k})_{\ast }: 
{\cal C}(|\Omega _{k+1}|)\rightarrow 
{\cal C}(|\Omega _{k}|)$ is surjective.
\end{lem} 
\begin{proof}
 Let $x^{1},x^{2}\in 
 \{1 ,\ldots ,m\} ^{k}$ and 
 suppose that $\{ x^{1},x^{2}\} $ makes a $1$-simplex 
 in $\Omega _{k}.$ Then 
 $h_{x^{1}}^{-1}(L)\cap h_{x^{2}}^{-1}(L)
 \neq \emptyset .$ Since 
 $L=\bigcup _{j=1}^{m}h_{j}^{-1}(L)$, 
 there exist $x^{1}_{k+1}$ and 
 $x^{2}_{k+1}$ in $\{ 1,\ldots ,m\} $ 
 such that $
 h_{x^{1}}^{-1}h_{x^{1}_{k+1}}^{-1}(L)
 \cap h_{x^{2}}^{-1}h_{x^{2}_{k+1}}^{-1}(L)
 \neq \emptyset .$ 
 Hence, 
 $\{ x^{1}x^{1}_{k+1},x^{2}x^{2}_{k+1}\} $ 
 makes a $1$-simplex in 
 $\Omega _{k+1}.$ Hence 
 the lemma holds.
\end{proof}
\begin{lem}
\label{1disconlem}
Let $m\geq 2$ and 
let $L$ be a backward self-similar 
set with respect to 
$\{ h_{1},\ldots ,h_{m}\} .$ 
Suppose 
that for each $j$ with 
$j\neq 1$, $h_{1}^{-1}(L)\cap 
h_{j}^{-1}(L)=\emptyset .$
For each $k$, let $C_{k}\in {\cal C}(|\Omega _{k}|)$ 
 be the element containing 
 $(1,\ldots ,1)\in \{ 1,\ldots ,m\} ^{k}.$ 
 Then, we have the following.
 \begin{enumerate}
 \item \label{1disconlem1}
 For each $k\in \NN $, 
 $C_{k}=\{ (1,\ldots ,1)\} .$
  
\item \label{1disconlem2} For each $k\in \NN $, 
$\sharp ({\cal C}(|\Omega _{k}|))
<\sharp ({\cal C}(|\Omega _{k+1}|)).$
\item \label{1disconlem3}
$L$ has infinitely many connected components. 
\item \label{1disconlem4}
Let $x:=(1,1,1,\ldots )\in \Sigma _{m}$ and let 
$x'\in \Sigma _{m}$ be an element with $x\neq x'.$ Then, 
for any $y\in L_{x}$ and $y'\in L_{x'}$, there exists no 
connected component $A$ of $L$ such that 
$y\in A$ and $y'\in A.$  
\end{enumerate}
\end{lem}
\begin{proof}
 We show statement \ref{1disconlem1} by induction on $k.$ 
 We have  $C_{1}=\{ 1\} .$ 
 Suppose $C_{k}=\{ (1,\ldots ,1)\} .$ 
 Let $w\in \{ 1,\ldots ,m\} ^{k+1}\cap C_{k+1}$ be 
 any element. 
 Since 
 $(\varphi _{k})_{\ast }(C_{k+1})=C_{k}$, 
 we have $\varphi _{k}(w)=(1,\ldots ,1)\in \{ 1,\ldots ,m\} ^{k}.$ 
 Hence, $w|k=(1,\ldots ,1)\in 
 \{ 1,\ldots ,m\} ^{k}.$ 
 Since $
 h_{1}^{-1}(L)\cap h_{j}^{-1}(L)=\emptyset $ 
for each $j\neq 1$, 
we obtain $w=(1,\ldots ,1)\in \{ 1,\ldots ,m\} ^{k+1}.$ 
Hence, the induction is completed. Therefore, 
we have shown statement \ref{1disconlem1}. 

 Since both $(1,\ldots 1,1)\in \{ 1,\ldots ,m\} ^{k+1}$ 
 and $(1,\ldots ,1,2)\in \{ 1,\ldots ,m\} ^{k+1}$ are 
 mapped to $(1,\ldots ,1)\in \{ 1,\ldots ,m\} ^{k}$ under 
 $\varphi _{k}$, by statement \ref{1disconlem1} and 
 Lemma~\ref{phionto}, 
 we obtain statement \ref{1disconlem2}. 
 For each $k\in \NN $, 
 we have 
\begin{equation}
\label{1disconlemeq1}
 L=\coprod _{C\in {\cal C}(|\Omega _{k}|)}
 \ \bigcup _{w\in \{ 1,\ldots ,m\} ^{k}\cap C}
 h_{w}^{-1}(L).
\end{equation} 
 Hence, by statement \ref{1disconlem2}, we conclude  
 that $L$ has infinitely many connected components.

 We now show statement \ref{1disconlem4}. 
 Let $k_{0}:=\min \{ l\in \NN \mid x_{l}'\neq 1\} . $ 
 Then, by (\ref{1disconlemeq1}) and statement 
 \ref{1disconlem1}, we get that there exist compact sets 
 $B_{1}$ and $B_{2}$ in $L$  such that 
 $B_{1}\cap B_{2}=\emptyset ,\ B_{1}\cup B_{2}=L,\ 
L_{x}\subset (h_{1}^{k_{0}})^{-1}(L)\subset  B_{1},$ and 
$L_{x'}\subset h_{x_{1}'}^{-1}\cdots h_{x_{k_{0}}'}^{-1}(L)\subset B_{2}.$ 
Hence, statement \ref{1disconlem4} holds.
\end{proof}
We now demonstrate Theorem~\ref{countthm}.

\noindent {\bf Proof of Theorem~\ref{countthm}:} 
By Theorem~\ref{mainth2}-\ref{mainth2-2} or Remark~\ref{hatjcptrem}, we have 
$\hat{{\cal J}}={\cal J}.$
Let $J_{1}\in \hat{{\cal J}}$ be the element containing 
$J(h_{m}).$ By Theorem~\ref{mainth0}, we must have 
$J_{0}\neq J_{1}.$ Then, by Theorem~\ref{mainth1}-\ref{mainth1-1}, 
we have the following two possibilities.\\ 
Case 1. $J_{0}<J_{1}.$\\ 
Case 2. $J_{1}<J_{0}.$

 Suppose we have case 1. Then, by Proposition~\ref{bminprop}, 
 we have that $J_{0}=J_{\min }$ and $J_{1}=J_{\max }.$ 
 Combining it with the assumption and Theorem~\ref{mainth1}-\ref{mainth1-3}, 
 we obtain 
\begin{equation}
\label{cteq1}
\bigcup _{j=1}^{m-1}h_{j}^{-1}(J_{\max })\subset J_{\min }. 
\end{equation}
By (\ref{cteq1}) and Theorem~\ref{mainth1}-\ref{mainth1-3}, 
we get 
\begin{equation}
\label{cteq1-1}
\bigcup _{j=1}^{m-1}h_{j}^{-1}(J(G))\subset J_{\min }.
\end{equation}
Moreover, since $J(h_{m})\cap J_{\min }=\emptyset $, 
Theorem~\ref{mainth2}-\ref{mainth2-4-2} implies that 
\begin{equation}
\label{cteq2}
h_{m}^{-1}(J(G))\cap J_{\min }=\emptyset .
\end{equation}  
Then, by (\ref{cteq1-1}) and (\ref{cteq2}), we get 
\begin{equation}
\label{cteq3}
h_{m}^{-1}(J(G))\cap \left(\bigcup _{j=1}^{m-1}h_{j}^{-1}(J(G))\right)
=\emptyset .
\end{equation}
We now consider the backward self-similar set $J(G)$ 
with respect to $\{ h_{1},\ldots ,h_{m}\} .$ 
By Lemma~\ref{hmslem}-\ref{bss}, we have 
\begin{equation}
\label{cteqast1}
J(G)=\bigcup _{w\in \Sigma _{m}}(J(G))_{w}.
\end{equation}
By Theorem~\ref{mainth2}-\ref{mainth2ast1} and Theorem~\ref{mainth2}-\ref{mainth2-4-2}, 
we obtain $(J(G))_{m^{\infty }}=J(h_{m})$, where $m^{\infty }=(m,m,\ldots )\in \Sigma _{m}.$ 
Combining this with (\ref{cteq3}), Lemma~\ref{1disconlem}, 
and (\ref{cteqast1}), 
we obtain 
\begin{equation}
\label{cteq4}
J_{\max }=(J(G))_{m^{\infty }}= J(h_{m}).
\end{equation}
Furthermore, by (\ref{cteq3}) and Lemma~\ref{1disconlem}, 
we get 
\begin{equation}
\label{cteq5}
\sharp (\hat{{\cal J}})\geq \aleph _{0}.
\end{equation}
Let $x=(x_{1},x_{2},\ldots )\in \Sigma _{m}$ be any element with $x\neq m^{\infty }$ 
and let $l:= \min \{ s\in \NN \mid x_{s}\neq m\} .$ 
Then, by (\ref{cteq1-1}), we have 
\begin{equation}
\label{cteq6}
(J(G))_{x}=\bigcap _{j=1}^{\infty }h_{x_{1}}^{-1}\cdots 
h_{x_{j}}^{-1}(J(G))\subset (h_{m}^{l-1})^{-1}(J_{\min }).
\end{equation}
Combining (\ref{cteqast1}) with (\ref{cteq4}) and (\ref{cteq6}),  
we obtain 
\begin{equation}
\label{cteq7}
J(G)=J_{\max }\cup \bigcup _{n\in \NN \cup \{ 0\} }
h_{m}^{-n}(J_{\min }).
\end{equation}
By (\ref{cteq5}) and (\ref{cteq7}), we get 
$\sharp (\hat{{\cal J}})=\aleph _{0}.$ 
Moreover, combining (\ref{cteq4}), (\ref{cteq7}), 
Theorem~\ref{mainth2}-\ref{mainth2ast1} and 
Theorem~\ref{mainth2}-\ref{mainth2-4-2}, we get that 
for each $J\in \hat{{\cal J}}$ with $J\neq J_{\max }$, 
there exists no sequence $\{ C_{j}\} _{j\in \NN }$ 
of mutually distinct elements of $\hat{{\cal J}}$ such that  
$\min _{z\in C_{j}}d(z,J)\rightarrow 0$ as $j\rightarrow \infty .$ 
Hence, all statements of Theorem~\ref{countthm} are true, provided that 
we have case 1.

 We now assume case 2: $J_{1}<J_{0}.$ 
Then, by Proposition~\ref{bminprop}, 
 we have that $J_{0}=J_{\max }$ and $J_{1}=J_{\min }.$ 
Since $J(h_{j})\subset J_{0}$ for each $j=1,\ldots ,m-1$, 
and since $J_{0}\neq J_{\min }$, 
Theorem~\ref{mainth2}-\ref{mainth2-4-2} 
implies that for each $j=1,\ldots ,m-1$, 
$h_{j}(J(h_{m}))\subset $ int$(K(h_{m})).$  
Hence, for each $j=1,\ldots ,m$, 
$h_{j}(K(h_{m}))\subset $ $K(h_{m}).$ 
Therefore, 
$ \mbox{int}(K(h_{m}))\subset F(G).$
Thus, we obtain 
$(J(G))_{m^{\infty }}=J(h_{m}).$ 
Combining this with the same method as that of case 1, we obtain  
\begin{equation}
\label{cteq8}
J_{\min }=(J(G))_{m^{\infty }}= J(h_{m}),
\end{equation}
\begin{equation}
\label{cteq9} 
J(G)=J_{\min }\cup \bigcup _{n\in \NN \cup \{ 0\} }
h_{m}^{-n}(J_{\max }),
\end{equation}
and 
\begin{equation}
\label{cteq10}
\sharp (\hat{{\cal J}})=\aleph _{0}. 
\end{equation}
Moreover, by (\ref{cteq8}) and (\ref{cteq9}), 
we get that for each $J\in \hat{{\cal J}}$ with 
$J\neq J_{\min }$, there exists no sequence $\{ C_{j}\}_{j\in \NN }$ 
of mutually distinct elements of $\hat{{\cal J}}$ 
such that $\min _{z\in C_{j}}d(z,J)\rightarrow 0$ as $j\rightarrow \infty .$ 
Hence, we have shown Theorem~\ref{countthm}.     
\qed 

\ 

We now demonstrate Proposition~\ref{countprop}.

\noindent {\bf Proof of Proposition~\ref{countprop}:}
Let $0<\epsilon <\frac{1}{2}$ and let 
$\alpha _{1}(z):=z^{2},\alpha _{2}(z):=(z-\epsilon )^{2}+\epsilon , $ 
and $\alpha _{3}(z):=\frac{1}{2}z^{2}.$ 
If we take a large $l\in \NN $, then there exists an open 
neighborhood $U$ of $\{ 0,\epsilon \} $ with 
$U\subset \{ |z|<1\} $ and a neighborhood $V$ of 
$(\alpha _{1}^{l}, \alpha _{2}^{l},\alpha _{3}^{l})$ 
in (Poly)$^{3}$ such that for each $(h_{1},h_{2},h_{3})\in V$, 
we have 
$\bigcup _{j=1}^{3}h_{j}(U)\subset U$
and $\bigcup _{j=1}^{3}C(h_{j})\cap \CC \subset U$, where 
$C(h_{j})$ denotes the set of all critical points of $h_{j}.$ 
Then, by Remark~\ref{pcbrem}, for each $(h_{1},h_{2},h_{3})\in V$,  
$\langle h_{1},h_{2},h_{3}\rangle \in {\cal G}.$ 
Moreover, if we take an $l$ large enough and $V$ so small, 
then for each $(h_{1},h_{2},h_{3})\in V$, we have 
that: 
\begin{enumerate}
\item \label{cp0}
$J(h_{1})<J(h_{3})$;
\item \label{cp1}
$J(h_{1})\cup J(h_{2})$ is connected; 
\item \label{cp2}
$h_{i}^{-1}(J(h_{3}))\cap (J(h_{1})\cup J(h_{2}))\neq \emptyset ,$
for each $i=1,2$;   
\item \label{cp3}
$\bigcup _{j=1}^{3}h_{j}^{-1}(A)\subset A$, where 
$A=\{ z\in \CC \mid \frac{1}{2}\leq |z|\leq 3\}$; and 
\item \label{cp4}
$h_{3}^{-1}(A)\cap (\bigcup _{j=1}^{2}h_{i}^{-1}(A))=
\emptyset .$

\end{enumerate}
Combining statements \ref{cp3} and \ref{cp4} above,  
Lemma~\ref{hmslem}-\ref{backmin}, and Lemma~\ref{hmslem}-\ref{bss}, 
we get that for each $(h_{1},h_{2},h_{3})\in V$, 
$J(\langle h_{1},h_{2},h_{3}\rangle )\subset A$ and 
$J(\langle h_{1},h_{2},h_{3}\rangle )$ is disconnected.
Hence, for each $(h_{1},h_{2},h_{3})\in V$, we have 
$\langle h_{1},h_{2},h_{3}\rangle \in {\cal G}_{dis}.$ 
Combining it with statements \ref{cp1} and 
\ref{cp2} above and Theorem~\ref{countthm}, it follows that 
$J(h_{1})\cup J(h_{2})\subset J_{0}$ for some 
$J_{0}\in \hat{{\cal J}}_{\langle h_{1},h_{2},h_{3}\rangle }$, 
$h_{j}^{-1}(J(h_{3}))\cap J_{0}\neq \emptyset $ for each 
$j=1,2,$ 
and 
$\sharp (\hat{{\cal J}}_{\langle h_{1},h_{2},h_{3}\rangle })
=\aleph _{0}$, for each 
$(h_{1},h_{2},h_{3})\in V.$ 
Since $J(h_{1})<J(h_{3})$, 
Theorem~\ref{countthm} implies that 
the connected component $J_{0}$ should be equal to 
$J_{\min }(\langle h_{1},h_{2},h_{3}\rangle )$, 
and that $J_{\max }(\langle h_{1},h_{2},h_{3}\rangle )=J(h_{3}).$ 
 
 Thus, we have proved Proposition~\ref{countprop}.
\qed 

\ 

 We now show Proposition~\ref{countcomp}.

\noindent {\bf Proof of 
Proposition~\ref{countcomp}:} 
In fact, we show the following claim:\\ 
Claim: There exists a polynomial semigroup 
$G=\langle h_{1},h_{2},h_{3}\rangle $ 
in ${\cal G}$ such that 
all of the following hold.
\begin{enumerate}
\item $\sharp (\hat{{\cal J}})=\aleph _{0} .$
\item $J_{\min }\supset J(h_{1})\cup  
J(h_{2})$ and 
there exists a superattracting fixed point 
$z_{0}$ of $h_{1}$ with $z_{0}\in $ int$(J_{\min }).$ 
\item $J_{\max }=J(h_{3}).$
\item There exists a sequence $\{ n_{j}\} _{j\in \NN }$ of 
positive integers such that 
$\hat{{\cal J}}=
\{ J_{\min }\} \cup \{ J_{j}\mid j\in \NN \} $, 
where $J_{j}$ denotes the element of $\hat{{\cal J}}$ with 
$h_{3}^{-n_{j}}(J_{\min })\subset J_{j}.$ 
\item For any $J\in \hat{{\cal J}}$ with $J\neq J_{\max }$, 
there exists no sequence $\{ C_{j}\} _{j\in \NN }$ 
of mutually distinct elements of $\hat{{\cal J}}$ 
such that $\min _{z\in C_{j}}d(z,J)\rightarrow 0$ as 
$j\rightarrow \infty .$   
\item $G$ is sub-hyperbolic: i.e.,  
$\sharp (P(G)\cap J(G))<\infty  $ and 
$P(G)\cap F(G)$ is compact.
\end{enumerate} 
To show the claim, 
let $g_{1}(z)$ be the second iterate of $z\mapsto z^{2}-1.$ 
Let $g_{2}$ be a polynomial 
such that $J(g_{2})=\{ z\mid |z|=1\} $ and 
$g_{2}(-1)=-1.$ Then, we have 
$g_{1}(\sqrt{-1})=3\in \CCI \setminus K(g_{1}).$ 
Take a large, positive integer $m_{1}$, and 
let $a:=g_{1}^{m_{1}}(\sqrt{-1}).$ 
Then, 
\begin{equation}
\label{countcompeq1}
J(\langle g_{1}^{m_{1}},g_{2}\rangle )\subset 
\{ z\mid |z|<a\} .
\end{equation}  
Furthermore, 
since   
$a>\frac{1}{2}+\frac{\sqrt{5}}{2}$, we have   
\begin{equation}
\label{countcompeq1-5}
\overline{(g_{1}^{m_{1}})^{-1}(\{ z\mid |z|<a\} )}\subset 
\{ z\mid |z|<a\} .
\end{equation}
Let $g_{3}$ be a polynomial such that 
$J(g_{3})=\{ z\mid |z|=a\} .$ 
Since $-1$ is a superattracting fixed point 
of $g_{1}^{m_{1}}$ and it belongs to $J(g_{2})$, 
by \cite[Theorem 4.1]{HM2}, we see that 
for any $m\in \NN $, 
\begin{equation}
\label{countcompeq2}
-1\in \mbox{int}(J(\langle g_{1}^{m_{1}} , g_{2}^{m}\rangle )).
\end{equation}
Since $J(g_{2})\cap \mbox{int}(K(g_{1}^{m_{1}}))\neq \emptyset $ and 
$J(g_{2})\cap F_{\infty }(g_{1}^{m_{1}})\neq \emptyset $, 
we can take an $m_{2}\in \NN $ such that 
\begin{equation}
\label{countcompeq3}
(g_{2}^{m_{2}})^{-1}(\{ z\mid |z|=a\} )
\cap J(\langle g_{1}^{m_{1}},g_{2}^{m_{2}}\rangle )\neq \emptyset 
\end{equation} 
and 
\begin{equation}
\label{countcompeq4}
\overline{(g_{2}^{m_{2}})^{-1}(\{ z\mid |z|<a\} )}\subset 
\{ z\mid |z|<a\} .
\end{equation}
Take a small $r>0$ such that
\begin{equation}
\label{countcompeq5} 
\mbox{ for each }j=1,2,3,\ g_{j}(\{ z\mid |z|\leq r\})\subset 
\{ z\mid |z|<r\} . 
\end{equation}
Take an $m_{3}$ such that 
\begin{equation}
\label{countcompeq6}
(g_{3}^{m_{3}})^{-1}(\{ z\mid |z|=r\} )
\cap (\bigcup _{j=1}^{2}(g_{j}^{m_{j}})^{-1}(\{ z\mid |z|\leq a\} )) 
=\emptyset 
\end{equation}
and 
\begin{equation}
\label{countcompeq4-5}
g_{3}^{m_{3}}(-1)\in \{ z\mid |z|<r\} .
\end{equation}
Let $K:=\{ z\mid r\leq |z|\leq a\} .$ 
Then, by (\ref{countcompeq1-5}), (\ref{countcompeq4}), 
 (\ref{countcompeq5}) and (\ref{countcompeq6}), we have 
\begin{equation}
\label{countcompeq7}
(g_{j}^{m_{j}})^{-1}(K)\subset K, \mbox{ for } 
j=1,2,3, \mbox{and }
(g_{3}^{m_{3}})^{-1}(K)\cap 
(\bigcup _{j=1}^{2}(g_{j}^{m_{j}})^{-1}(K))= \emptyset .
\end{equation}     
Let $h_{j}:=g_{j}^{m_{j}}$, for each $j=1,2,3$, 
and let $G=\langle h_{1},h_{2},h_{3}\rangle .$ 
Then, by (\ref{countcompeq7}) and Lemma~\ref{hmslem}-\ref{backmin}, 
we obtain: 
\begin{equation}
\label{countcompeq8}
J(G)\subset K \mbox{ and }
h_{3}^{-1}(J(G))\cap (\bigcup _{j=1}^{2}
h_{j}^{-1}(J(G)))=\emptyset .
\end{equation}
Combining it with Lemma~\ref{hmslem}-\ref{bss}, 
it follows that 
$J(G)$ is disconnected. 
Furthermore, combining (\ref{countcompeq5}) and 
(\ref{countcompeq4-5}),  
we see $G\in {\cal G}$,
$P(G)\cap J(G)=\{ -1\} $, and 
that $P(G)\cap F(G)$ is compact. 
By Proposition~\ref{bminprop}, there exists 
a $j\in \{ 1,2,3\} $ with 
$J(h_{j})\subset J_{\min }.$ 
Since $ J(G)\subset K\subset \{ z\mid |z|\leq a\} $ and 
$J(h_{3})=\{ z\mid |z|=a\} $, 
we have 
\begin{equation}
\label{countcompeqmax}
J(h_{3})\subset J_{\max }.
\end{equation} 
Hence, either $J(h_{1})\subset J_{\min }$ or 
$J(h_{2})\subset J_{\min }.$ Since 
$J(h_{1})\cup J(h_{2})$ is connected, 
it follows that 
\begin{equation}
\label{countcompeqmin}
J(h_{1})\cup J(h_{2})\subset J_{\min }.
\end{equation} 
Combining this with Theorem~\ref{mainth1}-\ref{mainth1-3},
we have $h_{j}^{-1}(J_{\min })\subset J_{\min }$, for each 
$j=1,2.$ Hence, 
\begin{equation}
\label{countcompeq8-5}
J(\langle h_{1},h_{2}\rangle )\subset  
J_{\min }.
\end{equation}
Since $\sqrt{-1}\in J(h_{2})$ and $h_{1}(\sqrt{-1})=a\in J(h_{3})$,
 we obtain 
\begin{equation}
\label{countcompeqin}
h_{1}^{-1}(J(h_{3}))\cap J_{\min }\neq 
 \emptyset .
\end{equation} 
  
Similarly, by (\ref{countcompeq3}) and 
(\ref{countcompeq8-5}), 
we obtain  
\begin{equation}
\label{countcompeqin2}
h_{2}^{-1}(J(h_{3}))\cap J_{\min }\neq \emptyset .
\end{equation} 
Combining (\ref{countcompeqmax}), 
(\ref{countcompeqin}), (\ref{countcompeqin2}), and Theorem~\ref{countthm}, 
we obtain  
$\sharp (\hat{{\cal J}})=\aleph _{0}$, 
$J_{\max }=J(h_{3})$, 
$J(G)=J_{\max }\cup \bigcup _{n\in \NN \cup \{ 0\} }
h_{3}^{-n}(J_{\min })$, and 
that for any $J\in \hat{{\cal J}}$ with 
$J\neq J_{\max}$, 
there exists no sequence $\{C _{j}\} _{j\in \NN }$ 
of mutually distinct elements of 
$\hat{{\cal J}}$ such that 
$\min _{z\in C_{j}}d(z,J)\rightarrow 0$ as $j\rightarrow \infty .$ 
 
Moreover, by (\ref{countcompeq2}) and (\ref{countcompeq8-5})
(or by Theorem~\ref{mainupthm}-\ref{mainupthm3}), 
the superattracting fixed point $-1$ of $h_{1}$ belongs to 
int$(J_{\min }).$ 

 Hence, we have shown the claim.

Therefore, we have proved Proposition~\ref{countcomp}.
\qed

\subsection{Proofs of results in \ref{Hyperbolicity}}
In this section, we prove the results in section~\ref{Hyperbolicity}.  

We now demonstrate Proposition~\ref{nonminnoncpt}. \\ 
{\bf Proof of Proposition~\ref{nonminnoncpt}:} 
Since $\G \setminus \G _{\min }$ 
is not compact, 
there exists a sequence $\{ h_{j}\} _{j\in \NN }$ 
in $\G \setminus \G _{\min }$ and an element $h_{\infty }\in \G _{\min }$ 
such that $h_{j}\rightarrow h_{\infty }$ as $j\rightarrow \infty .$ 
By Theorem~\ref{mainth2}-\ref{mainth2-4-2}, 
for each $j\in \NN $, 
$h_{j}(K(h_{\infty }))$ is included in a connected component 
$U_{j}$ of int$(\hat{K}(G)).$ 
Let $z_{1}\in $ int$(\hat{K}(G))$ $(\subset $ int$(K(h_{\infty })))$ 
be a point. Then, 
$h_{\infty }(z_{1})\in $ int$(\hat{K}(G))$ and 
$h_{j}(z_{1})\rightarrow h_{\infty }(z_{1})$ as 
$j\rightarrow \infty .$  
 Hence, we may assume that there exists a connected 
 component $U$ of int$(\hat{K}(G))$ such that 
 for each $j\in \NN $, 
 $h_{j}(K(h_{\infty }))\subset U.$ 
 Therefore, $K(h_{\infty })=h_{\infty }(K(h_{\infty }))
 \subset \overline{U}.$ 
Since $\overline{U}\subset K(h_{\infty })$, 
we obtain $K(h_{\infty })=\overline{U}.$ 
Since $U\subset $ int$(K(h_{\infty }))\subset \overline{U}$ and 
$U$ is connected, it follows that 
int$(K(h_{\infty }))$ is connected. 
Moreover, we have 
$U\subset $ int$(K(h_{\infty }))\subset $ int$(\overline{U})\subset 
$ int$(\hat{K}(G)).$  
Thus, 
\begin{equation}
\label{nonminnoncptpfeq1}
\mbox{ int}(K(h_{\infty }))=U.
\end{equation} 
Furthermore, since
\begin{equation}
\label{nonminnoncptpfeq2}
J(h_{\infty })<J(h_{j}) \mbox{ for each } j\in \NN ,
\end{equation}
and $h_{j}\rightarrow h_{\infty }$ as $j\rightarrow \infty $, 
we obtain 
\begin{equation}
\label{nonminnoncptpfeq3}
J(h_{j})\rightarrow J(h_{\infty }) \mbox{ as } j\rightarrow \infty ,
\end{equation} 
with respect to the Hausdorff metric.
Combining that $h_{j}\in \G \setminus \G _{\min }$ for each $j\in \NN $ 
with Theorem~\ref{mainth2}-\ref{mainth2ast1}, 
(\ref{nonminnoncptpfeq1}), (\ref{nonminnoncptpfeq2}), 
and (\ref{nonminnoncptpfeq3}), we see that 
for each $h\in \G _{\min }$, $K(h)=K(h_{\infty }).$ 
Combining it with (\ref{nonminnoncptpfeq1}), 
(\ref{nonminnoncptpfeq2}) and (\ref{nonminnoncptpfeq3}),  
it follows that statement~\ref{nonminnoncpt1} in Proposition~\ref{nonminnoncpt} holds. 
To prove statement~\ref{nonminnoncpt2}, let $h\in \Gamma _{\min }.$ 
Aplying the Riemann-Hurwitz formula to 
$h: \mbox{int}(K(h))\rightarrow \mbox{int}(K(h))$, 
we obtain that each finite critical point of $h$ belongs to 
int$(K(h)).$ If $h$ is hyperbolic, then by using quasiconformal surgery (\cite{Be1}), 
we can see that statement~\ref{nonminnoncpt2-1} holds. If $h$ is not hyperbolic, then 
statement~\ref{nonminnoncpt2-2} holds.  

Thus we have proved Proposition~\ref{nonminnoncpt}. 
\qed 

To demonstrate Theorem~\ref{shshprop}, 
we need the following.

\begin{lem}
\label{ksetmin}
Let $G$ be a  
polynomial semigroup generated by 
a non-empty compact set $\G $ in 
{\em Poly}$_{\deg \geq 2}.$ 
Suppose that $G\in {\cal G}_{dis}.$ 
Then, we have $\hat{K}(G_{\min ,\G })=\hat{K}(G).$ 
\end{lem}
\begin{proof}

Since $G_{\min ,\G }\subset G$, 
we have $\hat{K}(G)\subset \hat{K}(G_{\min ,\G }).$ 
Moreover, it is easy to see 
$\hat{K}(G_{\min ,\G })=\bigcap _{g\in G_{\min ,\G }}K(g).$ 
Let $g\in G_{\min ,\G }$ and 
$h\in \G \setminus \G _{\min }.$ 
For each $\alpha \in \G _{\min }$, 
we have $\alpha ^{-1}(J_{\min }(G))\subset J_{\min }(G).$ 
Since $\sharp (J_{\min }(G))\geq 3 $ (Theorem~\ref{mainth2}-\ref{mainth2-4-1}), 
Lemma~\ref{hmslem}-\ref{backmin} implies that 
$J(g)\subset J_{\min }(G).$ Hence, 
from Theorem~\ref{mainth2}-\ref{mainth2-4-2}, 
it follows that 
\begin{equation}
\label{ksetminpfeq1}
h(J(g))\subset \mbox{int}(\hat{K}(G))\subset \mbox{int}(\hat{K}(g)).
\end{equation}
Since $J(g)$ is connected and each connected component of int$(K(g))$ 
is simply connected, the above (\ref{ksetminpfeq1}) implies that 
$h(K(g))\subset K(g).$ 
Hence, we obtain 
$h(\hat{K}(G_{\min ,\G }))=
h(\bigcap _{g\in G_{\min ,\G }}K(g))
\subset \bigcap _{g\in G_{\min ,\G }}K(g)
=\hat{K}(G_{\min ,\G }).$
Combined with that $\alpha (\hat{K}(G_{\min ,\G }))\subset 
\hat{K}(G_{\min ,\G })$ for each $\alpha \in \G_{\min }$, 
it follows that for each $\beta \in G$, 
$\beta (\hat{K}(G_{\min ,\G }))\subset \hat{K}(G_{\min ,\G }).$ 
Therefore, we obtain $\hat{K}(G_{\min ,\G })\subset 
\hat{K}(G).$ Thus, it follows that  
 $\hat{K}(G_{\min ,\G })=\hat{K}(G).$
\end{proof}
\begin{df}
Let $G$ be a rational semigroup and $N$ a positive integer. 
We denote by $SH_{N}(G)$ the set of points $z\in \CCI $ 
satisfying that  
there exists a positive number $\delta $  
such that for each $g\in G$, 
$\deg (g:V\rightarrow B(z,\delta ))\leq N$, 
for each connected component $V$ of $g^{-1}(B(z,\delta )).$ 
Moreover, 
we set $UH(G):= \CCI \setminus \bigcup _{N\in \NN }SH_{N}(G).$ 

\end{df}
\begin{lem}
\label{060523lem1}
Let $G$ be a polynomial semigroup generated by a compact subset 
$\G $ of {\em Poly}$_{\deg \geq 2}.$ 
Suppose that $G\in {\cal G}_{dis }$ and that 
$\G \setminus \G _{\min }$ is not compact. 
Moreover, suppose that (a) in 
Proposition~\ref{nonminnoncpt}-\ref{nonminnoncpt2} holds.
Then, there exists an open neighborhood ${\cal U} $ of 
$\G _{\min }$ in $\G $ and an open set $U$ in {\em int}$(\hat{K}(G))$ 
with $\overline{U}\subset $ {\em int}$(\hat{K}(G))$ such that:
\begin{enumerate}
\item $\bigcup _{h\in {\cal U}}h(U)\subset U$; 
\item $\bigcup _{h\in {\cal U}}CV^{\ast }(h) \subset U$, 
and 
\item 
denoting by $H$ the polynomial semigroup generated by ${\cal U}$, 
we have that $P^{\ast }(H)\subset $ {\em int}$(\hat{K}(G))\subset F(H)$
and that $H$ is hyperbolic. 
\end{enumerate}  
\end{lem}
\begin{proof}
Let $h_{0}\in \G _{\min }$ be an element.  
Let 
${\cal E}:= 
\{ \psi (z)=az+b\mid a,b\in \CC, |a|=1, \psi (J(h_{0}))=J(h_{0})\} .$ 
Then, by \cite{Be}, ${\cal E}$ is compact in Poly. 
Moreover, by \cite{Be}, we have the following two claims:\\ 
Claim 1: If $J(h_{0})$ is a round circle with the center $b_{0}$ and 
radius $r$, then 
${\cal E}=\{ \psi (z)=a(z-b_{0})+b_{0}\mid |a|=r\} .$ \\ 
Claim 2: If $J(h_{0})$ is not a round circle, then 
$\sharp {\cal E}<\infty .$ 

 Let $z_{0}$ be the unique attracting fixed point of $h_{0}$ in $\CC .$ 
 Let $g\in G_{\min ,\G }.$ 
 By \cite{Be}, for each $n\in \NN $, 
 there exists an $\psi _{n}\in {\cal E}$ such that 
 $h_{0}^{n}g=\psi _{n}gh_{0}^{n}.$ Hence, 
 for each $n\in \NN $, 
 $h_{0}^{n}g(z_{0})=\psi _{n}gh_{0}^{n}(z_{0})=\psi _{n}g(z_{0}).$ 
 Combining it with Claim 1 and Claim 2, 
 it follows that there exists an $n\in \NN $ such that 
 $h_{0}^{n}(g(z_{0}))=z_{0}.$ For this $n$, 
 $g(z_{0})=\psi _{n}^{-1}(h_{0}^{n}(g(z_{0})))
 =\psi _{n}^{-1}(z_{0})\in 
 \bigcup _{\psi \in {\cal E}}\psi (z_{0}).$ 
 Combining it with Claim 1 and Claim 2 again, 
 we see that the set $C:= \overline{\bigcup _{g\in G_{\min ,\G}}\{ g(z_{0})\} } $ 
 is a compact subset of int$(\hat{K}(G)).$ 
Let $d_{H}$ be the hyperbolic distance on int$(\hat{K}(G)).$ 
Let $R>0$ be a large number such that 
setting $U:=\{ z\in \mbox{int}(\hat{K}(G))\mid 
\min _{a\in C}d_{H}(z,a)<R\} $, 
we have $\bigcup _{h\in \G _{\min }}CV^{\ast }(h)\subset U .$ 
 Then, for each $h\in \G _{\min }$, 
 $\overline{h(U)}\subset U.$ 
 Therefore, there exists an open neighborhood ${\cal U}$ of 
 $\G _{\min }$ in $\G $ such that 
 $\bigcup _{h\in {\cal U}}h(U)\subset U$, and such that 
$\bigcup _{h\in {\cal U}}CV^{\ast }(h) \subset U.$ 
Let $H$ be the polynomial semigroup generated by ${\cal U}.$ 
From the above argument, we obtain 
$P^{\ast }(H)=\overline{\bigcup _{g\in H}CV^{\ast }(g) }
\subset $ $\overline{\bigcup _{g\in H\cup \{ Id\} }g\left( \bigcup _{h\in {\cal U}}CV^{\ast }(h) \right) }   
\subset \overline{\bigcup _{g\in H\cup \{ Id\} }g(U)}\subset \overline{U}\subset 
$ int$(\hat{K}(G))\subset F(H).$ 
Hence, $H$ is hyperbolic. Thus, we have proved 
Lemma~\ref{060523lem1}. 
\end{proof}

We now demonstrate Theorem~\ref{shshprop}.\\ 
{\bf Proof of Theorem~\ref{shshprop}:}
Suppose that $G_{\min ,\G }$ is semi-hyperbolic.
We will consider the following two cases:\\ 
Case 1: $\G \setminus \G _{\min }$ is compact.\\ 
Case 2: $\G \setminus \G _{\min }$ is not compact.

 Suppose that we have Case 1.
Since $UH(G_{\min ,\G })\subset P(G_{\min ,\G })$, 
$G_{\min ,\G }\in {\cal G}$, and $G_{\min ,\G }$ is semi-hyperbolic, 
we obtain 
$UH(G_{\min ,\G })\cap \CC 
\subset F(G _{\min ,\G })\cap \hat{K}(G_{\min ,\G })$ 
=int$(\hat{K}(G_{\min ,\G })).$  
By Lemma~\ref{ksetmin}, we have 
$\hat{K}(G_{\min ,\G })=\hat{K}(G).$ 
Hence, we obtain 
\begin{equation}
\label{mainth3-2pfeq1}
UH(G_{\min ,\G })\cap \CC \subset  
\mbox{ int}(\hat{K}(G))\subset \CC \setminus J_{\min }(G).
\end{equation}
Therefore, there exists a positive integer $N$ and a 
positive number $\delta $ 
such that for each $z\in J_{\min }(G)$ and 
each $h\in G_{\min ,\G }$, 
we have 
\begin{equation}
\label{mainth3-2pfeq2}
\deg (h:V\rightarrow D(z,\delta ))\leq N,  
\end{equation}
for each connected component $V$ of $h^{-1}(D(z,\delta )).$ 
Moreover, combining Theorem~\ref{mainth2}-\ref{mainth2-4-2} and 
Theorem~\ref{mainth2}-\ref{mainth2-3}, we obtain  
$\bigcup _{\alpha \in \G \setminus \G_{\min }}
\alpha ^{-1}(J_{\min }(G))\cap P^{\ast }(G)=\emptyset .$ 
Hence, there exists a number $\delta _{1}$ such that 
for each $z\in \bigcup _{\alpha \in \G \setminus \G_{\min }}
\alpha ^{-1}(J_{\min }(G))$ and each $\beta \in G\cup \{ Id \} $, 
\begin{equation}
\label{mainth3-2pfeq3}
\deg (\beta :W\rightarrow D(z,\delta _{1}))=1,
\end{equation}
for each connected component $W$ of $\beta ^{-1}(D(z,\delta _{1})).$ 
For this $\delta _{1}$, there exists a number $\delta _{2}>0$ 
such that for each $z\in J_{\min }(G)$ and  each $\alpha \in \G 
\setminus \G _{\min }$, 
\begin{equation}
\label{mainth3-2pfeq4}
\mbox{diam }B\leq \delta _{1},\  
\deg (\alpha : B\rightarrow D(z,\delta _{2}))
\leq \max\{ \deg (\alpha )\mid \alpha \in \G \setminus \G _{\min }\} 
\end{equation}
for each connected component $B$ of $\alpha ^{-1}(D(z,\delta _{2})).$
Furthermore, 
by \cite[Lemma 1.10]{S1} (or \cite{S2}), 
we have that there exists a constant $0<c<1$ such that 
for each $z\in J_{\min }(G)$, each $h\in G_{\min ,\G }\cup \{ Id \} $, 
and each connected component $V$ of $h^{-1}(D(z,c\delta ))$, 
\begin{equation}
\label{mainth3-2pfeq5}
 \mbox{diam }V\leq \delta _{2}.
\end{equation}
Let $g\in G$ be any element. 

 Suppose that $g\in G_{\min ,\G}.$ Then, 
 by (\ref{mainth3-2pfeq2}), for each $z\in J_{\min }(G)$,  we have 
$\deg (g:V\rightarrow D(z, c\delta  ))\leq N$, for each 
connected component $V$ of $g^{-1}(D(z,c\delta )).$  

 Suppose that 
 $g$ is of the form $g=h\circ \alpha \circ g_{0}$, where 
$h\in G_{\min ,\G }\cup \{ Id \}$, $ \alpha \in \G \setminus \G _{\min }$, 
and $g_{0}\in G\cup \{ Id \} .$ 
 Then, 
 combining (\ref{mainth3-2pfeq3}), (\ref{mainth3-2pfeq4}), 
 and (\ref{mainth3-2pfeq5}),  
 we get that for each $z\in J_{\min }(G)$, 
 $\deg (g: W\rightarrow D(z,c\delta ))\leq N\cdot \max 
 \{ \deg (\alpha )\mid \alpha \in \G \setminus \G _{\min }\} $, 
 for each connected component $W$ of $g^{-1}(D(z,c\delta )).$  

 From the above argument, we see that 
 $J_{\min }(G)\subset  SH_{N'}(G)$, 
 where $N':=N\cdot \max \{\deg( \alpha )\mid 
 \alpha \in \G \setminus \G _{\min }\} .$ 
Moreover, by Theorem~\ref{mainth2}-\ref{mainth2-3},
we see that for any point $z\in J(G)\setminus J_{\min }(G)$, 
$z\in SH_{1}(G).$ 
Hence, we have shown that 
$J(G)\subset \CCI \setminus UH(G).$ Therefore, 
$G$ is semi-hyperbolic, provided that we have Case 1.  

 We now suppose that we have Case 2. 
 Then, by Proposition~\ref{nonminnoncpt}, 
 we have that for each $h\in \G _{\min }$, 
 $K(h)=\hat{K}(G)$ and int$(K(h))$ is non-empty and connected. 
 Moreover, for each $h\in \G _{\min }$, 
 int$(K(h))$ is an immediate basin of an attracting fixed point 
 $z_{h}\in \CC .$ 
Let ${\cal U}$ be the open neighborhood of 
$\G _{\min }$ in $\G $ as in 
Lemma~\ref{060523lem1}. 
Denoting by $H$ the polynomial semigroup generated by 
${\cal U}$, we have 
$P^{\ast }(H)\subset $ int$(\hat{K}(G)).$ 
Therefore, there exists a number $\delta >0$ such that 
\begin{equation}
\label{shshproppfeq02}
D(J(G),\delta )\subset \CC \setminus P(H).
\end{equation} 
Moreover, combining Theorem~\ref{mainth2}-\ref{mainth2-4-2} 
and that $\G \setminus {\cal U}$ is compact, 
we see that there exists a number 
$\epsilon >0$ such that 
\begin{equation}
\label{shshproppfeq03}
\overline{ \bigcup_{\alpha \in \G \setminus {\cal U}}
\alpha ^{-1}(D(J_{\min }(G),\epsilon ))}\subset A_{0},  
\end{equation}
where $A_{0}$ denotes the unbounded component of $\CC \setminus J_{\min }(G).$ 
Combining it with Theorem~\ref{mainth2}-\ref{mainth2-3}, 
it follows that there exists a number 
$\delta _{1}>0$ such that 
\begin{equation}
\label{shshproppfeq04}
D\left( \bigcup _{\alpha \in 
\G \setminus {\cal U}}\alpha ^{-1}
(D(J_{\min }(G),\epsilon )),\ \delta _{1}\right) 
\subset \CC \setminus P(G).
\end{equation}
For this $\delta _{1}$, 
there exists a number $\delta _{2}>0$ such that 
for each $\alpha \in \G \setminus {\cal U}$ and each 
$x\in D(J_{\min }(G),\epsilon )$, 
\begin{equation}
\label{shshproppfeq05}
\mbox{diam }B\leq \delta _{1},\ 
\deg (\alpha :B\rightarrow D(x,\delta _{2}))\leq 
\max \{ \deg (\beta )\mid \beta \in \G \setminus {\cal U}\} 
\end{equation}
for each connected component $B$ of 
$\alpha ^{-1}(D(x,\delta _{2})).$ 
By Lemma~\ref{invnormal2} and (\ref{shshproppfeq02}), 
there exists a constant $c>0$ such that 
for each $h\in H$ and each $z\in J_{\min }(G)$, 
\begin{equation}
\label{shshproppfeq06}
\mbox{diam }V\leq \min \{ \delta _{2},\epsilon \} , 
\end{equation}
for each connected component $V$ of $h^{-1}(D(z,c\delta )).$ 
Let $z\in J_{\min }(G)$ and $g\in G.$ We will show that 
$z\in \CC \setminus UH(G).$ 

 Suppose that $g\in H.$ Then, 
 (\ref{shshproppfeq02}) implies that 
 for each 
connected component $V$ of 
$g^{-1}(D(z,c\delta ))$, 
$\deg (g:V\rightarrow D(z,c\delta ))=1.$ 

 Suppose that $g$ is of the form $g=h\circ \alpha \circ 
 g_{0}$, where 
$h\in H\cup \{ Id\} , \alpha \in \G \setminus 
{\cal U}, g_{0}\in G\cup \{ Id \} .$ 
Let 
$W$ be a connected component of $g^{-1}(D(z,c\delta ))$ and 
let $W_{1}:=g_{0}(W)$ and $V:=\alpha (W_{1}).$ 
Let $z_{1}$ be the point such that 
$\{ z_{1}\} =V\cap h^{-1}(\{ z\} ).$ 
If $z_{1}\in \CC \setminus 
D(J_{\min }(G), \epsilon )$, then, 
by (\ref{shshproppfeq06}) and Theorem~\ref{mainth2}-\ref{mainth2-3}, 
$V\subset D(z_{1},\epsilon )\subset \CC \setminus P(G).$ 
 Hence, 
 $\deg (\alpha \circ g_{0}:W\rightarrow V)=1$, 
 which implies that 
 $\deg (g:W\rightarrow D(z,c\delta ))=1.$ 
If $z_{1}\in D(J_{\min }(G), \epsilon )$, 
then by (\ref{shshproppfeq06}), 
$V\subset D(z_{1},\delta _{2}).$ 
Combining it with 
(\ref{shshproppfeq04}) and (\ref{shshproppfeq05}), 
we obtain  
$\deg (\alpha \circ g_{0}:W\rightarrow V)=
\deg (\alpha :W_{1}\rightarrow V)\leq 
\max \{ \deg (\beta )\mid \beta \in \G \setminus {\cal U}\} .$ 
Therefore, 
$\deg (g:W\rightarrow D(z,c\delta ))\leq 
\max \{ \deg (\beta )\mid \beta \in \G \setminus {\cal U}\} .$ 
Thus, $J_{\min }(G)\subset \CC \setminus UH(G).$ 

 Moreover, Theorem~\ref{mainth2}-\ref{mainth2-3} 
 implies that 
 $J(G)\setminus J_{\min }(G)\subset 
 \CC \setminus P(G)\subset \CC \setminus UH(G).$ 
 Therefore, $J(G)\subset \CC \setminus UH(G)$, 
 which implies that $G$ is semi-hyperbolic.  

 Thus, we have proved Theorem~\ref{shshprop}.
\qed 

\ 

We now demonstrate Theorem~\ref{hhprop}.\\ 
{\bf Proof of Theorem~\ref{hhprop}:} 
We use the same argument as that in the proof of Theorem~\ref{shshprop}, 
but we modify it as follows: 
\begin{enumerate}
\item 
 In (\ref{mainth3-2pfeq1}), we replace 
$UH(G_{\min ,\G })\cap \CC $ by $P^{\ast }(G_{\min ,\G }).$  
\item  In (\ref{mainth3-2pfeq2}), 
we replace $N$ by $1.$  
\item  We replace (\ref{mainth3-2pfeq4}) by 
the following (\ref{mainth3-2pfeq4})' 
$\mbox{diam }B\leq \delta _{1},\  
\deg (\alpha : B\rightarrow D(z,\delta _{2}))=1.$  
\item  We replace (\ref{shshproppfeq05}) by the following 
(\ref{shshproppfeq05})' 
$\mbox{diam }B\leq \delta _{1},\ 
\deg (\alpha :B\rightarrow D(x,\delta _{2}))=1.$ 
(We take the number $\epsilon >0$ so small.)  
\end{enumerate} 
With these modifications, it is easy to see that $G$ is hyperbolic. 

 Thus, we have proved Theorem~\ref{hhprop}. 
\qed 

We now prove Proposition~\ref{p:2.312as}. 

{\bf Proof of Proposition~\ref{p:2.312as}:} 
Combining Lemma~\ref{060523lem1} and Theorems~\ref{shshprop}, \ref{hhprop}, 
it is easy to see that Proposition~\ref{p:2.312as} holds.   
\subsection{Proofs of results in \ref{Const}}
In this section, we prove the results in \ref{Const}. 

We now demonstrate Proposition~\ref{Constprop}.\\ 
{\bf Proof of Proposition~\ref{Constprop}:}
Conjugating $G$ by $z\mapsto z+b$, we may assume that $b=0.$ 
For each $h\in \G$, we set $a_{h}:=a(h)$ 
and $d_{h}:=\deg (h).$ Let $r>0$ be a number 
such that $\overline{D(0,r)}\subset 
$ int$(\hat{K}(G)).$ 

Let $h\in \G $ and let $\alpha >0$ be a number. 
Since $d\geq 2$ and $(d,d_{h})\neq (2,2)$, 
it is easy to see that 
$(\frac{r}{\alpha })^{\frac{1}{d}}>
2\left(\frac{2}{|a_{h}|}(\frac{1}{\alpha })
^{\frac{1}{d-1}}\right)^{\frac{1}{d_{h}}}
$ if and only if 
\begin{equation}
\label{Contproppfeq1}
\log \alpha <
\frac{d(d-1)d_{h}}{d+d_{h}-d_{h}d}
( \log 2-\frac{1}{d_{h}}\log \frac{|a_{h}|}{2}-\frac{1}{d}\log r) .
\end{equation} 
We set 
\begin{equation}
\label{Contproppfeq2}
c_{0}:=\min _{h\in \G }\exp \left(\frac{d(d-1)d_{h}}{d+d_{h}-d_{h}d}
( \log 2-\frac{1}{d_{h}}\log \frac{|a_{h}|}{2}-\frac{1}{d}\log r) \right)
\in (0,\infty ).
\end{equation}
Let $0<c<c_{0}$ be a small number and let $a\in \CC $ 
be a number with $0<|a|<c.$ 
Let $g_{a}(z)=az^{d}.$ 
Then, we obtain $K(g_{a})=\{ z\in \CC \mid 
|z|\leq (\frac{1}{|a|})^{\frac{1}{d-1}}\} $ and 
$g_{a}^{-1}(\{ z\in \CC \mid  |z|=r\} )=
\{ z\in \CC \mid |z|=(\frac{r}{|a|})^{\frac{1}{d}}\} .$
Let 
$D_{a}:=\overline{D(0,2(\frac{1}{|a|})^{\frac{1}{d-1}})}.$ 
Since $h(z)=a_{h}z^{d_{h}}(1+o(1))\ (z\rightarrow \infty )$ uniformly on 
$\G $, 
it follows that if $c$ is small enough, then 
for any $a\in \CC $ with $0<|a|<c$ and for any $h\in \G $, 
$h^{-1}(D_{a})\subset 
\left\{ z\in \CC \mid 
|z|\leq 2\left( \frac{2}{|a_{h}|}(\frac{1}{|a|})^{\frac{1}{d-1}}\right) 
^{\frac{1}{d_{h}}}\right\} .$  
This implies that for each $h\in \G $, 
\begin{equation}
\label{Contproppfeq3}
h^{-1}(D_{a})\subset g_{a}^{-1}(\{ z\in \CC \mid |z|<r\} ).
\end{equation} 
Moreover, if $c$ is small enough, then for any $a\in \CC $ with 
$0<|a|<c$ and any $h\in \G $,  
\begin{equation}
\label{Contproppfeq4}
\hat{K}(G)\subset g_{a}^{-1}(\{ z\in \CC \mid |z|<r\} ),\ 
\overline{h(\CCI \setminus D_{a})}\subset 
\CCI \setminus D_{a}.
\end{equation}
Let $a\in \CC $ with $0<|a|<c.$  
By (\ref{Contproppfeq3}) and (\ref{Contproppfeq4}), 
there exists a compact neighborhood $V$ of $g_{a}$ in Poly$_{\deg \geq 2}$,  
such that 
\begin{equation}
\label{Contproppfeq5} 
\hat{K}(G)\cup \bigcup _{h\in \G }h^{-1}
(D_{a})
\subset 
\mbox{int}\left( \bigcap _{g\in V}
g^{-1}(\{ z\in \CC \mid |z|<r\} )\right), \mbox{ and } 
\end{equation}
\begin{equation}
\label{Contproppfeq5-a}
\bigcup _{h\in \G \cup V}
\overline{h(\CCI \setminus D_{a})}
\subset \CCI \setminus D_{a},
\end{equation}
which implies that 
\begin{equation}
\label{Contproppfeq5-1}
\mbox{int}(\hat{K}(G))\cup 
(\CCI \setminus D_{a})
\subset F(H_{\G ,V}),
\end{equation}
where $H_{\G ,V}$ denotes the polynomial semigroup generated by 
the family $\G \cup V.$   

 By (\ref{Contproppfeq5}), we obtain that for any non-empty subset 
 $V'$ of $V$, 
\begin{equation}
\label{Contproppfeq7}
\hat{K}(G)=\hat{K}(H_{\G, V'} ),
\end{equation} 
where $H_{\G ,V'}$ denotes the polynomial semigroup generated by 
the family $\G \cup V'.$
If the compact neighborhood $V$ of $g_{a}$ is so small, then  
\begin{equation}
\label{Contproppfeq8}
\bigcup _{g\in V} CV^{\ast }(g)  \subset  \mbox{int}(\hat{K}(G)).
\end{equation} 
Since $P^{\ast }(G)\subset \hat{K}(G)$, 
combining it with (\ref{Contproppfeq7}) and (\ref{Contproppfeq8}),
we get  
that for any non-empty subset $V'$ of $V$, 
$P^{\ast }(H_{\G ,V'} )\subset 
\hat{K}(H_{\G ,V'} ).$ 
 Therefore, for any non-empty subset $V'$ of $V$,  
$H_{\G ,V'} \in {\cal G}.$ 

We now show that for any non-empty subset $V'$ of $V$, $J(H_{\G,V'})$ is disconnected and 
$(\Gamma \cup V')_{\min }\subset \Gamma.$  
Let $$U:=\left(\mbox{int}(\bigcap _{g\in V}g^{-1}(\{ z\in \CC \mid |z|<r\} ))\right)
\setminus 
\bigcup _{h\in \G }h^{-1}(D_{a}).$$ 
Then, for any $h\in \G $, 
\begin{equation}
\label{Contproppfeq6}
h(U)\subset \CCI \setminus D_{a}.
\end{equation} 
Moreover, for any $g\in V$, $g(U)\subset 
$ int$(\hat{K}(G)).$ 
Combining it with (\ref{Contproppfeq5-1}), 
(\ref{Contproppfeq6}), and Lemma~\ref{hmslem}-\ref{bss}, 
it follows that $U\subset 
F(H_{\G ,V} ).$ 
If the neighborhood $V$ of $g_{a}$ is  so small, then 
there exists an annulus $A$ in $U$ such that for any $g\in V$,  
$A$ separates $J(g)$ and $\bigcup _{h\in \G }h^{-1}(J(g)).$  
Hence, it follows that for any non-empty subset $V'$ of $V$, 
the polynomial semigroup $H_{\G, V'}$ generated by 
the family $\G \cup V'$ satisfies that $J(H_{\G ,V'})$ is 
disconnected and $(\G \cup V')_{\min }\subset \G .$  

 We now suppose that in addition to the assumption, 
 $G$ is semi-hyperbolic. Let $V'$ be any non-empty subset of 
 $V .$ Since $(\G \cup \overline{V'})_{\min }\subset \G $, 
 Theorem~\ref{shshprop} implies that the above 
 $H_{\G ,V'}$ is semi-hyperbolic. 

 We now suppose that in addition to the assumption, 
$G$ is hyperbolic.  Let $V'$ be any non-empty subset of 
$V.$ By (\ref{Contproppfeq7}) and 
(\ref{Contproppfeq8}), 
we have 
\begin{equation}
\label{Contproppfeq9}
\bigcup _{g\in \Gamma \cup \overline{V'}}CV^{\ast }(g)\subset 
\mbox{int}(\hat{K}(H_{\G ,\overline{V'}})).
\end{equation}  
Since $(\G \cup \overline{V'})_{\min }\subset \G $, 
combining it with (\ref{Contproppfeq9}) and 
Theorem~\ref{hhprop}, 
we obtain that $H_{\G ,V'}$ is hyperbolic. 

 Thus, we have proved Proposition~\ref{Constprop}. 
\qed  

\ 

 We now demonstrate Theorem~\ref{shshfinprop}.\\ 
 {\bf Proof of Theorem~\ref{shshfinprop}:} 
First, we show \ref{shshfinprop1}. 
Let $r>0$ be a number such that 
$D(b_{j},2r)\subset \mbox{int}(K(h_{1}))$ for each 
$j=1,\ldots ,m.$ 
If we take $c>0$ so small, then 
for each $(a_{2},\ldots ,a_{m})\in \CC ^{m-1}$ 
such that $0<|a_{j}|<c$ for each $j=2,\ldots ,m$, 
setting $h_{j}(z)=a_{j}(z-b_{j})^{d_{j}}+b_{j}$ 
($j=2,\ldots ,m$), we have 
\begin{equation}
\label{shshfinpropeq1}
h_{j}(K(h_{1}))\subset D(b_{j},r)\subset 
\mbox{int}(K(h_{1}))\ (j=2,\ldots ,m). 
\end{equation} 
Hence, $K(h_{1})=\hat{K}(G)$, 
where $G=\langle h_{1},\ldots ,h_{m}\rangle .$ 
Moreover, by (\ref{shshfinpropeq1}), 
we have $P^{\ast }(G)\subset K(h_{1}).$ 
Hence, $G\in {\cal G}.$ 

 If $\langle h_{1}\rangle $ is semi-hyperbolic, 
then using the same method as that of Case 1 in the proof of 
Theorem~\ref{shshprop}, we obtain that $G$ is semi-hyperbolic. 

 We now suppose that $\langle h_{1}\rangle $ is hyperbolic. 
By (\ref{shshfinpropeq1}), we have 
$\bigcup _{j=2}^{m}CV^{\ast }(h_{j})\subset 
\mbox{int}(\hat{K}(G)).$ Combining it with 
the same method as that in the proof of Theorem~\ref{hhprop}, 
we obtain that $G$ is hyperbolic. 
Hence, we have proved statement \ref{shshfinprop1}. 

 We now show statement \ref{shshfinprop2}. 
Suppose we have case (i). 
We may assume $d_{m}\geq 3.$ 
Then, by statement \ref{shshfinprop1}, 
there exists an element $a>0$ such that 
setting $h_{j}(z)=a(z-b_{j})^{d_{j}}+b_{j}$ ($j=2,\ldots ,m-1$), 
$G_{0}=\langle h_{1},\ldots ,h_{m-1}\rangle $ satisfies 
that $G_{0}\in {\cal G}$ and $\hat{K}(G_{0})=$ $K(h_{1})$ 
and if $\langle h_{1}\rangle $ is semi-hyperbolic (resp. hyperbolic), 
then $G_{0}$ is semi-hyperbolic (resp. hyperbolic). 
Combining it with Proposition~\ref{Constprop}, 
it follows that there exists an $a_{m}>0$ such that 
setting $h_{m}(z)=a_{m}(z-b_{m})^{d_{m}}+b_{m}$, 
$G=\langle h_{1},\ldots ,h_{m}\rangle $ satisfies that 
$G\in {\cal G}_{dis}$ and $\hat{K}(G)=\hat{K}(G_{0})=K(h_{1})$ and if 
$G_{0}$ is semi-hyperbolic (resp. hyperbolic), 
then $G$ is semi-hyperbolic (resp. hyperbolic).  

 Suppose now we have case (ii). 
 Then by Proposition~\ref{Constprop}, 
 there exists an $a_{2}>0$ such that 
setting $h_{j}(z)=a_{2}(z-b_{j})^{2}+b_{j}$ $(j=2,\ldots ,m)$, 
$G=\langle h_{1},\ldots ,h_{m}\rangle =\langle h_{1}, h_{2}\rangle $ 
satisfies that $G\in {\cal G}_{dis}$ and $\hat{K}(G)=K(h_{1})$ and if $\langle h_{1}\rangle $ 
is semi-hyperbolic (resp. hyperbolic), then 
$G$ is semi-hyperbolic (resp. hyperbolic). 

 Thus, we have proved Theorem~\ref{shshfinprop}.  
\qed  

We now demonstrate Theorem~\ref{sphypopen}.\\ 
{\bf Proof of Theorem~\ref{sphypopen}:} 
Statements \ref{sphypopen2} and \ref{sphypopen3} follow 
from Theorem~\ref{shshfinprop}. 

 We now show statement \ref{sphypopen1}. 
By \cite[Theorem 2.4.1]{S5}, 
${\cal H}_{m}$ and ${\cal H}_{m}\cap {\cal D}_{m}$ are open. 

 We now show that ${\cal H}_{m}\cap {\cal B}_{m}$ is open. 
In order to do that, 
let 
$(h_{1},\ldots ,h_{m})\in {\cal H}_{m}\cap {\cal B}_{m}.$ 
Let $\epsilon >0$ such that 
$D(P^{\ast }(\langle h_{1},\ldots ,h_{m}\rangle ),\ 3\epsilon )
\subset F(\langle h_{1},\ldots ,h_{m}\rangle ).$ 
By \cite[Theorem 1.35]{S1}, 
there exists an $n\in \NN $ such that 
for each $(i_{1},\ldots ,i_{n})\in 
\{ 1,\ldots ,m\} ^{n}$, 
$$h_{i_{n}}\cdots h_{i_{1}}
(D(P^{\ast }(\langle h_{1},\ldots ,h_{m}\rangle ),\ 2\epsilon ))
\subset D(P^{\ast }(\langle h_{1},\ldots ,h_{m}\rangle ),\ \epsilon /2).$$ 
Hence, there exists a neighborhood $U$ of $(h_{1},\ldots ,h_{m})$ 
in (Poly$_{\deg \geq 2})^{m}$  such that for each 
$(g_{1},\ldots ,g_{m})\in U$ and each 
$(i_{1},\ldots ,i_{n})\in \{ 1,\ldots ,m\} ^{n}$, 
$$g_{i_{n}}\cdots g_{i_{1}}
(D(P^{\ast }(\langle h_{1},\ldots ,h_{m}\rangle ),\ 2\epsilon ))
\subset D(P^{\ast }(\langle h_{1},\ldots ,h_{m}\rangle ),\ \epsilon ).$$ 
If $U$ is small, then for each 
$(g_{1},\ldots ,g_{m})\in U$, 
$\bigcup _{j=1}^{m}CV^{\ast }(g_{j})\subset 
D(P^{\ast }(\langle h_{1},\ldots ,h_{m}\rangle ),\ \epsilon ).$ 
Hence, if $U$ is small enough, 
then for each $(g_{1},\ldots ,g_{m})\in U$, 
$P^{\ast }(\langle g_{1},\ldots ,g_{m}\rangle )
\subset D(P^{\ast }(\langle h_{1},\ldots ,h_{m}\rangle ), \epsilon ).$ 
Hence, for each $(g_{1},\ldots ,g_{m})\in U$, 
$\langle g_{1},\ldots ,g_{m}\rangle \in {\cal G}.$ 
Therefore, ${\cal H}_{m}\cap {\cal B}_{m}$ is open. 

Thus, statement~\ref{sphypopen1} holds. 

  
 Thus, we have proved Theorem~\ref{sphypopen}.   
\qed

{\bf The author's E-mail address: sumi@math.sci.osaka-u.ac.jp}
\end{document}